\documentclass{amsart}
\usepackage[english]{babel}
\usepackage{graphicx}
\usepackage[hidelinks]{hyperref}
\usepackage{amssymb}
\usepackage{amsmath, mathtools}
\usepackage[foot]{amsaddr}
\usepackage{amsfonts}
\usepackage{bbm}
\usepackage{xcolor}
\usepackage{tikz}

\usepackage[top=1in, bottom=1.25in, left=1.25in, right=1.25in]{geometry}
\usepackage[bf, small]{caption}
\usepackage{subcaption}
\usepackage{multirow}
\usepackage{algorithm}
\usepackage{algpseudocode}
\usepackage{siunitx}

\theoremstyle{theorem} 
\newtheorem{Theorem}{Theorem}[subsection]
\newtheorem{Problem}[Theorem]{Problem}
\newtheorem*{Problem*}{Problem}
\newtheorem{Definition}[Theorem]{Definition}

\newtheorem{Assumption}[Theorem]{Assumption}

\newtheorem{Remark}[Theorem]{Remark}

\definecolor{DEblue}{rgb}{0.25, 0.0, 1.0}

\begin{document}
\title[Stabilized weighted ROMs for Advection-Dominated OCPs with Random Inputs]{Stabilized Weighted Reduced Order Methods for Parametrized Advection-Dominated Optimal Control Problems governed by Partial Differential Equations with Random Inputs}
\author{Fabio Zoccolan$^1$, Maria Strazzullo$^2$, and Gianluigi Rozza$^3$}
\address{$^1$ Institut de Mathématiques, École Polytechnique Fédérale de Lausanne, 1015 Lausanne, Switzerland, \newline email: fabio.zoccolan@epfl.ch}
\address{$^2$ DISMA, Politecnico di Torino, Corso Duca degli Abruzzi 24, Turin, Italy. \newline email: maria.strazzullo@polito.it}
\address{$^3$ mathlab, Mathematics Area, SISSA, via Bonomea 265, I-34136 Trieste, Italy. \newline email: gianluigi.rozza@sissa.it}

\begin{abstract}
In this work, we analyze Parametrized Advection-Dominated distributed Optimal Control Problems with random inputs in a Reduced Order Model (ROM) context. All the simulations are initially based on a finite element method (FEM) discretization; moreover, a \textit{space-time approach} is considered when dealing with unsteady cases. To overcome numerical instabilities that can occur in the optimality system for high values of the P\'eclet number, we consider a Streamline Upwind Petrov–Galerkin technique applied in an optimize-then-discretize approach. We combine this method with the ROM framework in order to consider two possibilities of stabilization: \textit{Offline-Only stabilization} and \textit{Offline-Online stabilization}. Moreover we consider random parameters and we use a \textit{weighted Proper Orthogonal Decomposition} algorithm in a partitioned approach to deal with the issue of uncertainty quantification. Several quadrature techniques are used to derive weighted ROMs: tensor rules, isotropic sparse grids, Monte-Carlo and quasi Monte-Carlo methods. We compare all the approaches analyzing relative errors between the FEM and ROM solutions and the computational efficiency based on the speedup-index.
\end{abstract}

\keywords{reduced order methods, time-dependent parametrized optimal control problem, stabilization, weighted proper orthogonal decomposition, random inputs, uncertainty quantification}

\maketitle

\section{Introduction}
\label{sec:intro}

Here we present a numerical study concerning stabilized Parametrized Advection-Dominated Optimal Control Problems (OCP($\boldsymbol{\mu}$)s) with random inputs in a Reduced Order Methods (ROMs) framework. As a matter of fact, engineering and scientific applications often need very fast evaluations of the
numerical solutions for many parameters that characterize the problem, for instance in \textit{real-time} scenarios. A solution to these many-query situations can be to exploit the parameter dependence of the OCP($\boldsymbol{\mu}$)s using ROMs \cite{benner2017model,hesthaven2016certified,quarteroni2011certified,quarteroni2014reduced,quarteroni2015reduced}. This process is divided in two stages. The former is the \textit{offline phase}, when many numerical solutions for different values of parameters are collected considering a first discretization of the OCP($\boldsymbol{\mu}$), such a finite element method (FEM) one, called the high-fidelity or truth approximation. Then all parameter-independent components are calculated and stored, and reduced spaces are built. The latter is the \textit{online phase}, when all parameter-dependent parts and, then, the reduced solutions are computed. More precisely, to deal with the randomness which is hidden in the parameters, we consider a modification of the Proper Orthogonal Decomposition (POD) that takes into account the probability distribution of the random inputs: the \textit{weighted POD} (wPOD) \cite{venturi2018weighted, venturi2019weighted}. We apply this procedure in a \textit{partitioned approach}, following good results shown in literature \cite{karcher2018certified,negri2015reduced,strazzullo2018model, zoccolan2024streamline}. As this algorithm aims to minimize the expectation of the square error between the truth and the ROM solutions, we can identify different types of \textit{weighted ROMs} (wROMs) \cite{carere2021weighted,chen2013weighted,chen2014weighted, chen2014comparison, chen2016multilevel,chen2017reduced,spannring2017weighted, torlo2018stabilized, venturi2018weighted, venturi2019weighted} based on the chosen quadrature rules. In this work, we will exploit Monte-Carlo and Quasi Monte-Carlo procedures, tensor rules based on Gauss-Jacobi and Clenshaw-Curtis quadrature techniques, and Smolyak isotropic sparse grids.

The optimization problem will always concern a linear-quadratic cost functional. We use FEM as the truth approximation, both for steady and unsteady problems. At a first level, FEM approximations of stochastic steady OCP($\boldsymbol{\mu}$)s have been already presented, for example, in \cite{hou2011finite} considering stochastic PDEs. In the parabolic case, we discretize time-dependency via a \textit{space-time approach} \cite{hinze2008hierarchical,stoll2010all, STOLL2013498}.
Concerning stabilization, we considered the Streamline Upwind Petrov–Galerkin (SUPG) \cite{brooks1982streamline, hughes1987recent, quarteroni2009numerical} suitably combined with the ROM setting in an \textit{optimize-then-discretize} approach. We exploit two possibilities: when stabilization only occurs in the offline phase, \textit{Offline-Only stabilization} or when SUPG is applied in both phases, \textit{Offline-Online stabilization}. 
Even if the choice of applying stabilization both in the offline and online stages seems to be expected, in many works was not a natural choice. We refer the interested reader to the introduction of \cite{Strazzullo20223148}.\\
{Furthermore, we decide to implement an ``optimize-then-discretize" framework so as to have a direct comparison with the results shown in the recent works on optimal control on similar topics \cite{strazzullo2020pod,strazzullo2021certified,strazzullo2022pod}. However, another framework is possible: the  ``discretize-then-optimize" approach. In this case, a numerical discretization is first applied to the governing PDEs with possible stabilization and, further, we build an optimal control problem based on the full order model. The main advantage of this technique is that the derived optimal system is symmetric, but, as a drawback, this method is not strongly consistent; i.e.\  the exact solution of the non-stabilized framework is not a solution of the stabilized system, too, unlike the ``optimize-then-discretize" case.}

Stabilized Advection-Dominated problems in a ROM framework without control are studied, for instance, in \cite{pacciarini2014stabilized,torlo2018stabilized}, both for steady and unsteady cases. Instead, in \cite{carere2021weighted} wROMs for generic OCP($\boldsymbol{\mu}$)s are applied to experiments concerning
environmental sciences, {whereas} SUPG Advection-Dominated distributed OCP($\boldsymbol{\mu}$)s are analyzed in a deterministic context in \cite{zoccolan2024streamline}, both for elliptic and parabolic experiments.
To the best of our knowledge, this is the first time that stabilized Advection-Dominated OCP($\boldsymbol{\mu}$)s with random inputs are analyzed in a ROM context, both for elliptic and parabolic problems.
This work is arranged as follows. In Section \ref{sec:problem}, we introduce linear-quadratic optimal control theory for PDEs and its FEM discretization. Section \ref{sec:supg} firstly concerns the basic theory of SUPG stabilization for Advection-Dominated PDEs in an \textit{optimize-then-discretize} approach \cite{collis2002analysis}, then the space-time procedure that will be used is presented. wROMs features will be illustrated in Section \ref{sec:wROMs}. Section \ref{sec:results} will concern numerical simulations. Two Advection-Dominated problems under distributed control and random inputs will be analyzed: the Graetz-Poiseuille Problem under geometrical parametrization and the Propagating Front in a Square Problem. We will compare the wPOD procedures through relative errors between the FEM and the ROM solutions and computational time considering the speedup-index. Finally, conclusions follow in Section \ref{sec:conclusions}.

\section{Problem Formulation for Random Input Optimal Control Problems}
\label{sec:problem}

\subsection{Mathematical Setting}
Let $\Omega$ be an open and bounded regular domain in $\mathbb{R}^{2}$, where $\Gamma_{N}$ and $\Gamma_{D}$ will indicate the portions of the boundary $\partial \Omega$ where Neumann and Dirichlet boundary conditions are imposed, respectively. {With the symbol $\Omega_{obs} \subseteq \Omega$ the \textit{observation domain} will be indicated, i.e.\ the subset of the domain where we seek the state variable to be similar to a \textit{desired solution profile} $y_d \in Y$, with $Y$ Hilbert space, in a sense that will be specified later.} For time-dependent problems we will also take into account the time interval $(0, T) \subset \mathbb{R}^{+}$. Let us consider a compact set $\mathcal{P} \subseteq \mathbb{R}^p$, for natural number $p$. We will call $\mathcal{P}$ and  as the \textit{parameter space} and a $p$-vector $\mathbf{\boldsymbol{\mu}} \in \mathcal{P}$ is the parameter of our Parametric OCP($\boldsymbol{\mu}$)s. As the setting is completely general, for instance, $\boldsymbol{\mu}$ can characterize our $y_d$ or geometrical and physical properties of the problem. Furthermore, we denote with $\mathcal{B}(Q,R)$ the space of linear continuous operators between Banach spaces $Q$ and $R$.

The triplet $(\mathfrak{A}, \mathcal{F}, P)$ will denote a complete probability space, composed by $\mathfrak{A}$, which is the set of outcomes $\omega \in \mathfrak{A}$,  $\mathcal{F}$, that is a $\sigma$-algebra of events, and $P: \mathcal{F} \rightarrow[0,1]$ with $P(\mathfrak{A})=1$, which is the chosen probability measure.
As dealing with random input OCP($\boldsymbol{\mu}$)s, the parameter $\boldsymbol{\mu}$ will be a
real-valued random vector. In detail, $\boldsymbol{\mu} :(\mathfrak{A}, \mathcal{F}) \rightarrow(\mathbb{R}^{p}, \mathfrak{B})$ is a measurable function, where $\mathfrak{B}$ is the Borel $\sigma$-algebra on $\mathbb{R}^{p}$.
The distribution function of $\boldsymbol{\mu}: \mathfrak{A} \rightarrow \mathcal{P} \subset \mathbb{R}^{p}$, being $\mathcal{P}$ the image of $\boldsymbol{\mu}$, is defined as $\mathbb{P}_{\boldsymbol{\mu}} : \mathcal{P} \rightarrow[0,1]$ such that 
\begin{equation}
    \forall \mu \in \mathcal{P}, \quad \mathbb{P}_{\boldsymbol{\mu}}(\mu)=P(\omega \in \mathfrak{A}: \boldsymbol{\mu}(\omega) \boldsymbol{\leq} \mu). 
\end{equation}
Let $\mathrm{d} \mathbb{P}_{\boldsymbol{\mu}}(\mu)$ denote the distribution measure of $\boldsymbol{\mu}$, i.e.,
\begin{equation}
    \forall H \subset \mathcal{P}, \quad  P( \boldsymbol{\mu} \in H)=\int_{H} \mathrm{d} \mathbb{P}_{\boldsymbol{\mu}}(\mu).
\end{equation} 
We assume that $\boldsymbol{\mu}$ admits a Lebesgue density, i.e.\ $\mathrm{d} \mathbb{P}_{\boldsymbol{\mu}}(\mu)$ is absolutely continuous with respect to the Lebesgue measure $\mathrm{d} \mu$. This practically means that there exists a probability density function $\rho_{\boldsymbol{\mu}}: \mathcal{P} \rightarrow \mathbb{R}^{+}$ such that $\rho_{\boldsymbol{\mu}}(\mu) \mathrm{d} \mu=\mathrm{d} \mathbb{P}_{\boldsymbol{\mu}}(\mu)$. It is worth to notice that the measure space $(\mathcal{P}, \mathfrak{B}(\mathcal{P}), \rho_{\boldsymbol{\mu}}(\mu) \mathrm{d} \mu)$ is isometric to $(\mathfrak{A}, \mathcal{F}, {P})$ under the random vector $\boldsymbol{\mu}$.

The aim of this work is to analyze random input OCP($\boldsymbol{\mu}$)s from the numerical point of view.

\begin{Problem}[Random Input Parametric Optimal Control Problem]
\label{OCP-general}
Consider the state equation $\mathcal{E} : Y \times U \to Q$, with  $Y,U$, and $Q$ real Banach spaces, satisfying a set of boundary and/or initial conditions, and a real functional $\mathcal{J}: Y \times U \to \mathbb{R}$. Then for $\mathbb{P}_{\boldsymbol{\mu}}$-a.e. find the pair $\big(y(\boldsymbol{\mu}), u(\boldsymbol{\mu})\big) \in X:=Y\times U$ that minimizes cost functional $\mathcal{J}(y(\boldsymbol{\mu}), u(\boldsymbol{\mu}); \boldsymbol{\mu})$ under the constraint $\mathcal{E}(y(\boldsymbol{\mu}), u(\boldsymbol{\mu}); \boldsymbol{\mu}) = 0$. 
\end{Problem}

Let $X_{ad}$ be the set of all couples $(y,u)$ solutions of $\mathcal{E}$: we will only consider the case of \textit{full admissibility}, i.e.\ when $X_{ad} = Y \times U$. Problem \ref{OCP-general} {looks for minimizers} among all state-control pairs such that:
$$\min_{(y(\boldsymbol{\mu}),u(\boldsymbol{\mu})) \in Y \times U} \mathcal{J}(y(\boldsymbol{\mu}),u(\boldsymbol{\mu}); \boldsymbol{\mu}) \ \text{s.t.} \ \mathcal{E}(y(\boldsymbol{\mu}),u(\boldsymbol{\mu}); \boldsymbol{\mu})= 0.$$

This can be achieved through the research of the critical points of the Lagrangian operator $\mathfrak{L} : Y \times U \times Q^{*} \to \mathbb{R}$ defined as:
\begin{equation}
  \label{lagrangian}
   \mathfrak{L}(y(\boldsymbol{\mu}),u(\boldsymbol{\mu}),p(\boldsymbol{\mu}); \boldsymbol{\mu}) = \mathcal{J}(y(\boldsymbol{\mu}),u(\boldsymbol{\mu}); \boldsymbol{\mu})+\langle p(\boldsymbol{\mu}),\mathcal{E}(y(\boldsymbol{\mu}),u(\boldsymbol{\mu}); \boldsymbol{\mu})\rangle_{Q^{*}Q},
\end{equation}
where {$p(\boldsymbol{\mu})$ is a Lagrange multiplier belonging to} the \textit{adjoint {space}} $Q^{*}$, the dual space of $Q$. For the sake of notation we write $y:=y(\boldsymbol{\mu})$, $u:=u(\boldsymbol{\mu})$ and $p:=p(\boldsymbol{\mu})$. In case that $\mathbb{P}_{\boldsymbol{\mu}}$ is the uniform distribution with support in $\mathcal{P}$, then Problem \ref{OCP-general} is called to be \textit{deterministic} problem. In this work, linear-quadratic problems will be involved.

\begin{Definition}[Linear-Quadratic OCP($\boldsymbol{\mu}$]\label{def:lin-quad-prob}
Let us consider the bilinear forms $m: Z \times Z \rightarrow \mathbb{R} $ and $n: U \times U \rightarrow \mathbb{R}$, which are symmetric and continuous, where $Z$ is a Banach space called the observation space. Fix $\alpha>0$, a constant called the \textit{penalization parameter} and consider a quadratic objective functional $\mathcal{J}$ of the form
\begin{equation}
\label{quadratic-functional}
   \mathcal{J}(y, u; \boldsymbol{\mu})=\frac{1}{2} m\left(\mathcal{O} y(\boldsymbol{\mu})-z_{d}(\boldsymbol{\mu}), \mathcal{O}y(\boldsymbol{\mu})-z_{d}(\boldsymbol{\mu})\right)+\frac{\alpha}{2} n(u(\boldsymbol{\mu}), u(\boldsymbol{\mu})),
\end{equation}

where $\mathcal{O} : Y \to Z$ is a linear and bounded operator called the \textit{observation map} and $z_{d}(\boldsymbol{\mu}) \in Z$ is the observed desired solution profile.
Consider an affine map $\mathcal{E}$ defined as \begin{equation}
\label{affine-state}
\mathcal{E}(y, u; \boldsymbol{\mu})=\mathcal{A}(\boldsymbol{\mu}) y+\mathcal{C}(\boldsymbol{\mu}) u-f(\boldsymbol{\mu}), \quad \forall \big(y(\boldsymbol{\mu}), u(\boldsymbol{\mu})\big) \in Y \times U,
\end{equation}
where  $\mathcal{A}(\boldsymbol{\mu}) \in \mathcal{B}(Y, Q)$, $\mathcal{C}(\boldsymbol{\mu}) \in \mathcal{B}(U, Q)$ and $f(\boldsymbol{\mu}) \in Q$.

Then an OCP($\boldsymbol \mu$)s with the above properties is said a \textit{Linear-Quadratic Optimal Control Problem}.
\end{Definition}

For Linear-Quadratic OCP($\boldsymbol{\mu}$)s well-posedness of Problem \ref{def:lin-quad-prob} yields {when $m$ is positive definite and $n$ is coercive} \cite{brezzi1974existence,brezzi2012mixed}. More precisely, the reader can refer to \cite{carere2019thesis} to a comparison between the Lagrangian approach for the full-admissibility case and the adjoint one. Via the functional derivative of  $\mathfrak{L}$, we obtain an optimality system to be solved to find the unique solution. In this case, this reads as finding $(y, u, p) \in Y \times U \times Q^{*}$ that satisfies \cite{carere2019thesis},
\begin{equation}\label{Opt-lin_system}
\begin{aligned}
\begin{cases}
D_y\mathfrak{L}(y, u, p; \boldsymbol{\mu})(\bar{y})=0 \ \Longrightarrow \ m(\mathcal{O} y, \mathcal{O} \bar{y}; \boldsymbol{\mu})+\left\langle \mathcal{A}^{*}(\boldsymbol{\mu}) p, \bar{y}\right\rangle_{Y^{*} Y} =m\left(\mathcal{O} \bar{y}, z_{d}; \boldsymbol{\mu}\right), & \forall \bar{y} \in Y, \\
D_u\mathfrak{L}(y, u, p; \boldsymbol{\mu})(\bar{u})=0 \ \Longrightarrow \ \alpha n(u, \bar{u}; \boldsymbol{\mu})+\left\langle \mathcal{C}^{*}(\boldsymbol{\mu}) p, \bar{u}\right\rangle_{U^{*} U} =0, & \forall \bar{u} \in U, \\
D_p\mathfrak{L}(y, u, p; \boldsymbol{\mu})(\bar{p})=0 \ \Longrightarrow \ \langle\bar{p}, \mathcal{A}(\boldsymbol{\mu}) y+\mathcal{C}(\boldsymbol{\mu}) u\rangle_{Q^{*} Q} =\langle\bar{p}, f(\boldsymbol{\mu})\rangle_{Q^{*} Q}, & \forall \bar{p} \in Q^{*}.
\end{cases}
\end{aligned}
\end{equation}

In system \eqref{Opt-lin_system}, the first equation is called the \textit{adjoint equation}, the second one is the \textit{gradient equation} and the last one is \textit{state equation}. 

\begin{Remark}[Notation] \label{remark-notation}
  For the sake of notation, when Hilbert spaces will be taken into account, bilinear forms $A,B$, and their adjoint counterparts will be indicated uniquely as:
  \begin{equation*}
      \begin{aligned}
           \langle \mathcal{A}(\boldsymbol{\mu}) y, p\rangle_{Q Q^{*}} := a(y,p; \boldsymbol{\mu}) \qquad & \langle \mathcal{C}(\boldsymbol{\mu}) u, p\rangle_{Q Q^{*}} := c(u,p; \boldsymbol{\mu}).
      \end{aligned}
\end{equation*}

\end{Remark}

\begin{Remark}[Parabolic Problems]
  Concerning unsteady problems, one must add more hypotheses to the mathematical setting of Linear-Quadratic OCP($\boldsymbol{\mu}$)ss to reach well-posedness. We will consider the following Hilbert spaces $\mathcal{Y}=L^{2}(0, T ; Y),$ $\mathcal{Y}^{*}= L^{2}\left(0, T ; Y^{*}\right)$, $\mathcal{Z}= L^{2}\left(0, T ; Z\right)$, $\mathcal{U}=L^{2}(0, T ; U)$
with respective norms given by
\begin{equation}
\|y\|_{\mathcal{Y}}^{2}:=\int\limits_{0}^{T}\|y\|_{Y}^{2} \mathrm{dt}, \quad  \|y\|_{\mathcal{Y}^{*}}^{2}:=\int\limits_{0}^{T}\|y\|_{Y^{*}}^{2} \mathrm{dt}, \quad \|z\|_{\mathcal{Z}}^{2}:=\int\limits_{0}^{T}\|z\|_{Z}^{2} \mathrm{dt}, \quad \text{and} \quad \|u\|_{U}^{2}:=\int\limits_{0}^{T}\|u\|_{U}^{2} \mathrm{dt}.
\end{equation}
\end{Remark}
Then we define the Hilbert space $\mathcal{Y}_{t}:=\left\{ y \in \mathcal{Y} \quad \text{s.t.} \quad \partial_t y \in \mathcal{Y}^{*} \right\}$.
For parabolic problems we will also consider the case of full-admissibility as $X_{a d} = \mathcal{Y}_{t} \times U$ \cite{BALLARIN2022307, strazzullo2020pod,strazzullo2021certified}. 

\subsection{High-Fidelity Discretization}
\label{sec:truth_discretization}

In this work, the discretization of the optimality sistem \eqref{Opt-lin_system} follows an \textit{one shot} or \textit{all-at-once approach} \cite{hinze2008hierarchical,stoll2010all, STOLL2013498}. When we will consider Advection-Dominated OCP($\boldsymbol{\mu}$)s, a stabilization technique will be also needed. Therefore, a SUPG method will be applied in a \textit{optimize-then-discretize} approach, as we will see in Section \ref{sec:supg}. This means that firstly the optimality conditions are computed by obtaining system \eqref{Opt-lin_system} and then we discretize and stabilize it.

Concerning numerical implementation, we use a FEM discretization for all three variables, where $\mathcal{T}_h$ is {a regular triangularization} on $\Omega$. Its elements $K$ are triangles and the parameter $h$ denotes the \textit{mesh size}, i.e.\ the maximum diameter of an element of the chosen grid. In addition, we define
$$
\Omega_{h}:=\operatorname{int}\left(\bigcup_{K \in \mathcal{T}_{h}}{K}\right),
$$
as a quasi-uniform mesh for $\Omega$. Considering $\mathbb{P}^{r}(K)$ as the space of polynomials of degree at most equal to $r$ defined on $K$ and defining 
$$
X^{\mathcal{N}, r}=\left\{q^{\mathcal{N}} \in C(\bar{\Omega}): q^{N}_{|_{K}} \in \mathbb{P}^{r}(K), \forall K \in \mathcal{T}_{h}\right\}
$$
we set
$Y^{\mathcal{N}}=Y \cap X^{\mathcal{N}, r}$, $U^{\mathcal{N}}=U \cap X^{\mathcal{N}, r}$ and $\big(Q^{\mathcal{N}}\big)^{*}={Q^{*}} \cap X^{\mathcal{N}, r}$. In this work, the numerical implementation will always be made by a $\mathbb{P}^{1}$-FEM approximation and the same triangulation $\mathcal{T}_{h}$ for $Y^{\mathcal{N}},U^{\mathcal{N}}$, and $\big(Q^{\mathcal{N}}\big)^{*}$. A similar discussion can be made for time-dependent problems, as we will see in Section \ref{time-dep-discr}. This first discretization procedure will be indicated as the \textit{truth} or \textit{high-fidelity} approximation.

From now on, $Y,U,Q$ will be always Hilbert spaces and we will consider the Identity restricted to our observation domain $\Omega_{obs}$ as the Observation function $\mathcal{O}$ for both steady and unsteady problems. Therefore, 
${Z}={Y}$ for steady problems and $\mathcal{Z}=\mathcal{Y}$ for unsteady ones are assumed. Our desired state will be denoted by $y_d$ and with the same symbol will also indicate its FEM discretization.

\section{SUPG stabilization for Advection-Dominated OCP($\boldsymbol{\mu}$)s}
\label{sec:supg}

In this work, we only deal with \textit{Advection-Diffusion equations}.

\begin{Definition}[Advection-Diffusion Equations]\label{advection-diffusion}
 Let us take into account the following problem:
\begin{equation}\label{advection-diffusion problem}
T (\boldsymbol{\mu})y:={-\mathrm{div}(\gamma(\boldsymbol{\mu}) \nabla y)}+\mathbf{\eta} (\boldsymbol{\mu}) \cdot \nabla y=f(\boldsymbol{\mu})  \text { in } \Omega \subset \mathbb{R}^2,
\end{equation}
where suitable boundary conditions are applied on $ \partial \Omega$. In addition, we require that:
 \begin{itemize}
     \item the diffusion coefficient $\gamma: \Omega \rightarrow \mathbb{R}$ is uniformly bounded, i.e.\ there exists $\gamma_{\max}, \gamma_{\min} >0$ such that
 \begin{equation}
     P\big( \omega \in \mathfrak{A} \ : \ \gamma_{\min} < \gamma(x; \boldsymbol{\mu}) < \gamma_{\max} \ \forall x \in \overline{\Omega} \big) = 1.
 \end{equation} 
 {Moreover, we request that $\gamma$ is a.e. differentiable in $\Omega$.}
\item the advection field $\mathbf{\eta}: \Omega \rightarrow \mathbb{R}^{2}$ belongs to $\left(L^{\infty}(\Omega)\right)^{2}$ for a.e. $\boldsymbol{\mu} \in \mathcal{P}$. More precisely, for a.e. $\boldsymbol{\mu}$ the following inequality holds:
$
0 \geq \operatorname{div} \mathbf{\eta}(x) \geq-\vartheta, \ \forall x \in \Omega,
$
with $\vartheta \in \mathbb{R}^{+}_{0}$;
\item $f: \Omega \rightarrow \mathbb{R}$ is an $L^{2}(\Omega)$-function for a.e. $\boldsymbol{\mu}$; in addition, $f$ has bounded second moments with respect to the integral along $\mathfrak{A}$ and $\Omega$.
 \end{itemize}

With these hypotheses, Problem \eqref{advection-diffusion problem} is called Advection-Diffusion problem and the operator $T (\boldsymbol{\mu})y:={-\mathrm{div}(\gamma(\boldsymbol{\mu}) \nabla y)}+\mathbf{\eta} (\boldsymbol{\mu}) \cdot \nabla y$ is said the Advection-Diffusion operator.
\end{Definition}

For more details regarding the well-posedness and theoretical results of Stochastic Advection-Diffusion OCP($\boldsymbol{\mu}$)s, we refer to \cite{chen2013stochastic,chen2014weighted}.

\begin{Definition}[P\'eclet number and Advection-Dominated problem]\label{adv-problem-def}
Let us consider  the FEM discretization related to an Advection-Diffusion problem and its regular triangulation $\mathcal{T}_{h}$. For any element $K \in \mathcal{T}_{h}$, the \textit{local P\'eclet number} is defined as \cite{quarteroni2008numerical, quarteroni2009numerical}:
\begin{equation}\label{Péclet number}
\mathbb{P}\mathrm{e}_{K}(x):=\frac{|\mathbf{\eta}(x)| h_{K}}{2 \gamma(x)} \quad \forall x \in K,
\end{equation}
where $h_{K}$ is the diameter of $K$. If
$
\mathbb{P}\mathrm{e}_{K}(x)>1 \ \forall x \in K, \ \forall K \in \mathcal{T}_{h},
$
we say to study an Advection-Dominated problem.
\end{Definition}

\subsection{Setting for Stabilized Advection-Dominated OCP($\boldsymbol{\mu}$)s - Steady case}

Numerical instabilities can appear, when dealing with Advection-Dominated OCP($\boldsymbol{\mu}$)s, i.e.\ when $|\mathbf{\eta}(\boldsymbol{\mu})| \gg \gamma$. In order to adjust this unpleasant behaviour without modifying the mesh size, we use the \textit{Streamline upwind/Petrov Galerkin} (SUPG) method \cite{brooks1982streamline,hughes1979multidimentional,hughes1987recent, quarteroni2008numerical} in a \textit{optimize-then-discretize approach} \cite{collis2002analysis}. This assures the strong consistency of the optimality system \cite{collis2002analysis}. For the sake of simplicity, we define our Advection-Dominated problem on $H_{0}^{1}(\Omega)$ and we do not indicate the parameter dependence. We denote with $T^{*}$ the adjoint operator related to $T$. This last operator can be split into its symmetric and skew-symmetric parts as $T=T_{S}+T_{S S}$ \cite{quarteroni2008numerical}, where:
\begin{equation}\label{sym-skewsym}
\begin{aligned}
    \text{symmetric part: }&  T_{S} y={-\mathrm{div}(\gamma(\boldsymbol{\mu}) \nabla y)}-\frac{1}{2}(\operatorname{div} \mathbf{\eta}) y, \quad \forall y \in H_{0}^{1}(\Omega),\\ 
    \text{skew-symmetric part: }&  T_{S S} y=\mathbf{\eta} \cdot \nabla y+\frac{1}{2}(\operatorname{div} \mathbf{\eta})y, \quad \forall y \in H_{0}^{1}(\Omega).
\end{aligned}
\end{equation}

This two parts can be immediately recovered using the formulae:
\begin{equation}\label{sks-s-parts}
T_{S} =\frac{T+T^{*}}{2}, \qquad T_{S S} =\frac{T-T^{*}}{2}.
\end{equation}
After having considered FEM spaces, the stabilization occurs in the bilinear and linear terms involved in the state and the adjoint equations. Instead, the gradient equation is left unstabilized \cite{collis2002analysis}. We recall that we use a distributed control.

We defined the {stabilized} bilinear form $a_{s}$ and $c_s$, and the {stabilized} forcing term $F_s$ as 
\begin{equation}\label{supg-form}
\begin{aligned}
    a_{s}\left(y^{\mathcal{N}}, q^{\mathcal{N}} ; \boldsymbol{\mu}\right)&:=a\left(y^{\mathcal{N}}, q^{\mathcal{N}} ; \boldsymbol{\mu}\right)+\sum_{K \in \mathcal{T}_{h}} \delta_{K}\left(T y^{\mathcal{N}}, \frac{h_{K}}{|\mathbf{\eta}|} T_{SS}q^{\mathcal{N}}\right)_{K}, \quad y^{\mathcal{N}}, q^{\mathcal{N}} \in Y^{\mathcal{N}}, \\
    c_{s}\left(u^{\mathcal{N}},q^{\mathcal{N}}; \boldsymbol{\mu}\right)&:= - \int\limits_{\Omega} u^{\mathcal{N}} q^{\mathcal{N}} - \sum_{K \in \mathcal{T}_{h}} \delta_{K}\left( u^{\mathcal{N}}, \frac{h_{K}}{|\mathbf{\eta}|} T_{SS}q^{\mathcal{N}}\right)_{K}, \quad u^{\mathcal{N}} \in U^{\mathcal{N}}, q^{\mathcal{N}} \in Y^{\mathcal{N}}, \\
 F_s\left(q^{\mathcal{N}} ; \boldsymbol{\mu}\right) &:=F\left(q^{\mathcal{N}} ; \boldsymbol{\mu}\right)+\sum_{K \in \mathcal{T}_{h}} \delta_{K}\left( f(\boldsymbol{\mu}), \frac{h_{K}}{|\mathbf{\eta}|} T_{SS}q^{\mathcal{N}}\right)_{K}, \quad \forall q^{\mathcal{N}} \in Y^{\mathcal{N}}.
 \end{aligned}
\end{equation}
where $\delta_{K}$ is a local positive dimensionless parameter related to the element $K \in \mathcal{T}_{h}$, consequently it can be different for each triangle, and $(\cdot, \cdot)_{K}$ is the inner scalar product in $L^{2}(K)$.
In \eqref{supg-form}, $a\left(y^{\mathcal{N}}, q^{\mathcal{N}} ; \boldsymbol{\mu}\right) = \big(T y^{\mathcal{N}}, q^{\mathcal{N}}\big)_{L^2(\Omega)}$ and $F\left(q^{\mathcal{N}} ; \boldsymbol{\mu}\right) =\big(f, q^{\mathcal{N}}\big)_{L^2(\Omega)}$, where $f$ collects all forcing and lifting terms of the problem.

For the remaining conditions of the optimality system, we will always consider $m$ and $n$ form as the $L^2(\Omega_{obs})$ and the $L^2(\Omega)$ products for steady problems.
The adjoint equation is an Advection-Dominated equation, too, where the advective term has opposite sign with respect to the state one: indeed, $T^{*} = T_{S}-T_{SS}$ from \eqref{sks-s-parts}. We use the next SUPG forms for $z^{\mathcal{N}} \in Y^{\mathcal{N}}$:
\begin{equation}
\begin{aligned}
a_{s}^{*}\left(z^{\mathcal{N}}, p^{\mathcal{N}}; \boldsymbol{\mu}\right) &:=a^{*}\left(z^{\mathcal{N}}, p^{\mathcal{N}}; \boldsymbol{\mu}\right)+\sum_{K \in \mathcal{T}_{h}} \delta^{a}_{K}\left( (T_{S}-T_{SS}) p^{\mathcal{N}}, \frac{h_{K}}{|\mathbf{\eta}|}\left(-T_{SS}\right) z^{\mathcal{N}}\right)_{K}, \\
\big( y^{\mathcal{N}}-y_d, z^{\mathcal{N}}; \boldsymbol{\mu}\big)_{s} &:= \int\limits_{\Omega_{obs}} (y^{\mathcal{N}}-
    y_d)z^{\mathcal{N}} \ \mathrm{dx} +\sum_{K \in {\mathcal{T}_{h}}_{\vert_{\Omega_{obs}}}} \delta^{a}_{K}\left( y^{\mathcal{N}}-y_d, \frac{h_{K}}{|\mathbf{\eta}|}\left(-T_{SS}\right) z^{\mathcal{N}}\right)_{K},
\end{aligned}
\end{equation}
where $a^{*}$ is the adjoint form of $a$, $\delta^{a}_{K}$ is the positive stabilization parameter of the stabilized adjoint equation. In our numerical experiments, we will always consider $\delta_{K}= \delta^{a}_{K}$. Finally, the SUPG optimality system for a steady OCP($\boldsymbol{\mu}$) reads as:
\begin{equation}\label{supg-system}
    \begin{aligned}
    &\textit{discretized adjoint equation:} \ &a_{s}^{*}\left(z^{\mathcal{N}}, p^{\mathcal{N}}; \boldsymbol{\mu}\right)+\big( y^{\mathcal{N}}-y_d, z^{\mathcal{N}}; \boldsymbol{\mu}\big)_{s}=0, & \  \forall z^{\mathcal{N}} \in Y^{\mathcal{N}},\\
   &\textit{discretized gradient equation:} \ &c^{*}\big(v^{\mathcal{N}}, p^{\mathcal{N}}; \boldsymbol{\mu}\big)+\alpha n\big(u^{\mathcal{N}}, v^{\mathcal{N}}; \boldsymbol{\mu}\big)=0, & \  \forall v^{\mathcal{N}} \in U^{\mathcal{N}}, \\
         &\textit{discretized state equation:} \ &a_{s}\left(y^{\mathcal{N}}, q^{\mathcal{N}}; \boldsymbol{\mu}\right)+c_{s}\left(u^{\mathcal{N}},q^{\mathcal{N}}; \boldsymbol{\mu}\right)=F_{s}(q^{\mathcal{N}}; \boldsymbol{\mu}), & \  \forall q^{\mathcal{N}} \in \left(Q^{\mathcal{N}}\right)^{*}.
    \end{aligned}
\end{equation}

We denote with  $K_{s}$ and $K_{s}^{T}$ the stiffness matrices related to the stabilized forms $a_{s}$ and $a_{s}^{*}$, respectively, $M$ is the not-stabilized mass matrix related to $n$, instead,  $M_{s}$ is the stabilized mass matrix related to $m$ after stabilization, $C_s$ is the matrix linked to stable form $c_s$, the block $C^{T}$ refers to $c$, and $f_{s}$ is the vector that contains the coefficients of the stabilized force term as components. Moreover, we consider with the symbol $\mathbf{y}$, $\mathbf{u}$ and $\mathbf{p}$ as
the {vectors of coefficients} of $y^{\mathcal{N}}$, $u^{\mathcal{N}}$ and $p^{\mathcal{N}}$, expressed in terms of the nodal basis of $Y^{\mathcal{N}}$, $U^{\mathcal{N}}$, $(Q^{\mathcal{N}})^{*}$, respectively.
Finally, the discretized block system related to \eqref{supg-system} is:
\begin{equation} \label{stab-block-system}
\left(\begin{array}{ccc}
M_{s} & 0 & K_{s}^{T} \\
0 & \alpha M &C^{T} \\
K_{s} &C_{s} & 0
\end{array}\right)\left(\begin{array}{l}
\mathbf{y} \\
\mathbf{u} \\
\mathbf{p}
\end{array}\right)=\left(\begin{array}{c}
M_{s} \mathbf{y}_{d} \\
0 \\
\mathbf{f}_{s}
\end{array}\right).
\end{equation}

\subsection{Setting for Stabilized Advection-Dominated OCP($\boldsymbol{\mu}$)s - Unsteady case} \label{time-dep-discr}
We show the SUPG approach for time-dependent OCP($\boldsymbol{\mu}$)s proposed in \cite{zoccolan2024streamline}. A classical implicit Euler discretization is applied to all forms including time-derivatives \cite{akman2017streamline,hinze2008hierarchical,stoll2010all,strazzullo2020pod,strazzullo2021certified,strazzullo2022pod}. We divide the time interval $(0,T)$ in $N_t$ sub-intervals of equal length $\Delta t := t_j - t_{j-1}$, $j \in \{1, \dots , N_t \}$.
Starting from this framework, a discretization along time is done, where each discrete instant of time is considered as a steady-state Advection-Dominated equation in a \textit{space-time approach} \cite{hinze2008hierarchical,stoll2010all, STOLL2013498, strazzullo2020pod, strazzullo2021certified,strazzullo2022pod}. In addition, SUPG stabilization occurs for time-dependent forms, too. The general scheme is described as follows.

Let us firstly define the discrete vectors $\boldsymbol{y}=\left[y_{1}^{T}, \ldots, y_{N_{t}}^{T}\right]^{T}$, $\boldsymbol{u}=\left[u_{1}^{T}, \ldots, u_{N_{t}}^{T}\right]^{T}$ and $\boldsymbol{p}=\left[p_{1}^{T}, \ldots, p_{N_{t}}^{T}\right]^{T}$, where $y_{i} \in Y^{\mathcal{N}}, u_{i} \in U^{\mathcal{N}}$ and $p_{i} \in (Q^{\mathcal{N}})^{*}$ for $1 \leq j \leq N_{t}$. Also here, $y_{j}, u_{j}$ and $p_{j}$ indicate the column vectors containing the coefficients of the FEM discretization for state, control and adjoint, respectively (unlike the steady case, there are not denoted in bold style). This implies $\mathcal{N}_{tot}=3 \times N_{t} \times \mathcal{N}$ as the global dimension of the block system. We express all other terms in based of the respective nodal basis.
The vector representing the initial condition for the state variable is $\boldsymbol{y}_{0}=\left[y_{0}^{T}, 0^{T}, \ldots, 0^{T}\right]^{T}$, where $0$ is the zero vector in $\mathbb{R}^{\mathcal{N}}$, $\boldsymbol{y}_{d}=$ $\left[y_{d_{1}}^{T}, \ldots, y_{d_{N_{t}}}^{T}\right]^{T}$ is the vector including discrete time components of the discretized desired solution profile; instead, $\boldsymbol{f}_s=\left[f_{s_{1}}^{T}, \ldots, f_{s_{N_{t}}}^{T}\right]^{T}$ corresponds to the stabilized forcing term. We recall that $Y,U,Q$ are Hilbert Spaces and, for the sake of simplicity, we assume $Y^{\mathcal{N}} \equiv (Q^{\mathcal{N}})^{*}$.
So now we can see locally the time block discretization. 

\begin{itemize}
    \item \textit{Adjoint equation}: this equation is discretized backward in time using the forward Euler method, which is equal to a backward Euler with respect to time $T-t$, for $t \in(0, T)$ \cite{eriksson1987error}. Firstly, we add a stabilized term to the form related to $\partial_t p$ and $a^{*}$ defined as: 
$$
s^{*}\left(z^{\mathcal{N}},p^{\mathcal{N}}(t); \boldsymbol{\mu}\right)=\sum_{K \in \mathcal{T}_{h}} \delta_{K}\left(-\partial_{t} p^{\mathcal{N}}(t)+ T^{*} p^{\mathcal{N}}(t), -\frac{h_{K}}{|\mathbf{\eta}|} T_{S S} z^{\mathcal{N}}\right)_{K},
$$
where we define the form
\begin{equation}
    m^{*}_{s}\left(p^{\mathcal{N}}, z^{\mathcal{N}}; \boldsymbol{\mu}\right)=\left(p^{\mathcal{N}}, z^{\mathcal{N}}\right)_{L^{2}(\Omega)}-\sum_{K \in \mathcal{T}_{h}} \delta_{K}\left(p^{\mathcal{N}}, \frac{h_{K}}{|\mathbf{\eta}|} T_{S S} z^{\mathcal{N}}\right)_{K}.
\end{equation}
Then, the time discretization is: for each $j \in \{N_t-1, N_t-2,..., 1\}$, find $p_{j}^{\mathcal{N}} \in Y^{\mathcal{N}}$ { s.t. }
\begin{equation}
\begin{aligned}
\frac{1}{\Delta t} m^{*}_{s}\left(p_{j}^{\mathcal{N}}(\boldsymbol{\mu})-p_{j+1}^{\mathcal{N}}( \boldsymbol{\mu}), z^{\mathcal{N}} ; \boldsymbol{\mu} \right)+a^{*}_{s}\left(z^{\mathcal{N}}, p_{j}^{\mathcal{N}}( \boldsymbol{\mu}); \boldsymbol{\mu}\right)= - \big( y^{\mathcal{N}}_{j}-y_{d_j}, z^{\mathcal{N}}; \boldsymbol{\mu}\big)_{s} \quad \forall z^{\mathcal{N}} \in Y^{\mathcal{N}},
\end{aligned}
\end{equation}
Considering $M_{s}^{T}$ as the matrix inherent to $m_{s}^{*}$, the block subsystem reads
$$
M_{s}^{T} p_{j}=M_{s}^{T} p_{j+1}+\Delta t\left(-M_{s}^{T} y_{j}-K_s^{T} p_{j}+M_{s}^{T} y_{d_{j}}\right) \quad \text { for } j \in\left\{N_{t}-1, N_{t}-2, \ldots, 1\right\}.
$$
Finally, we derive the following block system:
\begin{equation*}
\underbrace{\left[\begin{array}{cccc}
M_{s}^{T}+\Delta t K_s^{T} & -M_{s}^{T} & & \\
& \ddots & \ddots & \\
& & M_{s}^{T}+\Delta t K_s^{T} & -M_{s}^{T} \\
& & & M_{s}^{T}+\Delta t K_s^{T} 
\end{array}\right]}_{\mathcal{A}^{T}_s}
\boldsymbol{p}
+\left[\begin{array}{c}
\Delta t M_{s}^{T} y_{1} \\
\vdots \\
\vdots \\
\Delta t M_{s}^{T} y_{N_{t}}
\end{array}\right]  \\
{=\left[\begin{array}{c}
\Delta t M_{s}^{T} y_{d_{1}} \\
\vdots \\
\vdots \\
\Delta t M_{s}^{T} y_{d_{N_{t}}}
\end{array}\right]}.
\end{equation*}
Setting the diagonal block matrix $\mathcal{M}^{T}_s\in \mathbb{R}^{\mathcal{N} \cdot N_{t}} \times \mathbb{R}^{\mathcal{N} \cdot N_{t}}$ with diagonal entries $[ M_{s}^{T}, \ldots,  M_{s}^{T}]$, the adjoint system to be solved is:
$
\Delta t \mathcal{M}^{T}_s \boldsymbol{y}+\mathcal{A}^{T}_s \boldsymbol{p}=\Delta t \mathcal{M}^{T}_s \boldsymbol{y}_{d}.
$
\item \textit{Gradient equation}. We seek the solution of $
    \alpha \Delta t M u_{j} {+} \Delta t C^{T} p_{j}=0, \ \forall j \in\left\{1,2, \ldots, N_{t}\right\},$
which is equal to the following block system:
\begin{equation}
\alpha \Delta t \underbrace{\left[\begin{array}{cccc}
 M & &  & \\
 & M  &  & \\
 & \ddots & \ddots & \\
 & & & M
\end{array}\right]}_{\mathcal{M}}
\left[\begin{array}{c}
u_{1} \\
u_{2} \\
\vdots \\
u_{N_{t}}
\end{array}\right] 
{+}\Delta t \underbrace{\left[\begin{array}{cccc}
C^{T} & 0 & \cdots & \\ 
& C^{T} &  & \\ 
& & \ddots & \\ 
& & & C^{T}\end{array}\right]}_{\mathcal{C^{T}}}
\left[\begin{array}{c}
p_{1} \\ 
p_{2} \\ 
\vdots \\ 
p_{N_{t}}\end{array}\right]
=\left[\begin{array}{c}
0\\ 
0\\ 
\vdots \\ 
0\end{array}\right] .
\end{equation}

More compactly, we solve $
    \alpha \Delta t \mathcal{M} \boldsymbol{u}{+}\Delta t \mathcal{C}^{T} \boldsymbol{p}=0.$

\item \textit{State equation}. A backward Euler method is used for a discretization forward in time. The stabilized term related to $\partial_t y$ and the bilinear form $a$ is \cite{john2011error, pacciarini2014stabilized, rozza2016stabilized}: 
$$
s\left(y^{\mathcal{N}}(t), q^{\mathcal{N}}; \boldsymbol{\mu}\right)=\sum_{K \in \mathcal{T}_{h}} \delta_{K}\left(\partial_{t} y^{\mathcal{N}}(t)+T y^{\mathcal{N}}(t), \frac{h_{K}}{|\mathbf{\eta}|} T_{S S} q^{\mathcal{N}}\right)_{K},
$$
where $y^{\mathcal{N}}(t) \in Y^{\mathcal{N}}$ for each $t \in (0,T)$ and $q^{\mathcal{N}} \in Y^{\mathcal{N}}$. Defining the stabilized term $m_s$ as
\begin{equation}
    m_{s}\left(y^{\mathcal{N}}, q^{\mathcal{N}}; \boldsymbol{\mu}\right)=\left(y^{\mathcal{N}}, q^{\mathcal{N}}\right)_{L^{2}(\Omega)}+\sum_{K \in \mathcal{T}_{h}} \delta_{K}\left(y^{\mathcal{N}}, \frac{h_{K}}{|\mathbf{\eta}|} T_{S S} q^{\mathcal{N}}\right)_{K},
\end{equation}
then the backward Euler approach reads as: {for each} $j \in \{1, 2, \cdots, N_t\}$, find $y_{j}^{\mathcal{N}} \in Y^{\mathcal{N}}$  s.t.
\begin{equation}
\begin{aligned}
\frac{1}{\Delta t} m_{s}\left(y_{j}^{\mathcal{N}}(\boldsymbol{\mu})-y_{j-1}^{\mathcal{N}}( \boldsymbol{\mu}), q^{\mathcal{N}} ; \boldsymbol{\mu} \right)+a_{s}\left(y_{j}^{\mathcal{N}}( \boldsymbol{\mu}), q^{\mathcal{N}} ; \boldsymbol{\mu}\right)+c_{s}\left(u^{\mathcal{N}}_j,q^{\mathcal{N}}; \boldsymbol{\mu}\right) = F_s\left(q^{\mathcal{N}}; \boldsymbol{\mu}\right),
\end{aligned}
\end{equation}
given the initial condition $y^{\mathcal{N}}_{0}$ which satisfies $
\left(y^{\mathcal{N}}_{0}, q^{\mathcal{N}}\right)_{L^{2}(\Omega)}=\left(y_{0}, q^{\mathcal{N}}\right)_{L^{2}(\Omega)}, \ \forall q^{\mathcal{N}} \in Y^{\mathcal{N}}.$
The matrix state equation to be solved becomes
\begin{equation}
   M_s y_{j}+\Delta t K_s y_{j} + \Delta t C_s u_{j}=M_s y_{j-1}+f_{s_{j}} \Delta t \quad \text { for } j \in\left\{1,2, \ldots, N_{t}\right\},
\end{equation}
where the stabilized mass matrix of $m_s$ is $M_s$. Thus, we have
\begin{equation*}
\underbrace{\left[\begin{array}{cccc}
M_s+\Delta t K_s & 0 & & \\
-M_s & M_s+\Delta t K_s & 0 & \\
 &  \ddots & \ddots & 0 \\
 & 0 & -M_s & M_s+\Delta t K_s
\end{array}\right]}_{\mathcal{A}_s} 
\boldsymbol{y}
{+}\Delta t \underbrace{\left[\begin{array}{ccc}
C_s &0 & \\ 
0 &\ddots& 0\\ 
& 0& C_s\end{array}\right]}_{\mathcal{C}_s}
\boldsymbol{u}
=
\mathcal{M}_s \boldsymbol{y}_{0}+ \Delta t \boldsymbol{f}_{s},
\end{equation*}
where $\mathcal{M}_s \in \mathbb{R}^{\mathcal{N} \cdot N_{t}} \times \mathbb{R}^{\mathcal{N} \cdot N_{t}}$ is a block diagonal matrix which diagonal entries are $[M_s, \ldots, M_s]$. In a more compact notation, we have
$
\mathcal{A}_{s} \boldsymbol{y}{+}\Delta t \mathcal{C}_s \boldsymbol{u}=\mathcal{M}_s \boldsymbol{y}_{0}+\Delta t \boldsymbol{f}_s.
$

\end{itemize}
 
The final system considered and solved through an one-shot approach is the following:
\begin{equation}\label{stab-par-block}
\left[\begin{array}{ccc}
\Delta t \mathcal{M}^{T}_s & 0 & \mathcal{A}^{T}_s \\
0 & \alpha \Delta t \mathcal{M} &\Delta t \mathcal{C}^{T} \\
\mathcal{A}_s &\Delta t \mathcal{C}_s & 0
\end{array}\right]\left[\begin{array}{l}
\boldsymbol{y} \\
\boldsymbol{u} \\
\boldsymbol{p}
\end{array}\right]=\left[\begin{array}{c}
\Delta t \mathcal{M}^{T}_s \boldsymbol{y}_{d} \\
0 \\
\mathcal{M}_s \boldsymbol{y}_{0}+\Delta t \boldsymbol{f}_s
\end{array}\right].
\end{equation}

\section{Weighted ROMs for random inputs advection-dominated OCP($\boldsymbol \mu$)s}
\label{sec:wROMs}

Numerical simulations for OCP($\boldsymbol{\mu}$)s can be very expensive in relation to computational time and storage. To overcome this problem, in this work we will consider ROMs \cite{benner2017model,hesthaven2016certified,quarteroni2011certified,quarteroni2014reduced,quarteroni2015reduced}. We will study the case when the parameter $\boldsymbol{\mu}$ can be affected by randomness, i.e.\ it can follow a particular probability distribution. That is the case of \emph{random inputs OCP($\boldsymbol{\mu}$)s}. 
In this scenario, a suitable modification of the ROMs, the wROM\cite{carere2021weighted,chen2013weighted,chen2014weighted,chen2014comparison,chen2016multilevel,chen2017reduced,spannring2017weighted, torlo2018stabilized, venturi2018weighted, venturi2019weighted}, takes into account the uncertainty quantification (UQ) of the problems and shows efficient results concerning errors and computational time. For the sake of notation, from now on we denote a generic probability distribution with the symbol $\rho$. ROM procedure is divided in two stages: an \textit{offline phase} and an \textit{online phase}. 

To exploit the potential of the ROMs setting, we assume an \textit{affine decomposition} of the forms in \eqref{supg-system} \cite{hesthaven2016certified}. Therefore, Assumption \ref{linear-expansion} is required here.

\begin{Assumption}\label{linear-expansion}
We request that all the forms in \eqref{supg-system} are affine in $\boldsymbol{\mu}=(\mu_1, \ldots, \mu_p) \in \mathcal{P}$. More precisely, we request that \cite{chen2013weighted, chen2014weighted}:
\begin{enumerate}
    \item the random diffusivity $\gamma: \Omega \times \mathcal{P} \to \mathbb{R}$ is of the form
\begin{equation}
    \gamma(\boldsymbol{\mu},x)=\gamma_{0}(x)+\sum_{k=1}^{p} \theta^{\gamma}_k(\mu_{k}) \gamma_{k}(x),
\end{equation}
with $\gamma_{k} \in L^{\infty}(\Omega)$, for $k=0, \ldots, p$ and $\theta^{\gamma}_k$ depending only on $\mu_{k}$;
   \item the random advection field $\mathbf{\eta}: \Omega \times \mathcal{P} \to \mathbb{R}^{2}$ 
   is of the form
\begin{equation}
    \mathbf{\eta}(\boldsymbol{\mu}, x)=\mathbf{\eta}_{0}(x)+\sum_{k=1}^{p} \theta^{\mathbf{\eta}}_k(\mu_{k}) \mathbf{\eta}_{k}(x),
\end{equation}
with $\mathbf{\eta}_{k} \in (L^{\infty}(\Omega))^{2}$, for $k=0, \ldots, p$ and $\theta^{\mathbf{\eta}}_k$ depending only on $\mu_{k}$;
    \item  the random forcing term $f: \Omega \times \mathcal{P} \to \mathbb{R}$ is of the form
\begin{equation}    \label{force-ass}
f(\boldsymbol{\mu}, x)=f_{0}(x)+\sum_{k=1}^{p} \theta^{f}_k(\mu_{k}) f_{k}(x),
\end{equation}
with $f_{k} \in L^{2}(\Omega)$, for $k=0, \ldots, p$ and $\theta^{f}_k$ depending only on $\mu_{k}$.
\end{enumerate}
\end{Assumption}
For example, Assumption \ref{linear-expansion} can be satisfied by truncating a
Karhunen–Loève expansion \cite{schwab2006karhunen}.

\subsection{Offline phase} \label{offline phase}

The \textit{offline} phase is the most expensive stage of the wROMs, which usually depends on $\mathcal{N}$. However, this should be done only once. The aim of this procedure is to build \textit{reduced spaces} $Y^{N}$, $U^{N}$ and $(Q^{N})^{*}$ that are good approximations of the high-fidelity ones and to compute all block matrix components that are $\boldsymbol{\mu}$-independent. Then, everything is memorized in order to be ready to be used in the online phase. The construction of the reduced basis is achieved through a modified version of the POD algorithm: the wPOD \cite{carere2021weighted,venturi2019weighted,venturi2018weighted}, described in Section \ref{wpod}. Here, we firstly compute high-fidelity evaluation of optimal solutions $\left(y^{\mathcal{N}}(\boldsymbol{\mu}), u^{\mathcal{N}}(\boldsymbol{\mu}), p^{\mathcal{N}}(\boldsymbol{\mu})\right)$ for different parameters $\boldsymbol{\mu}$, the so-called snapshots, to build the bases. Because of good performance presented in literature \cite{karcher2018certified,negri2015reduced,strazzullo2018model}, this process will go through a \textit{partitioned approach}, i.e.\ the wPOD is executed separately for all three variables. After this step, the three reduced spaces for state, control, and adjoint are constructed as, respectively, 
\begin{equation}\label{basis-red}
\begin{aligned}
Y^{N} =\operatorname{span}\left\{\xi_{n}^{y}, \ n=1,  \ldots, N\right\}&,  \qquad
U^{N} =\operatorname{span}\left\{\xi_{n}^{u}, \ n=1, \ldots, N\right\}, \\
(Q^{N})^{*} =\operatorname{span}&\left\{\xi_{n}^{p},\  n=1, \ldots, N\right\}.
\end{aligned}
\end{equation}
In order to ensure well-posedness for the reduced space approximation, we need to implement an enriched space for state and adjoint variables. This means to impose $G^N \equiv Y^N\equiv (Q^{N})^{*}$, where $G^{N}=\operatorname{span}\left\{\sigma_{n}, n=1, \ldots, 2 N\right\}$ and
$\left\{\sigma_{n}\right\}_{n=1}^{2 N}=\left\{\xi_{n}^{y}\right\}_{n=1}^{N} \cup\left\{\xi_{n}^{p}\right\}_{n=1}^{N}$
\cite{doi:10.1137/090760453,doi:10.1137/110854084,karcher2018certified,kunisch2008proper, negri2013reduced,negri2015reduced}. This whole discussion holds true for parabolic problems in a space-time context, too. As a matter of fact, when dealing with time-dependent OCP($\boldsymbol{\mu}$)s in a space-time approach, the time instances are not separated in the wPOD algorithm. Therefore, each snapshot carries {all the time instances}.

\subsubsection{Weighted Proper Orthogonal Decomposition}\label{wpod}
The peculiarity of wPOD is to take into account the probability distribution that characterizes $\boldsymbol{\mu}$ to create reduced spaces with less number of basis with respect to the deterministic case without losing in accuracy \cite{carere2021weighted,venturi2019weighted,venturi2018weighted}. We will notice that there will be different ways to consider randomness in the wPOD: the general idea is to suitably attribute a larger weight to those samples that are more significant according to the distribution of $\boldsymbol{\mu}$. From now we will refer to the POD algorithm based on the Monte-Carlo procedure in a deterministic context, i.e.\ when the distribution $\rho$ is the uniform one, as \emph{Standard POD} to distinguish it from the wPOD.
As we will consider a partitioned approach, we show the procedure for the state space: adjoint and control variables will follow the same process. 

To consider stochasticity, wPOD needs to find the $N$-dimensional space $Y^{N}$, with $N \ll \mathcal{N}$, such that it minimizes the following estimate:
\begin{equation} \label{s-error-estimator}
E = \sqrt{\int_{\mathcal{P}} \inf_{\zeta^{y} \in Y^{N} }\left\|y^{\mathcal{N}}(\boldsymbol{\mu})-\zeta^{y}\right\|_{\mathcal{Y}}^{2} \rho(\boldsymbol{\mu})\mathrm{d}\boldsymbol{\mu} }.
\end{equation}
Let us consider a set of $N_{\text{train}}$ ordered parameters $\boldsymbol{\mu}_{1}, \ldots, \boldsymbol{\mu}_{N_{\text{train}}} \in \mathcal{P}_{N_{\text{train}}}$, where $\mathcal{P}_{N_{\text{train}}} \in \mathcal{P}$ is a discretization of $\mathcal{P}$ called the \emph{training set} and its cardinality is $|\mathcal{P}_{N_{\text{train}}} | = N_{\text{train}}$. One can choose $N_{\text{train}}$ so that $\mathcal{P}_{N_{\text{train}}}$ is a good approximation of $\mathcal{P}$. We can relate $\boldsymbol{\mu}_{1}, \ldots, \boldsymbol{\mu}_{N_{\text{train}}}$ to the ordered snapshots $y^{\mathcal{N}}\left(\boldsymbol{\mu}_{1}\right), \ldots, y^{\mathcal{N}}\left(\boldsymbol{\mu}_{N_{\text{train }}}\right)$. Considering $w: \mathcal{P} \rightarrow \mathbb{R}^{+}$ a weight function, a discretization of problem \eqref{s-error-estimator} is meant to find the $N$-dimensional space $Y^{N}$ which minimize the quantity
\begin{equation} \label{w-quadr}
    \frac{1}{N_{\text{train}}} \sum_{k=1}^{N_{\text{train}}} w\left(\boldsymbol{\mu}_{k}\right)\left\|y^{\mathcal{N}}\left(\boldsymbol{\mu}_{k}\right)-y^{N}\left(\boldsymbol{\mu}_{k}\right)\right\|_{Y}^{2}.
\end{equation}

One could think that the natural choice can be $w(\boldsymbol{\mu})=\rho(\boldsymbol{\mu})$ and in a UQ context this means to just discretize the expectation of the square error
\begin{equation}\label{expectation}
\mathbb{E}\left[\left\|y^{\mathcal{N}}-y^{N}\right\|_{Y}^{2}\right]:=\int_{\mathcal{P}}\left\|y^{\mathcal{N}}(\boldsymbol{\mu})-y^{N}(\boldsymbol{\mu})\right\|^{2} \rho(\boldsymbol{\mu}) d \boldsymbol{\mu},
\end{equation}
which is the argument of the square root in \eqref{s-error-estimator}. However, this is not the unique choice in this scenario: therefore it will be interesting to understand which method is better to approximate \eqref{expectation}.
Here we illustrate different techniques that we use in the numerical tests in Section \ref{sec:results} to approximate \eqref{expectation}.
Considering the training set $\mathcal{P}_{N_{\text{train}}}=\left\{\boldsymbol{\mu}_{1}, \ldots, \boldsymbol{\mu}_{ N_{\text{train}}}\right\}$, which can be composed by the nodes of the chosen quadrature formula that approximates \eqref{expectation}, we indicate with $\omega=\left(\omega_{1}, \ldots, \omega_{N_{\text{train}}}\right)$ the standard weights of a chosen quadrature rule, with $\rho_{1}, \ldots, \rho_{N_{\text{train}}}$ the values of the density $\rho$ in the nodes in $\mathcal{P}_{N_{\text{train}}}$, and with $w=\left(w_{1}, \ldots, w_{N_{\text{train}}}\right)$ the definitive weights used in wPOD algorithm. For a node $\mu_k$, we have the correspondent quantities $\omega_k,\rho_k,$ and $w_k$. As a final result of this first step, the wPOD furnished the following sum to minimize
\begin{equation}\label{w-montecarlo}
\frac{1}{N_{\text{train}}} \sum_{k=1}^{N_{\text{train}}}w_{k}\left\|y^{\mathcal{N}}\left(\boldsymbol{\mu}_k\right)-y^{N}\left(\boldsymbol{\mu}_k\right)\right\|_{Y}^{2},
\end{equation}
which is achieved here through the following algorithms:
\begin{itemize}
 \item \textit{Weighted Monte-Carlo method}, 
where $\boldsymbol{\mu}_{1}, \ldots, \boldsymbol{\mu}_{N_{\text{train}}}$ are $N_{\text{train}}$ parameters extracted from the random variable $\mu$ according to its distribution $\rho$ and $\rho_{i}$ are the values of the density $\rho$ in these points. For this approximation, we have  $\mathcal{P}_{N_{\text{train}}}=\left\{\boldsymbol{\mu}_{1}, \ldots, \boldsymbol{\mu}_{N_{\text{train}}}\right\}$ and $w_{k}=\rho(\boldsymbol{\mu}_k)$, for all $k=1, \ldots, N_{\text{train}}$;
\item \textit{Pseudo-Random method based on a Halton Sequence}, 
where $\boldsymbol{\mu}_{1}, \ldots, \boldsymbol{\mu}_{N_{\text{train}}}$ are the nodes extracted by a sampling completely based on the Halton sequence \cite{sullivan2015introduction} and $\rho_{k}=\rho(\boldsymbol{\mu}_k)$. Also in this case, $\mathcal{P}_{N_{\text{train}}}=\left\{\boldsymbol{\mu}_{1}, \ldots, \boldsymbol{\mu}_{N_{\text{train}}}\right\}$ and $w_{k}=\rho_{k}$, for all $k=1, \ldots, N_{\text{train}}$;
\item \textit{Tensor product Gauss-Jacobi rule},
where $\boldsymbol{\mu}_{1}, \ldots, \boldsymbol{\mu}_{N_{\text{train}}}$ are the nodes of the tensor product Gauss-Jacobi quadrature rule and $\omega_{1}, \ldots, \omega_{N_{\text{train}}}$ are the correspondent quadrature weights. We can use this formula when the distribution is a Beta$\left(\alpha_{k}, \beta_{k}\right)$, as suitable Jacobi polynomials are orthogonal to this distribution \cite{ralston2001first}. As a matter of fact, simulations in Section \ref{sec:results} will consider different Beta distributions for all components of $\boldsymbol{\mu}$. Therefore, we implement a Gauss-Jacobi formula using $\left(\alpha_{k}, \beta_{k}\right)$ as its parameters in each dimension \cite{ralston2001first}, accordingly to the distribution of $\boldsymbol{\mu}$. For this approximation, we have $\mathcal{P}_{N_{\text{train}}}=\left\{\boldsymbol{\mu}_{1}, \ldots, \boldsymbol{\mu}_{N_{\text{train}}}\right\}$ and $w_{k}=\omega_{k}$, for all $k=1, \ldots, N_{\text{train}};$
\item \textit{Tensor product Clenshaw-Curtis rule},
where $\boldsymbol{\mu}_{1}, \ldots, \boldsymbol{\mu}_{N_{\text{train}}}$ are the nodes of the tensor product Clenshaw-Curtis quadrature rule and $\omega_{1}, \ldots, \omega_{N_{\text{train}}}$ are the correspondent quadrature weights \cite{sullivan2015introduction}. In this case we obtain $\mathcal{P}_{N_{\text{train}}}=\left\{\boldsymbol{\mu}_{1}, \ldots, \boldsymbol{\mu}_{N_{\text{train}}}\right\}$  and $w_{k}=\rho_{k} \omega_{k}$, for all $k=1, \ldots, N_{\text{train}}$.
\end{itemize}

In numerical tests of Section \ref{sec:results},  we will respectively call as Weighted Monte-Carlo, Pseudo-Random, Gauss-Jacobi, and Clenshaw-Curtis wPOD algorithms the rules just specified. As it is known, tensor rules can be efficient, but their structure implies huge computational costs for elevate cardinality of the training set $\mathcal{P}_{\text{train}}$ or high-dimensional parameter space $\mathcal{P}$. For this purpose, when we will use Clenshaw-Curtis or Gauss-Jacobi methods, we will consider  \textit{sparse grid} techniques based on a Smolyak algorithm, too \cite{smolyak1963quadrature,xiu2005high}: we will implement \textit{isotropic} ones \cite{nobile2008sparse}.

We remark that in convection-dominated problems the choice of the distribution and of the quadrature rule can represent an issue. It is important to have enough \emph{stabilized information} and a peculiar way to pick the parameters may lead to unstable snapshots. Indeed, instabilities in the snapshots may affect the reliability of the basis functions causing a poor online approximation of the phenomenon. This, however, was not the case in our numerical results, where the numerical stabilization is suited to tackle the whole parametric space. Therefore, the performance of the different weighted POD algorithms are exclusively based on how they sample from the parameter space and how the related snapshots are properly weighted.

Once chosen the rule \eqref{w-quadr} to approximate \eqref{expectation}, the procedure to minimize \eqref{w-quadr} is described as follows.
 Let us define the deterministic correlation matrix of the snapshots of the state variable ${D}^{y} \in \mathbb{R}^{N_{\text{train}} \times N_{\text{train}}}$ in the following way:
\begin{equation}\label{correlation-matrix}
{D}_{kl}^{y}:=\frac{1}{N_{\text{train}}}\left(y^{\mathcal{N}}\left(\boldsymbol{\mu}_{k}\right), y^{\mathcal{N}}\left(\boldsymbol{\mu}_{l}\right)\right)_{Y}, \quad 1 \leq k,l \leq N_{\text{train}}.
\end{equation}
Firstly, we define the \emph{weighted correlation matrix} as
\begin{equation}
\boldsymbol{W}^{y}:=W \cdot D^{y},
\end{equation}
where $W:=\mathrm{diag}(w_{1},
\cdots, w_{N_{\text{train}}})$ 
is the diagonal matrix whose elements are the weights of \eqref{w-quadr}. 
The matrix $\boldsymbol{W}^{y}$ is not symmetric in the usual matrix sense, but with respect to the scalar product induced by $W^{y}$, hence $\boldsymbol{W}^{y}$ is diagonalizable anyway \cite{venturi2019weighted}.
Therefore, we seek the solution of the eigenvalue problem
$
\boldsymbol{W}^{y} g_{n}^{y}=\lambda_{n}^{y} g_{n}^{y}, \
1 \leq n \leq N_{\text{train}},
$
where $\left\|g_{n}^{y}\right\|_{{Y}}=1$, i.e.\ we pursue to find an eigenvalue $\lambda_{n}^{y}$ with the relative eigenvector of norm equal to one. We will indicate with $\left(g_{n}^{y}\right)_{t}$ the  $t$-th component of the eigenvector $g_{n}^{y} \in \mathbb{R}^{N_{\text{train}}}$. For the sake of simplicity, we rearrange the eigenvalues $\lambda_{1}^{y}, \ldots, \lambda_{N_{\text{train}}}^{y}$ in a decreasing layout. Then, let us look at the first $N$ eigenvalue-eigenvector pairs $(\lambda_{1}^{y},g_{1}^{y}), \ldots, (g_{N}^{y},\lambda_{N}^{y})$. The basis functions $\chi_{n}^{y}$ for the state equation are constructed through the following relation:
\begin{equation}\label{eigen-const}
\zeta_{n}^{y}=\frac{1}{\sqrt{\lambda_{n}^{y}}} \sum_{t=1}^{N_{\text{train}}}\left(g_{n}^{y}\right)_{t} y^{\mathcal{N}}\left(\boldsymbol{\mu}_{k}\right), \quad 1 \leq n \leq N.
\end{equation}
In order to choose $N$, one can refer to same study of eigenvalues of $\boldsymbol{W}^{y}$  \cite{hesthaven2016certified, quarteroni2015reduced, venturi2018weighted}.
At the end, our reduced spaces are built as \eqref{basis-red} and, then, enriched spaces are constructed.

We summarise the whole wPOD procedure for OCP($\boldsymbol{\mu}$)s in Algorithm \ref{alg:POD}. 
\begin{algorithm}
 \caption{Weighted POD algorithm for OCP($\boldsymbol{\mu}$) problems through a partitioned approach}\label{alg:POD}
\begin{flushleft}
\textbf{Input}: FEM spaces $Y^{\mathcal{N}}$, $U^{\mathcal{N}}$, and $(Q^{\mathcal{N}})^{*}$ parameter domain $\mathcal{P}$, and $N_{\text{train}}$. \\
\end{flushleft}
\begin{flushleft}
\textbf{Output:} reduced spaces $Y^{N}$, $U^{N}$ and $(Q^{N})^{*}$. \\
\end{flushleft}
\begin{flushleft}
Considering the high-fidelity spaces $Y^{\mathcal{N}}$, $U^{\mathcal{N}}$ and $(Q^{\mathcal{N}})^{*}$:
\end{flushleft}
\begin{algorithmic}[1]
\State Choose a quadrature rule \eqref{w-quadr} to approximate \eqref{expectation}. This step defines a sample $\mathcal{P}_{\text{train}} \subset \mathcal{P}$ and the respective weights $w_{1}, \ldots, w_{N_{\text{train}}}$. Define the matrix $W:=\mathrm{diag}(w_{1},
\cdots, w_{N_{\text{train}}})$ ;
\ForAll {$\boldsymbol{\mu} \in \mathcal{P}_{\text{train}}$} 
\State Solve the high-fidelity SUPG OCP($\boldsymbol{\mu}$) system \eqref{supg-system};
\EndFor
\State Calculate the matrices ${D}_{kl}^{y}:=\frac{1}{N_{\text{train}}}\left(y^{\mathcal{N}}\left(\boldsymbol{\mu}_{k}\right), y^{\mathcal{N}}\left(\boldsymbol{\mu}_{l}\right)\right)_{Y}, \ 1 \leq k,l \leq N_{\text{train}}$ and $\boldsymbol{W}^{y}:=W \cdot D^{y}$. Do the same for the control $u$ and the adjoint $p$;
\State Compute eigenvalues $\lambda_{1}^{y}, \ldots, \lambda_{N_{\text{train}}}^{y}$ and the corresponding orthonormalized eigenvectors $g_{1}^{y}, \ldots, g_{N_{\text{train}}}^{y}$ of $\boldsymbol{W}^{y}$. Do the same procedure for $u$ and $p$ variables;
\State {Fix $N$ according to a certain criterion} and construct $Y^{N}= \operatorname{span}\left\{\xi_{n}^{y}, n=1, \ldots, N\right\}$, where $\xi_{n}^{y}=\frac{1}{\sqrt{\lambda_{n}^{y}}} \sum_{t=1}^{N_{\text{train}}}\left(g_{n}^{y}\right)_{t} y^{\mathcal{N}}\left(\boldsymbol{\mu}_{k}\right)$. Do the same for $u$ and $p$ variables. 
\State Build the aggregated space $
G^{N}=\operatorname{span}\big\{\left\{\xi_{n}^{y}\right\}_{n=1}^{N} \cup\left\{\xi_{n}^{p}\right\}_{n=1}^{N}\big\}$ and set $G^N \equiv Y^N\equiv (Q^{N})^{*}$.
\end{algorithmic}
\end{algorithm}

\subsection{Online phase} \label{online phase}
 In this stage, all operations have usually a  $\mathcal{N}$-independent cost. This process reflects to be computationally cheap and, therefore, it can be recalled multiple times using small machine resources. Firstly, we choose a parameter $\boldsymbol{\mu}$. We get all the $\boldsymbol{\mu}$-independent quantities and reduced spaces back from the storage. Immediately, we {combined} parameter independent part with the $\boldsymbol{\mu}$-dependent ones, which are rapidly calculated here. Then a Galerkin projector onto $Y^{N},U^{N}$ and $(Q^{N})^{*}$ is performed, computing the reduced solution $y^{N}$, $u^{N}$ and $p^{N}$ through a reduced block matrix system. As previously seen in Section \ref{sec:supg}, a stabilization is needed in the truth approximation. However, it could also not be the case for the online stage. This scenario leads to two possibilities: we do not use SUPG in the online phase,  \textit{Offline-Only stabilization}, or, on the contrary, stabilization occurs also here \textit{Offline-Online stabilization}.

\section{Numerical Results}
\label{sec:results}

In this last part, we illustrate numerical simulations concerning two Advection-Dominated OCP($\boldsymbol{\mu}$)s under random inputs: the Graetz-Poiseuille Problem and the Propagating Front in a Square Problem. In both experiments, the parameter $\boldsymbol{\mu}$ will be a random vector and it will follow a prescribed probability density function that will be specified. The deterministic version of both experiments can be founded in \cite{zoccolan2024streamline}.

The Offline approximation will be always based on a $\mathbb{P}^{1}-$FEM, which means to consider a finite element method characterized by polynomials of degree less or equal than $1$.
In steady and unsteady simulations, the same stabilization parameter $\delta_K$ will be employed for both stabilization in the high-fidelity approximation and in the Online phase: namely in Offline-Online stabilization, $\delta_K$ is the same for both phases. 

For each simulation, \emph{relative errors} between the FEM and the reduced solutions, i.e.\
\begin{equation} \label{rel-error}
e_{y, N}(\boldsymbol{\mu}):=\frac{\left\|y^{\mathcal{N}}(\boldsymbol{\mu})-y^{N}(\boldsymbol{\mu})\right\|_{Y}}{\left\|y^{\mathcal{N}}(\boldsymbol{\mu})\right\|_{Y}}, \
e_{u, N}(\boldsymbol{\mu}):=\frac{\left\|u^{\mathcal{N}}(\boldsymbol{\mu})-u^{N}(\boldsymbol{\mu})\right\|_{U}}{\left\|u^{\mathcal{\mathcal{N}}}(\boldsymbol{\mu})\right\|_{U}}, \
e_{p, N}(\boldsymbol{\mu}):=\frac{\left\|p^{\mathcal{N}}(\boldsymbol{\mu})-p^{N}(\boldsymbol{\mu})\right\|_{Q^{*}}}{\left\|p^{\mathcal{N}}(\boldsymbol{\mu})\right\|_{Q^{*}}},
\end{equation}
for the state, the control and the adjoint, respectively, will be shown.
Due to the parametric nature of the problems, for each quantity in \eqref{rel-error} a simple average is computed for $\boldsymbol{\mu}$ {distributed} according to its probability density in a testing set $\mathcal{P}_{\text{test}} \subseteq \mathcal{P}$  of size $N_{\text{test}}$, for every dimension $N=1, \ldots, N_{\max}$ of the reduced space built through a chosen wPOD procedure. In every graph, the base-$10$ logarithm of these averages will be shown. When we will specify to use a POD procedure based on a Monte-Carlo sampling \cite{sullivan2015introduction} of a uniform density distribution, we will talk about \textit{Standard POD}. {In order to compare the different wPOD possibilities, we use the same testing set for all of them: it will be taken using a Monte-Carlo method according to the distribution of $\boldsymbol{\mu}$. Moreover, the same testing set will be used to compute the projection error of the FEM space into the POD one for all the three variables}. Obviously, the performance of the Standard POD will be based on a testing set of uniform density.
The sum of the errors with respect to each discretized instant of time $t$ will be taken into account in the unsteady versions.

In order to compare the computational cost between the FEM solution with that of the reduced one for any possible dimension $N$, we use the \textit{speedup-index}, i.e.\ 
\begin{equation}\label{speed-up}
\text{speedup-index}=\frac{\text{ computational time of the high-fidelity solution }}{\text{ computational time of the reduced solution}},
\end{equation}

which will be calculated for any $\boldsymbol{\mu}$ in the testing set. Again, we will show its sample average for any dimension $N$. For each test case, we will use the same $\mathcal{P}_{\text{test}}$ to compute relative errors and the speedup-index. The steady experiments are run using a machine with $16$GB of RAM and Intel Core i7-7500U Dual Core, $2.7$GHz for the CPU; whereas all parabolic simulations are computed considering $16$GB of RAM and Intel Core i$7-7700$ Quad Core, $3.60$GHz for the CPU. 

The code concerning steady experiments is implemented using the RBniCS library \cite{RBniCS}; instead, the unsteady ones are provided using both RBniCS and multiphenics \cite{multiphenics} libraries. These are Python-based libraries, built on FEniCS \cite{logg2012automated}.

\subsection{Numerical Tests for the Graetz-Poiseuille Problem}
\label{sec:c graetz}

The Graetz-Poiseuille problem is an Advection-Diffusion problem that represents the heat conduction in a rectilinear pipe. Here the transfer of heat can be regulated through the walls of the duct, which can be held at a certain fixed temperature \cite{gelsomino2011comparison,pacciarini2014stabilized, rozza2009reduced, torlo2018stabilized}. 

Firstly, we present simulation concerning the stationary case, where a distributed control is employed all over the whole domain. The parameter $\boldsymbol{\mu} = (\mu_1, \mu_2)$ is composed by the diffusion component $\mu_1$ and the geometrical one $\mu_2$, which characterizes the length of the plate.
\begin{figure}[h!]
   \centering
        \begin{tikzpicture}[scale=3.0]
                 \draw[color=DEblue!100, fill=DEblue!10] (1,1) -- (2.7, 1) -- (2.7,0.8) -- (1,0.8) --(1,1)node[midway, left, scale=1.2]{};	
                 \draw[color=DEblue!100, fill=DEblue!10] (1,0) -- (2.7, 0) -- (2.7,0.2) -- (1,0.2) --(1,0)node[midway, left, scale=1.2]{};
                 \draw[black] (1.5,0.9) node[scale=1.]{{$\Omega_{obs}$}};
                 \draw[black] (1.5,0.1) node[scale=1.]{{$\Omega_{obs}$}};
                 \draw[black] (1.35,0.5) node[scale=1.5]{{$\Omega_o$}};
                  \draw[blue] (0,0) -- (1,0) node[midway, below, scale=1.2]{$\Gamma_{o,1}$};
                 \draw[red] (1,0) -- (2.7,0) node[midway, below, scale=1.2]{$\Gamma_{o,2}$};
                 \draw[dashed,red] (2.7,0) -- (2.7,1) node[midway, right, scale=1.2]{$\Gamma_{o,3}$};
                 \draw[red] (2.7,1) -- (1,1) node[midway, above, scale=1.2]{$\Gamma_{o,4}$};
                 \draw[blue] (1,1) -- (0,1) node[midway, above, scale=1.2]{$\Gamma_{o,5}$};
                 \draw[blue] (0,1) -- (0,0) node[midway, left, scale=1.2]{$\Gamma_{o,6}$};
                 \filldraw[black] (0,0) circle (0.3pt) node[left]{(0,0)};
                 \filldraw[black] (1,0) circle (0.3pt) node[below]{(1,0)};
                 \filldraw[black] (2.7,0) circle (0.3pt) node[below]{(1+$\mu_2$,0)};
                 \filldraw[black] (2.7,0.2) circle (0.3pt) node[right]{(1+$\mu_2$,0.2)};
                 \filldraw[black] (2.7,0.8) circle (0.3pt) node[right]{(1+$\mu_2$,0.8)};
                 \filldraw[black] (2.7,1) circle (0.3pt) node[above]{(1+$\mu_2$,1)};
                 \filldraw[black] (1,1) circle (0.3pt) node[above]{(1,1)};
                 \filldraw[black] (0,1) circle (0.3pt) node[left]{(0,1)};
        \end{tikzpicture}
        \caption{Geometry of the Graetz-Poiseuille Problem.}
       \label{fig:Geometry-Graetz}
\end{figure}
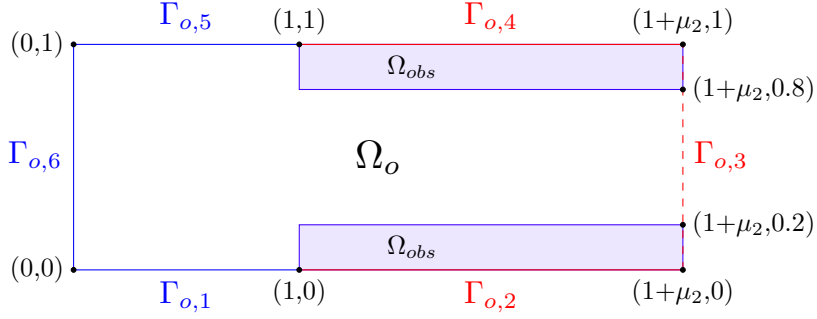

{The problem is studied using $(x_0, x_1)$ as spatial coordinates.}
{$\Omega_o$ is the domain observed for a certain value $\mu_2$ with boundary $\Gamma_o$.} We deal with homogeneous Neumann boundary conditions (BC) on $\Gamma_{o,3}:=\{1+\mu_2 \} \times [0,1]$ considering Figure \ref{fig:Geometry-Graetz}. Instead, 
Dirichlet conditions are set on sides $\Gamma_{o,1}:=[0,1]\times \{0\}$, $\Gamma_{o,5}:= [0,1]\times \{1\}$, $\Gamma_{o,6}:=\{0\}\times [0,1]$ {by imposing} $y=0$ and $\Gamma_{o,2}:=[1,1+\mu_2]\times \{0\}$ and $\Gamma_{o,4}:=[1,1+\mu_2]\times \{1\}$ {by imposing} $y=1$. \\

The formulation of the problem is the following:
given $\boldsymbol{\mu} \in \mathcal{P}$, find $(y, u) \in \tilde{Y} \times U$ which solves
\begin{equation*}
\min_{(y,u)} \displaystyle \frac{1}{2} \int_{\Omega_{obs}(\boldsymbol{\mu})}(y(\boldsymbol{\mu} )- y_d)^2 \; d \Omega_o(\boldsymbol{\mu}) \vspace{.5cm} +
\frac{\alpha}{2} \int_{\Omega_o(\boldsymbol{\mu})}u(\boldsymbol{\mu} )^2 \; d \Omega_o(\boldsymbol{\mu}), \ \text{ such that }
\vspace{-5mm}
\end{equation*}
\begin{equation} \label{graetz-system}
    \begin{cases}
       \displaystyle -\frac{1}{{\mu_{1}}} \Delta y(\boldsymbol{\mu})+4 x_1(1-x_1) \partial_{x_0} y(\boldsymbol{\mu})=u(\boldsymbol{\mu}), & \text { in } \Omega_{o}(\boldsymbol{\mu}), \\
\displaystyle y(\boldsymbol{\mu})=0, & \text { on } \Gamma_{o, 1}(\boldsymbol{\mu}) \cup \Gamma_{o, 5}(\boldsymbol{\mu}) \cup \Gamma_{o, 6}(\boldsymbol{\mu}), \\
\displaystyle y(\boldsymbol{\mu})=1, & \text { on } \Gamma_{o, 2}(\boldsymbol{\mu}) \cup \Gamma_{o, 4}(\boldsymbol{\mu}), \\
\displaystyle \frac{\partial y(\boldsymbol{\mu})}{\partial \nu}=0, & \text { on } \Gamma_{o, 3}(\boldsymbol{\mu}),
    \end{cases}
\end{equation}
where $\tilde{Y}:= \big\{v \in H^{1} \big(\Omega_o\big) \text{ s.t. it satisfies the } \mathrm{BC} \text{ in }  (\ref{graetz-system}) \big\}$ and $U=L^2(\Omega_o)$. 
For the sake of clarity, a lifting function $R_{y} \in H^{1}(\Omega)$  that fulfills the $\mathrm{BC} \text{ in }  (\ref{graetz-system})$ is used. Consequently, the variable $\bar{y} := y - R_y$, with $\bar{y} \in Y$, is used, where
$$Y:= \big\{ v \in H^{1}_{0}\big(\Omega\big) \text{ s.t. } \frac{\partial \bar{y}}{\partial \nu}=0, \text { on } \Gamma_{3} \text{ and } \bar{y}=0 \text{ on } \Gamma \setminus \Gamma_{3} \big\}.$$ Furthermore, we settle $Q := Y^{*}$ without any loss of generality. Therefore, the adjoint variable $p$ is null on $ \Gamma \setminus \Gamma_{3}$. The observation domain is $\Omega_{obs}:= [1,1+\mu_2] \times [0.8,1] \cup [1,1+\mu_2] \times [0,0.2]$ as illustrated in Figure \ref{fig:Geometry-Graetz}.  The value $\mu_2$ can change the domain under study. Having that the domain $\Omega_o$ is $\boldsymbol{\mu}$-dependent itself, in the Offline Phase snapshots are based on different domains due to the sampling of the geometrical parameter components \cite{hesthaven2016certified,quarteroni2011certified,rozza2008reduced,rozza2011reduced}. To deal with the \textit{geometrical parametrization} of the problem, we set a reference domain $\Omega$ and we build affine maps that transform $\Omega$ in $\Omega_o$ for a defined $\boldsymbol{\mu}$. This procedure implies an automatic modification of some bilinear and linear forms involved in the weak formulation of Problem \eqref{graetz-system}.

We choose $\Omega = (0,2) \times (0,1)$ as reference domain, that is the original one $\Omega_o(\boldsymbol{\mu})$ corresponding to ${\mu_2}=1$. We assume that $\mu_2$ is positive for the sake of simplicity.
Considering Figure \ref{fig:Geometry-Graetz}, we {divide} this into $2$ subdomains, which are defined as $\Omega_{1}=(0,1) \times(0,1)$ and $\Omega_{2}=(1,2) \times(0,1)$. Then, we build two affine transformations:
\begin{equation}\label{V1}
V_{1}(\boldsymbol{\mu}): \Omega_{1} \rightarrow  \Omega_{o, 1}(\boldsymbol{\mu}) \subset \mathbb{R}^{2}, \ \text{ such that } \
V_{1}\left(\left(\begin{array}{l}
x \\
y
\end{array}\right) ; \boldsymbol{\mu}\right):=\left(\begin{array}{l}
x \\
y
\end{array}\right),
\end{equation}
which is the identity map defined on the first subdomain $\Omega_{1}$ and $V_{2}(\boldsymbol{\mu}): \Omega_{2} \rightarrow \Omega_{o, 2}(\boldsymbol{\mu}) \subset \mathbb{R}^{2}$ as
\begin{equation}\label{V2}
V_{2}\left(\left(\begin{array}{l}
x \\
y
\end{array}\right) ; \boldsymbol{\mu}\right)=\left(\begin{array}{c}
\mu_{2} x \\
y
\end{array}\right)+\left(\begin{array}{c}
1-\mu_{2} \\
0
\end{array}\right)=R_{2}\left(\begin{array}{l}
x \\
y
\end{array}\right)+\left(\begin{array}{c}
1-\mu_{2} \\
0
\end{array}\right),
\end{equation}
where we have
\begin{equation}
R_{2}:=\left(\begin{array}{cc}
\mu_{2} & 0 \\
0 & 1
\end{array}\right).
\end{equation}
Glueing together $V_1$ and $V_2$ for each $\boldsymbol{\mu} \in \mathcal{P}$, we manage to build a one-to-one transformation $V(\boldsymbol{\mu})$ defined all over $\Omega$. We denote the restrictions of $\mathcal{T}_{h}$ to $\Omega_{1}$ and $\Omega_{2}$ with $\mathcal{T}_{h}^{1}$ and $\mathcal{T}_{h}^{2}$, respectively.
Therefore, we can express all the forms of the weak formulation under the effect of this transformation. For instance, after possible lifting, we have $a_s=a+s$ and $a^{*}_s=a^{*}+s^{*}$  as
\begin{equation*}
\begin{aligned}
a\left(y^{\mathcal{N}}, q^{\mathcal{N}} ; \boldsymbol{\mu}\right):&= \int_{\Omega^{1}} \frac{1}{\mu_{1}} \nabla y^{\mathcal{N}} \cdot \nabla q^{\mathcal{N}} +4  x_1(1-x_1)  \partial_{x_0} y^{\mathcal{N}} q^{\mathcal{N}} \\
&+\int_{\Omega^{2}} \frac{1}{\mu_{1} \mu_{2}} \partial_{x_0}y^{\mathcal{N}} \partial_{x_0} q^{\mathcal{N}}+\frac{\mu_{2}}{\mu_{1}} \partial_{x_1} y^{\mathcal{N}} \partial_{x_1} q^{\mathcal{N}}+4 x_1(1-x_1) \partial_{x_0} y^{\mathcal{N}} q^{\mathcal{N}}, \\
s\left( y^{\mathcal{N}}, q^{\mathcal{N}} ; \boldsymbol{\mu}\right):&=  \sum_{K \in \mathcal{T}_{h}^{1}}\delta_K h_{K} \int_{K}\left(4 x_1(1-x_1) \partial_{x_0}  y^{\mathcal{N}}\right) \partial_{x_0} q^{\mathcal{N}} \\
&+ \sum_{K \in \mathcal{T}_{h}^{2}} \delta_K \frac{h_{K}}{\sqrt{\mu}_{2}} \int_{K}\left(4 x_1(1-x_1) \partial_{x_0}  y^{\mathcal{N}}\right) \partial_{x_0} q^{\mathcal{N}}, 
\end{aligned}
\end{equation*}
\begin{equation*}
\begin{aligned}
a^{*}\left(z^{\mathcal{N}}, p^{\mathcal{N}} ; \boldsymbol{\mu}\right):&= \int_{\Omega^{1}} \frac{1}{\mu_{1}} \nabla p^{\mathcal{N}} \cdot \nabla z^{\mathcal{N}} -4  x_1(1-x_1)  \partial_{x_0} p^{\mathcal{N}} z^{\mathcal{N}} \\
&-\int_{\Omega^{2}} \frac{1}{\mu_{1} \mu_{2}} \partial_{x_0}p^{\mathcal{N}} \partial_{x_0} z^{\mathcal{N}}-\frac{\mu_{2}}{\mu_{1}} \partial_{x_1} p^{\mathcal{N}} \partial_{x_1} z^{\mathcal{N}}-4 x_1(1-x_1) \partial_{x_0} p^{\mathcal{N}} z^{\mathcal{N}}, \\
s^{*}\left( z^{\mathcal{N}}, p^{\mathcal{N}} ; \boldsymbol{\mu}\right):&=  \sum_{K \in \mathcal{T}_{h}^{1}}\delta_K h_{K} \int_{K}\left(4 x_1(1-x_1) \partial_{x_0}  p^{\mathcal{N}}\right) \partial_{x_0} z^{\mathcal{N}} \\
&+ \sum_{K \in \mathcal{T}_{h}^{2}} \delta_K \frac{h_{K}}{\sqrt{\mu}_{2}} \int_{K}\left(4 x_1(1-x_1) \partial_{x_0}  p^{\mathcal{N}}\right) \partial_{x_0} z^{\mathcal{N}},
\end{aligned}
\end{equation*}
\emph{for all} $y^{\mathcal{N}}, q^{\mathcal{N}}, z^{\mathcal{N}}, p^{\mathcal{N}},  \in Y^{\mathcal{N}}.$
In order to take into account the possible bad effect on stabilized forms due to an extension or shortening of our domain $\Omega_o$, we choose the stabilization parameter for $K \in \mathcal{T}_{h}^{2}$ as $\delta_K \frac{h_{K}}{\sqrt{\mu_{2}}}$, where $\sqrt{\mu_{2}}= \sqrt{|\operatorname{det}(R_{2})|}$ \cite{negri2013reduced, pacciarini2014stabilized, rozza2016stabilized}.

For the FEM discretization, a quite coarse mesh of size $h=0.034$ is used and the total dimension of the numerical problem is $13146$. We take $\delta_K =1.0$ \emph{for all} $K \in \mathcal{T}_{h}$.
The parameter space is set as $\mathcal{P}:=\big[1,10^5\big]\times \big[0.5, 1.5\big]$, from which we want to extract a
training set $\mathcal{P}_{\text{train}}$ with cardinality $N_{\text{train}}=100$. For the $n$ bilinear form, we consider a penalization $\alpha=0.01$. Our aim is to minimize the $L^2$-error between the state and the desired solution profile $y_d(x)=1.0$, function defined on $\Omega_{obs}$ of Figure \ref{fig:Geometry-Graetz}. 
Each wPOD procedure is computed until a $N_{\max}=20$ in a partitioned approach and then all algorithms are compared using a testing set $\mathcal{P}_{\text{test}}$ of $100$ elements in $\mathcal{P}$. 

\begin{figure}
        \centering
        \includegraphics[scale=0.23]{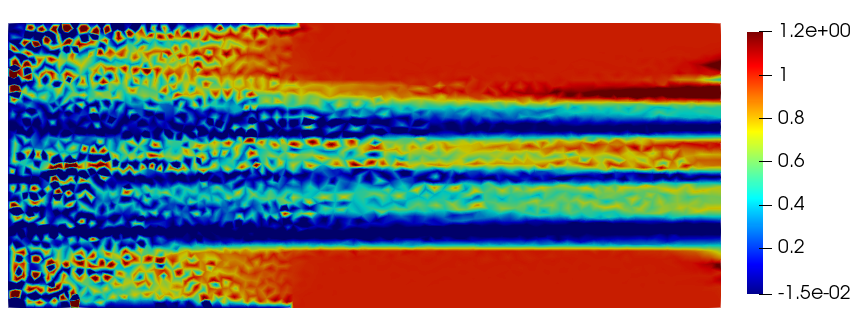}
        \includegraphics[scale=0.23]{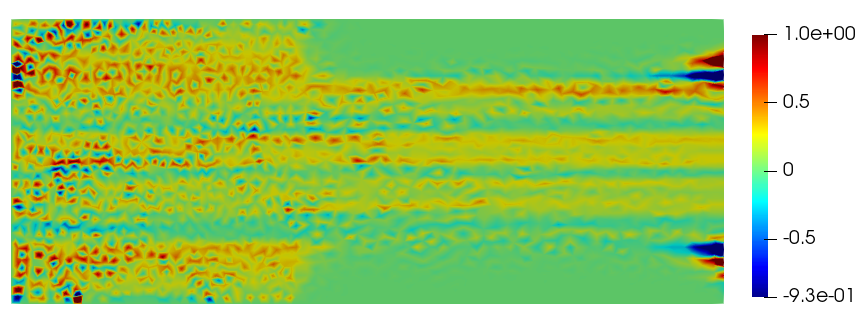} \\
        \includegraphics[scale=0.23]{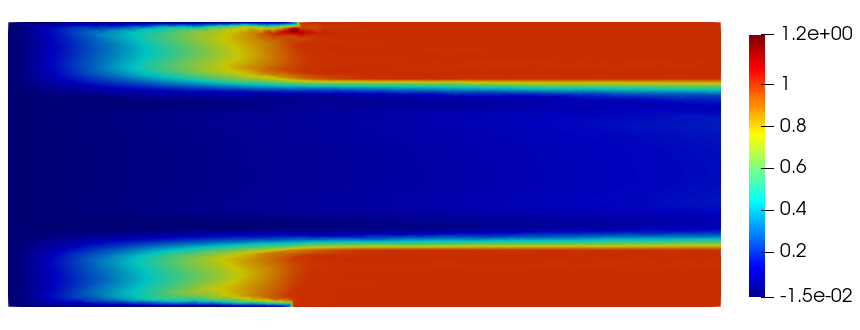}
        \includegraphics[scale=0.23]{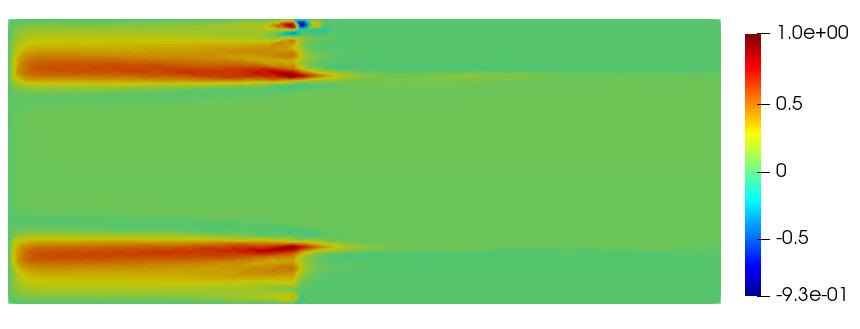} 
        \caption{(\underline{top}) FEM not stabilized and (\underline{bottom}) FEM stabilized solution, $y$ (right) and $u$ (left), $\boldsymbol{\mu} = (10^{5},1.5)$,  $h=0.034$, $\alpha=0.01$, $\delta_K=1.0$.}
        \label{fig:Graetz_FEM}
\end{figure}

We suppose that $\boldsymbol{\mu}$ follows a $\operatorname{Beta}(5,3)$ distribution for both parameter $\mu_1$ and $\mu_2$, i.e.\
\begin{equation}\label{beta-graetz}
\begin{aligned}
    {\mu_1} \sim 1 + \big(10^5 - 1\big) X_1, \text{ where } X_1 \sim \text{Beta}(5,3), \\
     \mu_2 \sim 0.5 + \big( 1.5 - 0.5\big) X_2, \text{ where } X_2 \sim \text{Beta}(5,3),
\end{aligned}
\end{equation}
where $\mu_1$ and $\mu_2$ are independent random variables.
This implies that we consider more probable the parameters for which the Graetz-Poiseuille Problem has high values of the P\'eclet number. In Figure \ref{fig:Graetz_FEM}, we highlight how the FEM solutions of the state and the control are for $\boldsymbol{\mu} = (10^{5},1.5)$. The adjoint solution is not shown here because it is proportional to the control due to the gradient equation \cite{collis2002analysis}. From Figure \ref{fig:Graetz_FEM}, one can see that a stabilization is necessarily needed.  \sloppy

\begin{figure}
        \centering
        \includegraphics[scale=0.2]{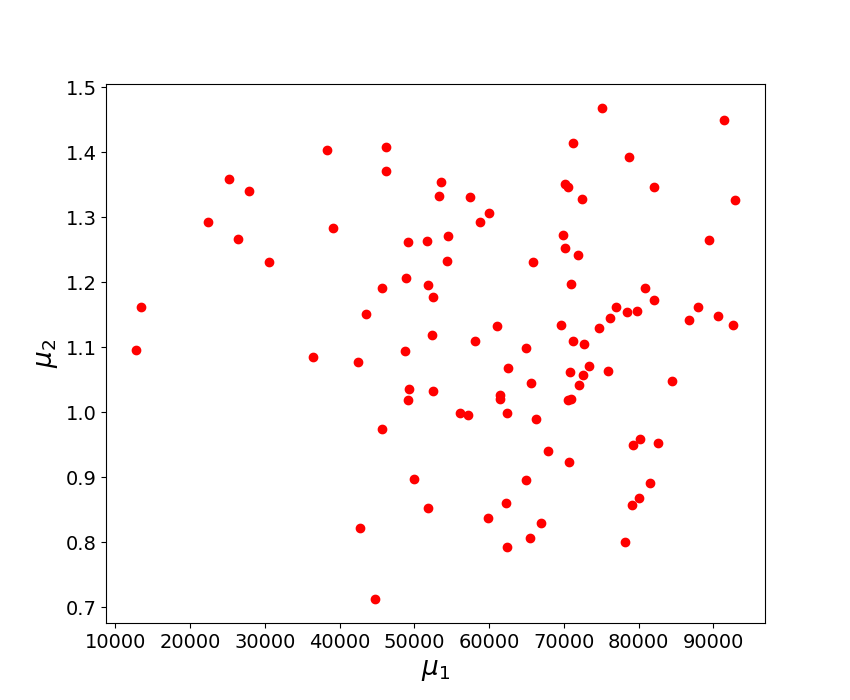}  
        \includegraphics[scale=0.2]{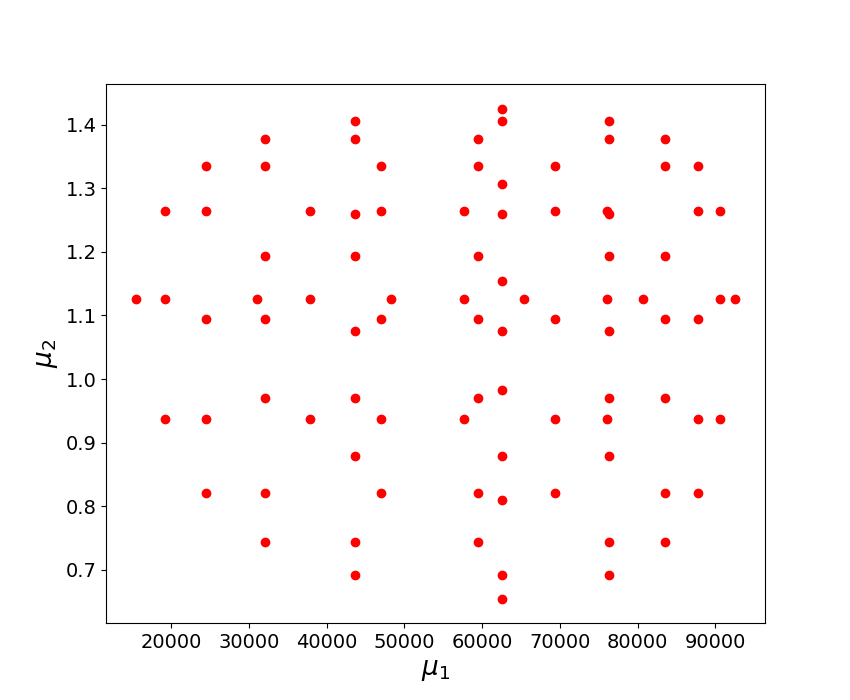} 
         \includegraphics[scale=0.2]{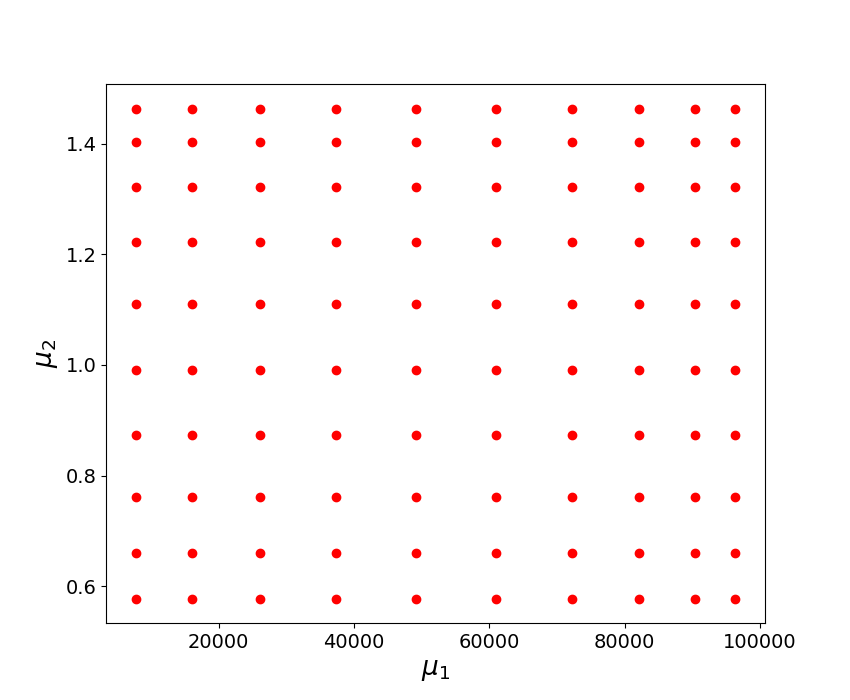} \\
        \includegraphics[scale=0.2]{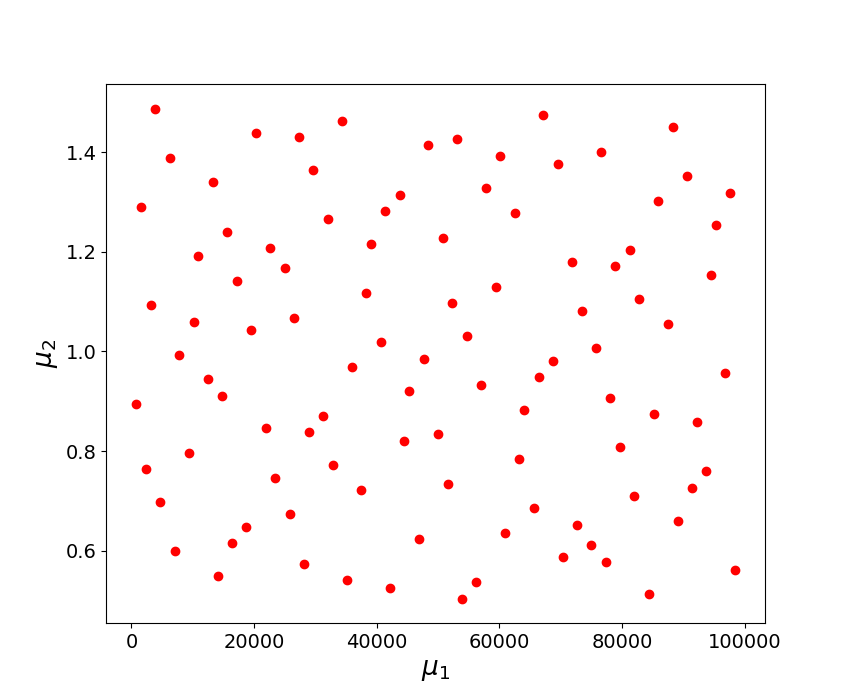} 
        \includegraphics[scale=0.2]{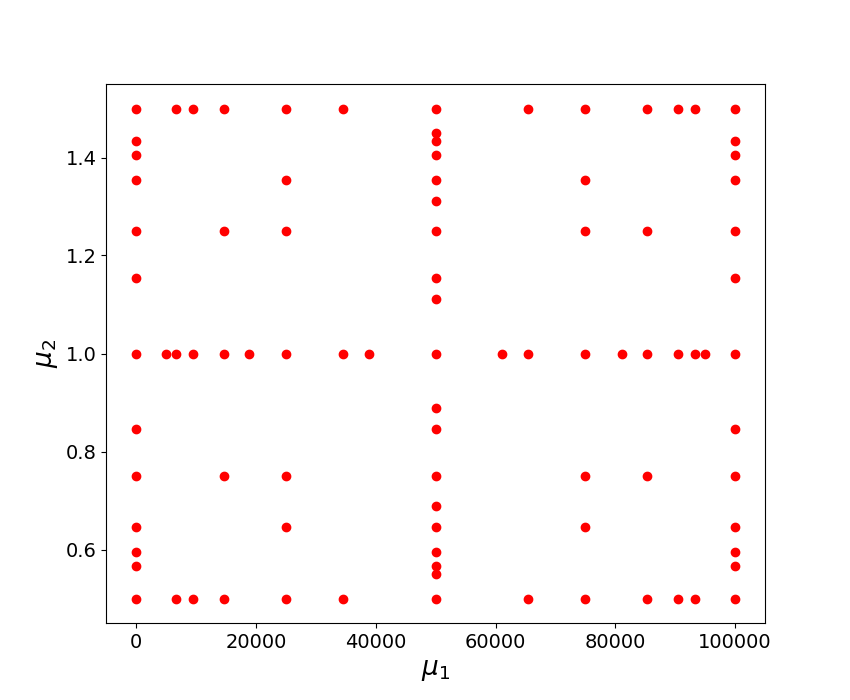} 
        \includegraphics[scale=0.2]{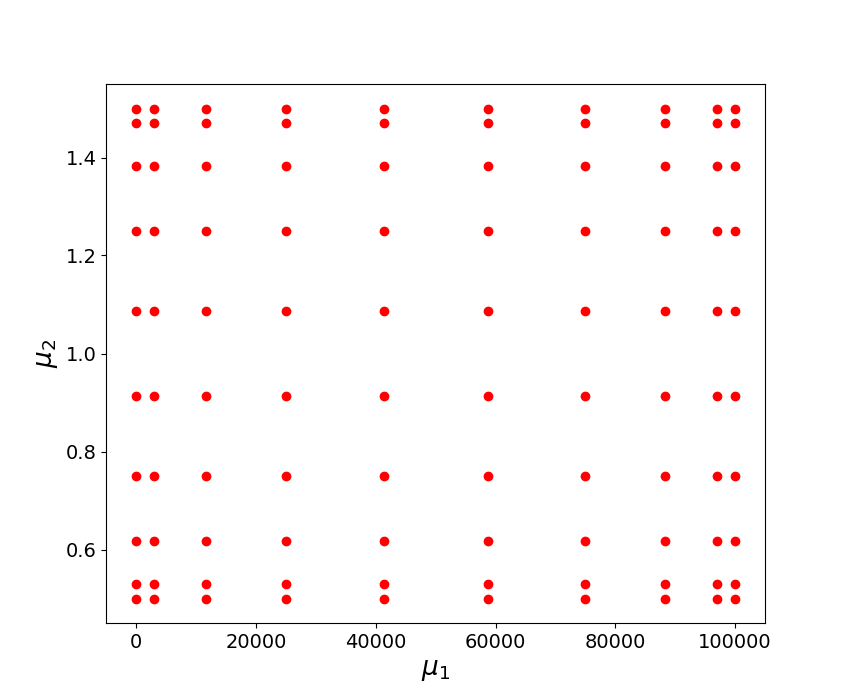} 
        \caption{Grid points for the quadrature formulae of the weighted POD regarding the Graetz-Poiseuille Problem; (\underline{top}) Monte-Carlo method  with $\boldsymbol{\mu}$ following distribution \eqref{beta-graetz} on the parameter space $\mathcal{P}$ (\underline{left}), Smolyak grid based on a Gauss-Jacobi rule (\underline{center}), Tensor grid based on a Gauss-Jacobi rule (\underline{right}); (\underline{bottom}) Pseudo-Random method based on a Halton sequence (\underline{left}), Smolyak grid based on a Clenshaw-Curtis rule (\underline{center}), Tensor grid based on a Clenshaw-Curtis rule (\underline{right}).}
        \label{fig:graetz-grids}
\end{figure}

\begin{figure}
        \centering
        \includegraphics[scale=0.16]{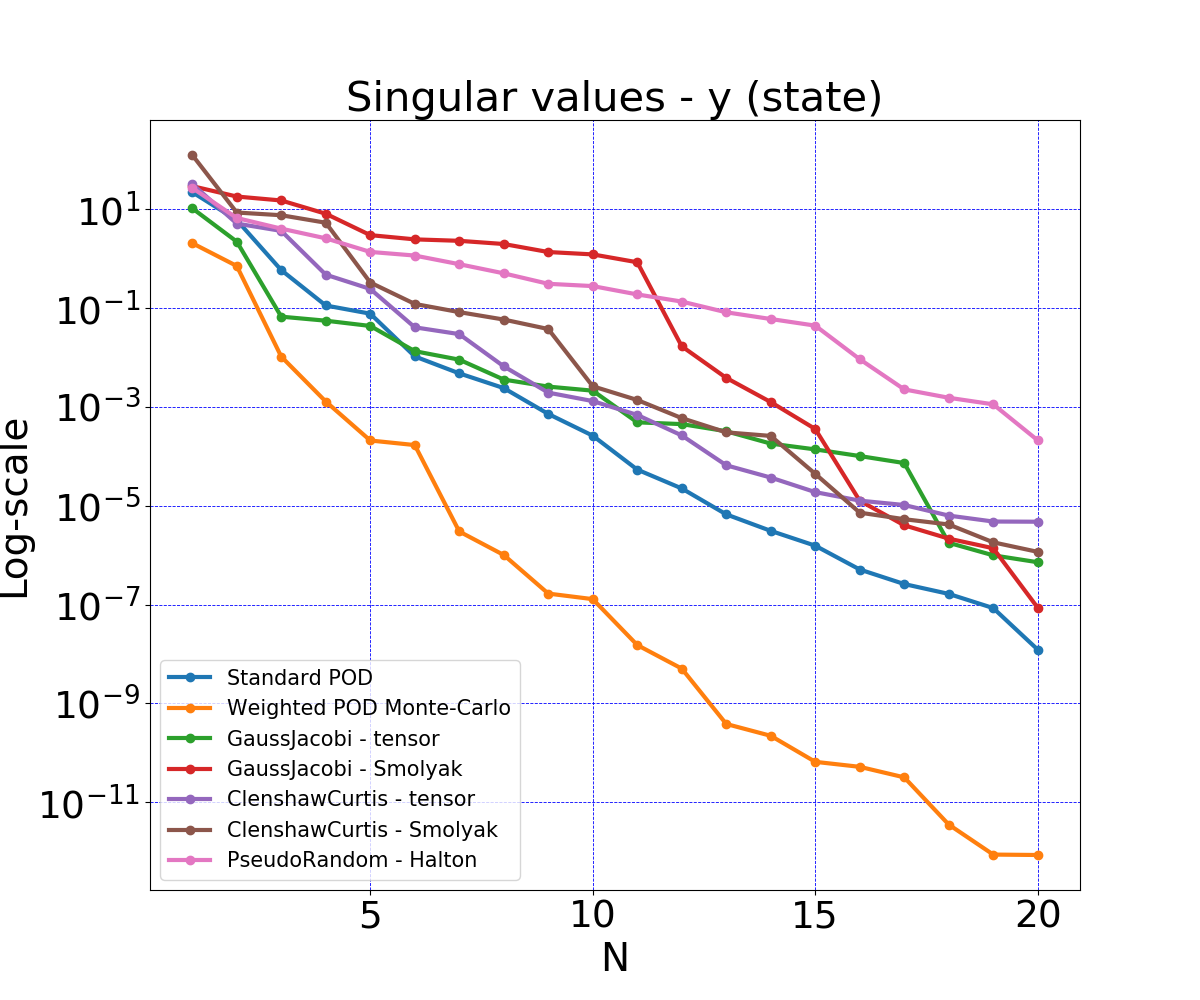} 
        \includegraphics[scale=0.16]{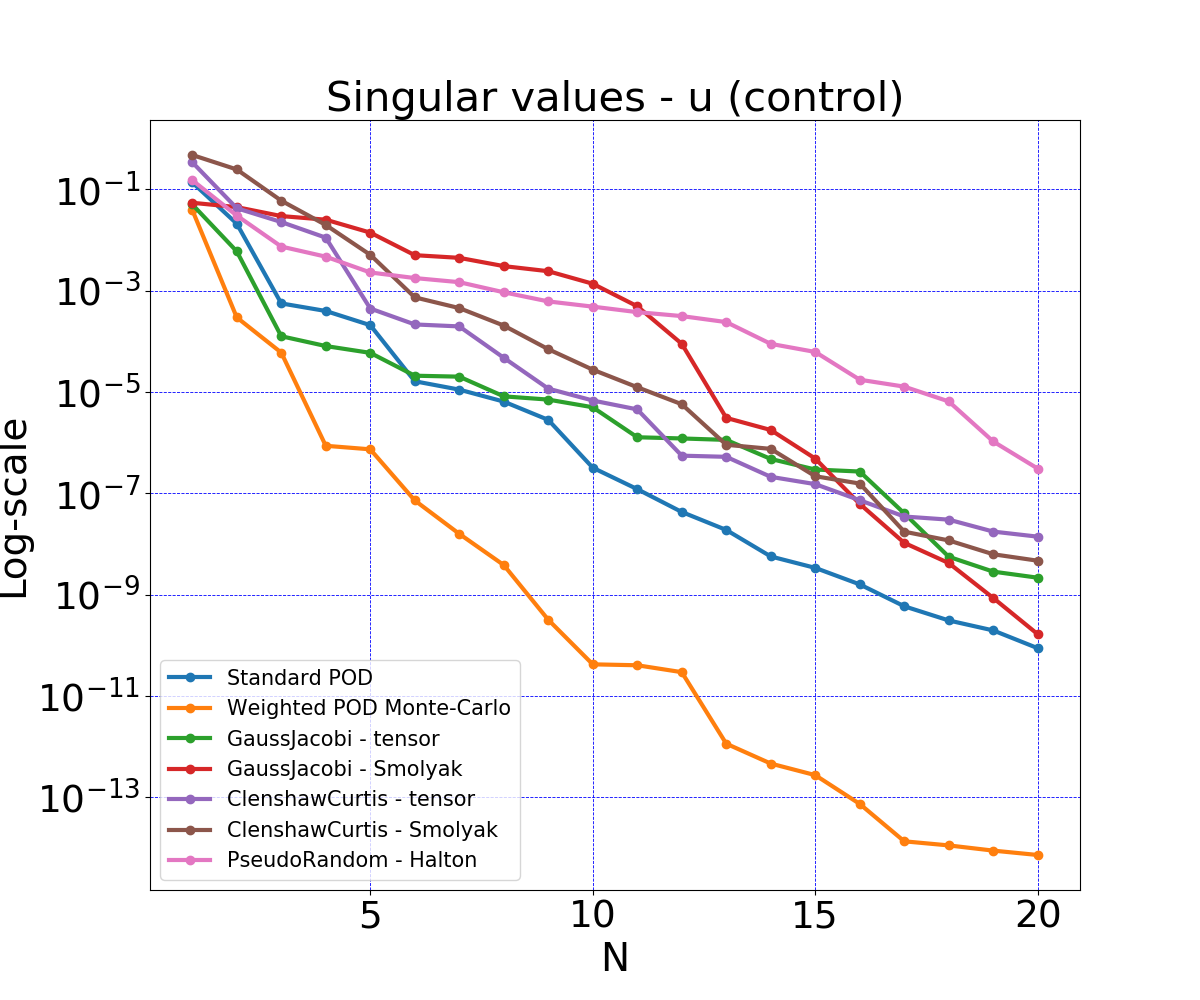} 
        \includegraphics[scale=0.16]{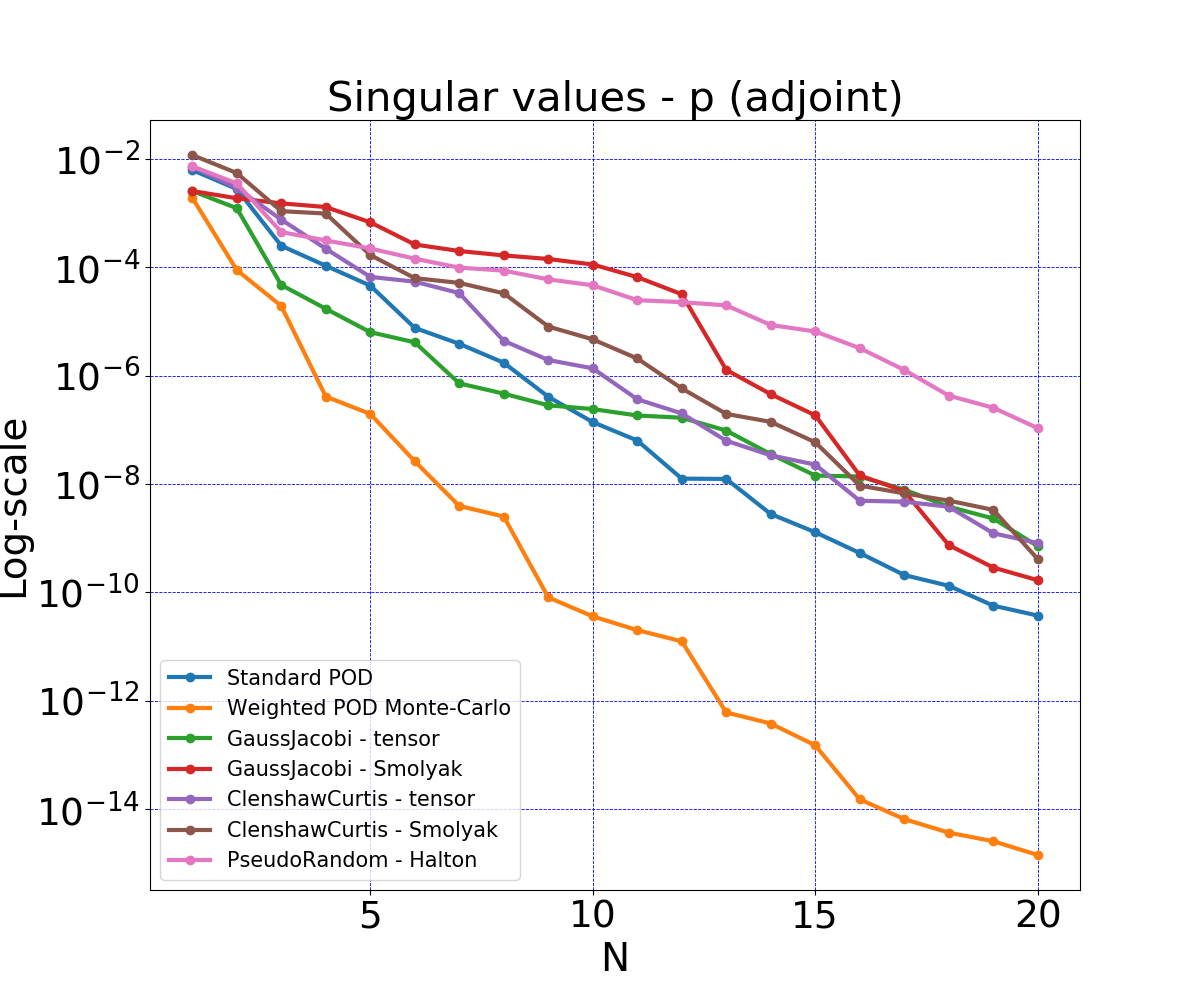} 
        \caption{Singular values decay for the snapshot matrices for the Graetz-Poiseuille Problem with $\boldsymbol{\mu}$ following distribution \eqref{beta-graetz} on the parameter space $\mathcal{P}$; State (\underline{left}), Control (\underline{center}), Adjoint (\underline{right}); Standard POD (blue), wPOD Monte-Carlo (orange), Gauss-Jacobi tensor rule (green), Gauss-Jacobi Smolyak grid (red), Clenshaw-Curtis tensor rule (cyan), Clenshaw-Curtis Smolyak grid (dark green), Pseudo-Random based on Halton numbers (pink).}
        \label{fig:plot_Graetz Geom-sing values}
\end{figure}

\begin{figure}
        \centering
        \includegraphics[scale=0.16]{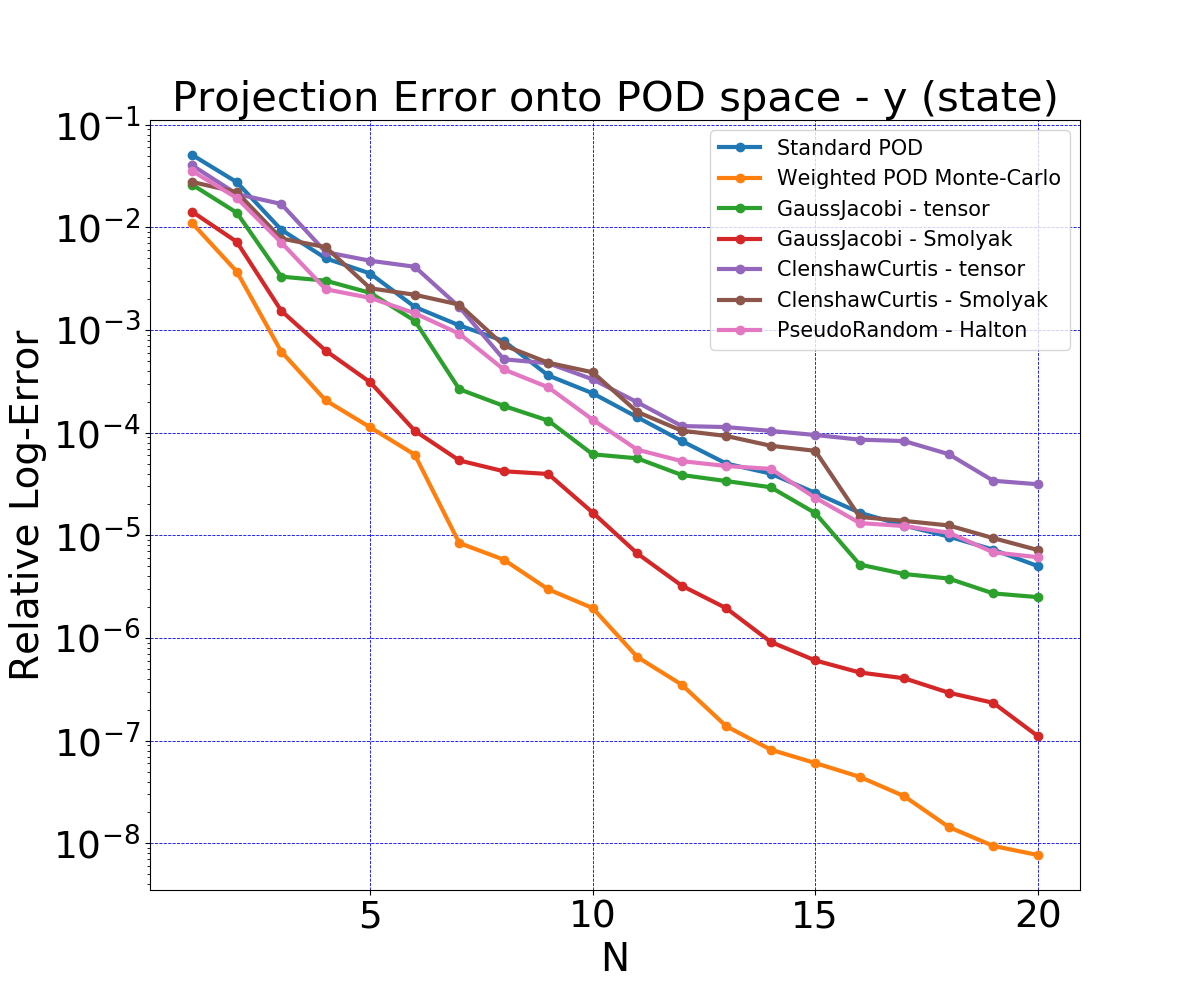} 
        \includegraphics[scale=0.16]{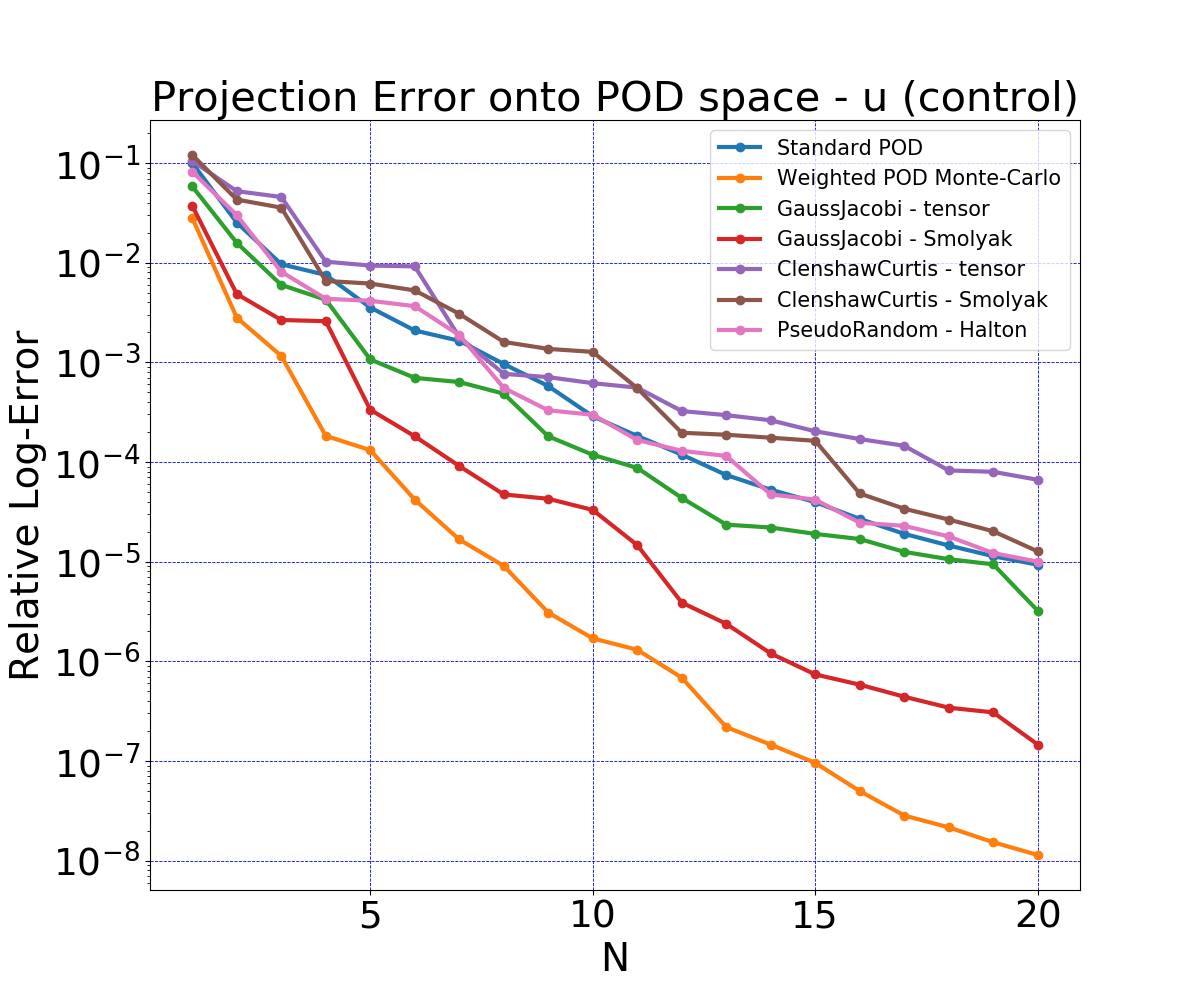} 
        \includegraphics[scale=0.16]{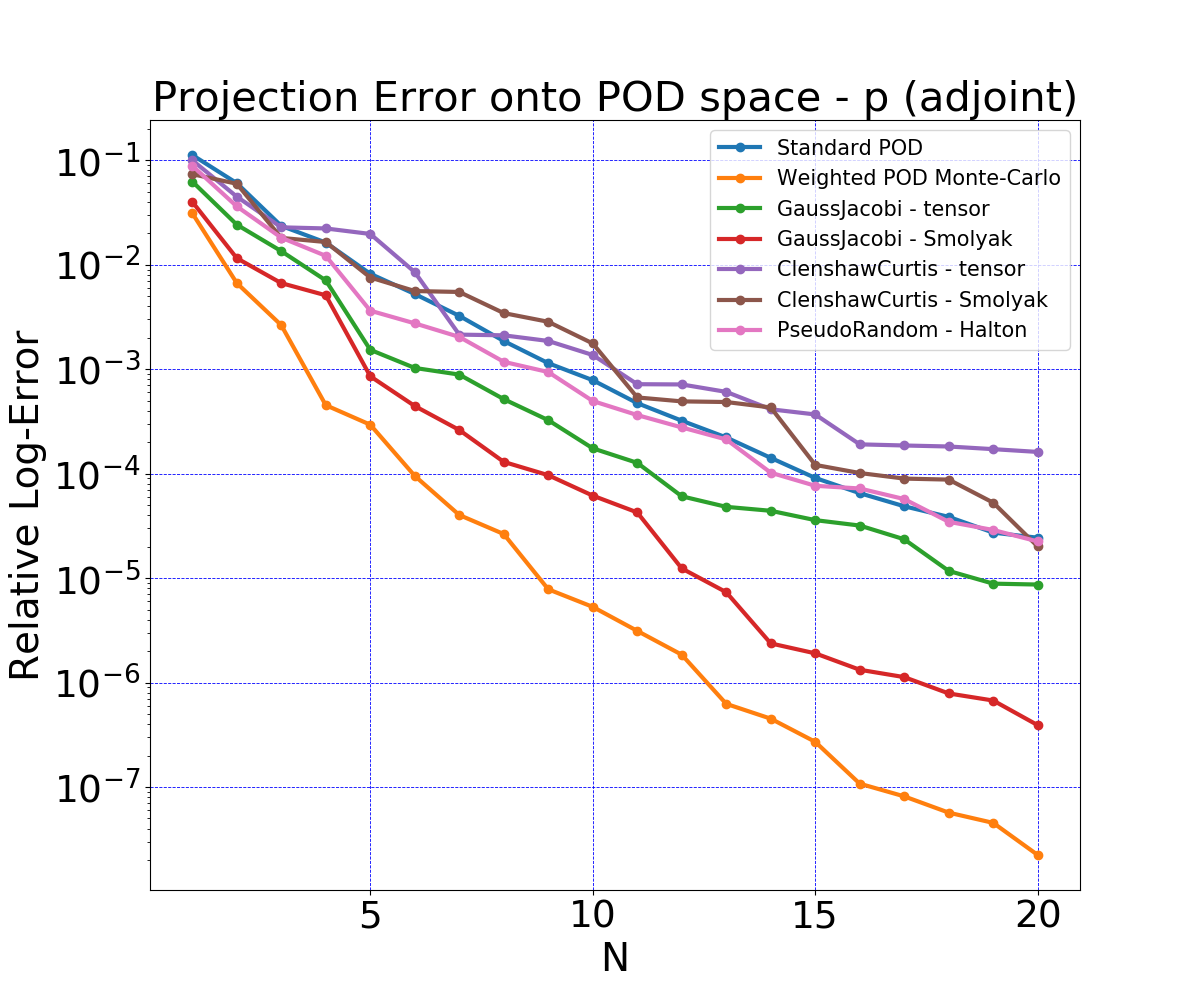} 
        \caption{Projection Errors onto the POD space for the Graetz-Poiseuille Problem with $\boldsymbol{\mu}$ following distribution \eqref{beta-graetz} on the parameter space $\mathcal{P}$; State (\underline{left}), Control (\underline{center}), Adjoint (\underline{right}); Standard POD (blue), wPOD Monte-Carlo (orange), Gauss-Jacobi tensor rule (green), Gauss-Jacobi Smolyak grid (red), Clenshaw-Curtis tensor rule (cyan), Clenshaw-Curtis Smolyak grid (dark green), Pseudo-Random based on Halton numbers (pink).}
        \label{fig:plot_Graetz Geom-proj error}
\end{figure}

\begin{figure}
        \centering
        \includegraphics[scale=0.121]{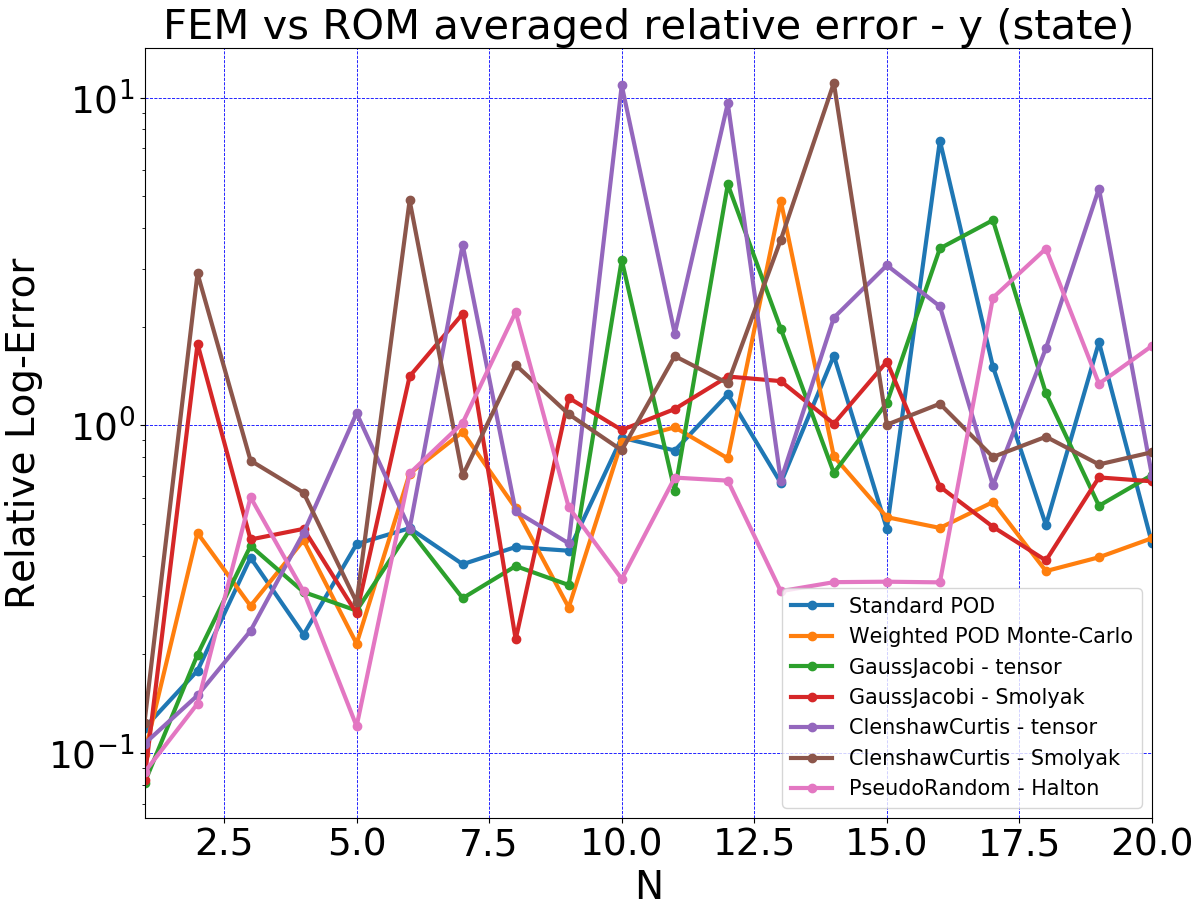} 
        \includegraphics[scale=0.121]{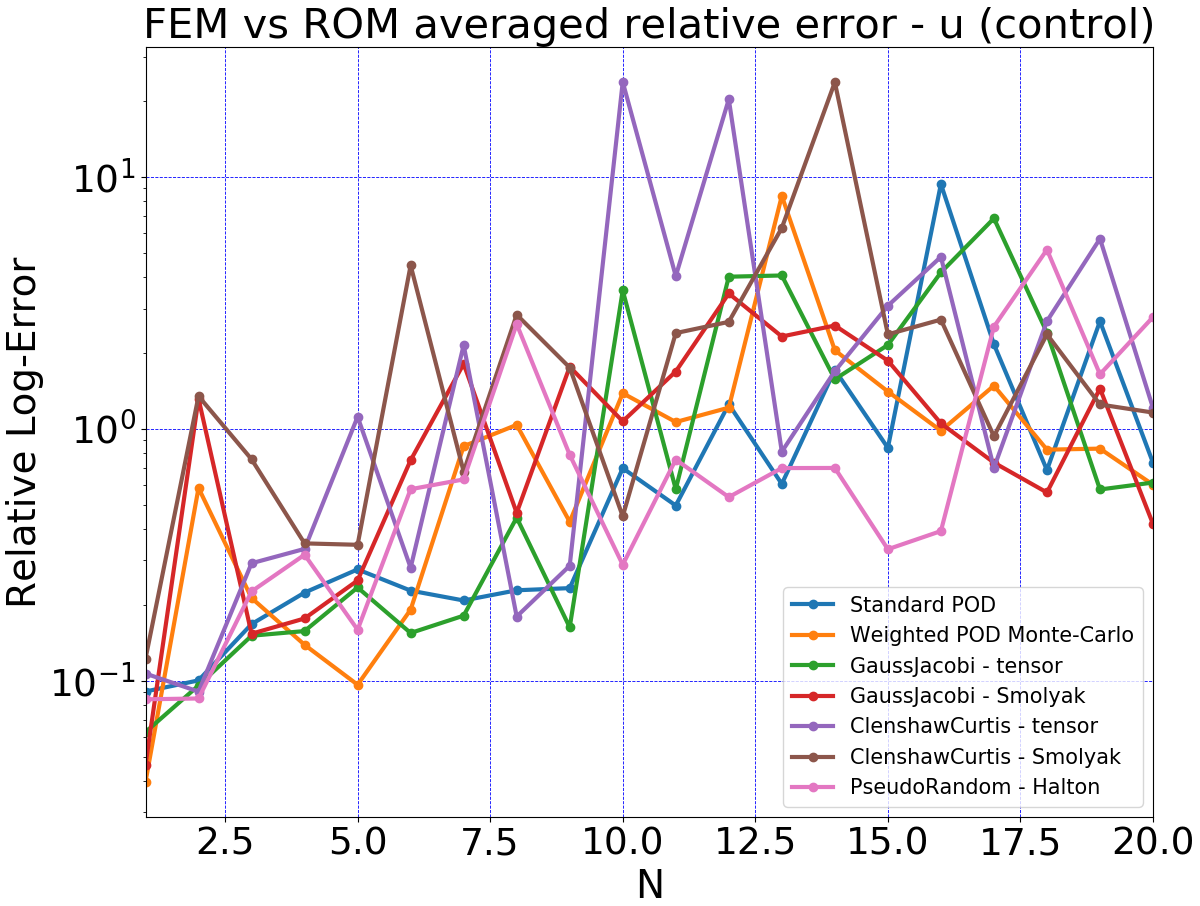} 
        \includegraphics[scale=0.121]{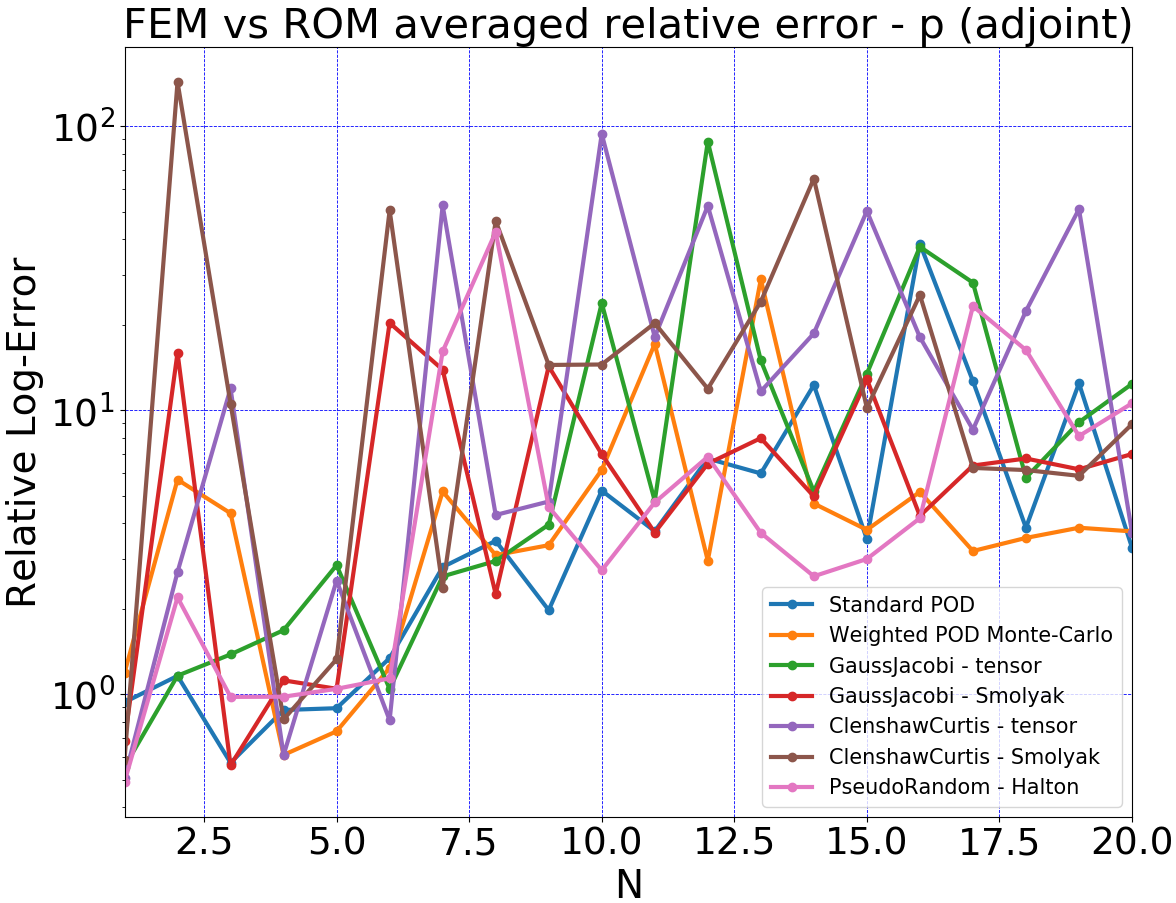} 
        \caption{Relative Errors for the Graetz-Poiseuille Problem with $\boldsymbol{\mu}$ following distribution \eqref{beta-graetz} on the parameter space $\mathcal{P}$ - \textit{Offline-Only} Stabilization; State (\underline{left}), Control (\underline{center}), Adjoint (\underline{right}); Standard POD (blue), wPOD Monte-Carlo (orange), Gauss-Jacobi tensor rule (green), Gauss-Jacobi Smolyak grid (red), Clenshaw-Curtis tensor rule (cyan), Clenshaw-Curtis Smolyak grid (dark green), Pseudo-Random based on Halton numbers (pink).}
        \label{fig:plot_Graetz Geom-not stab}
\end{figure}

\begin{figure}
    \centering
    \includegraphics[scale=0.12]{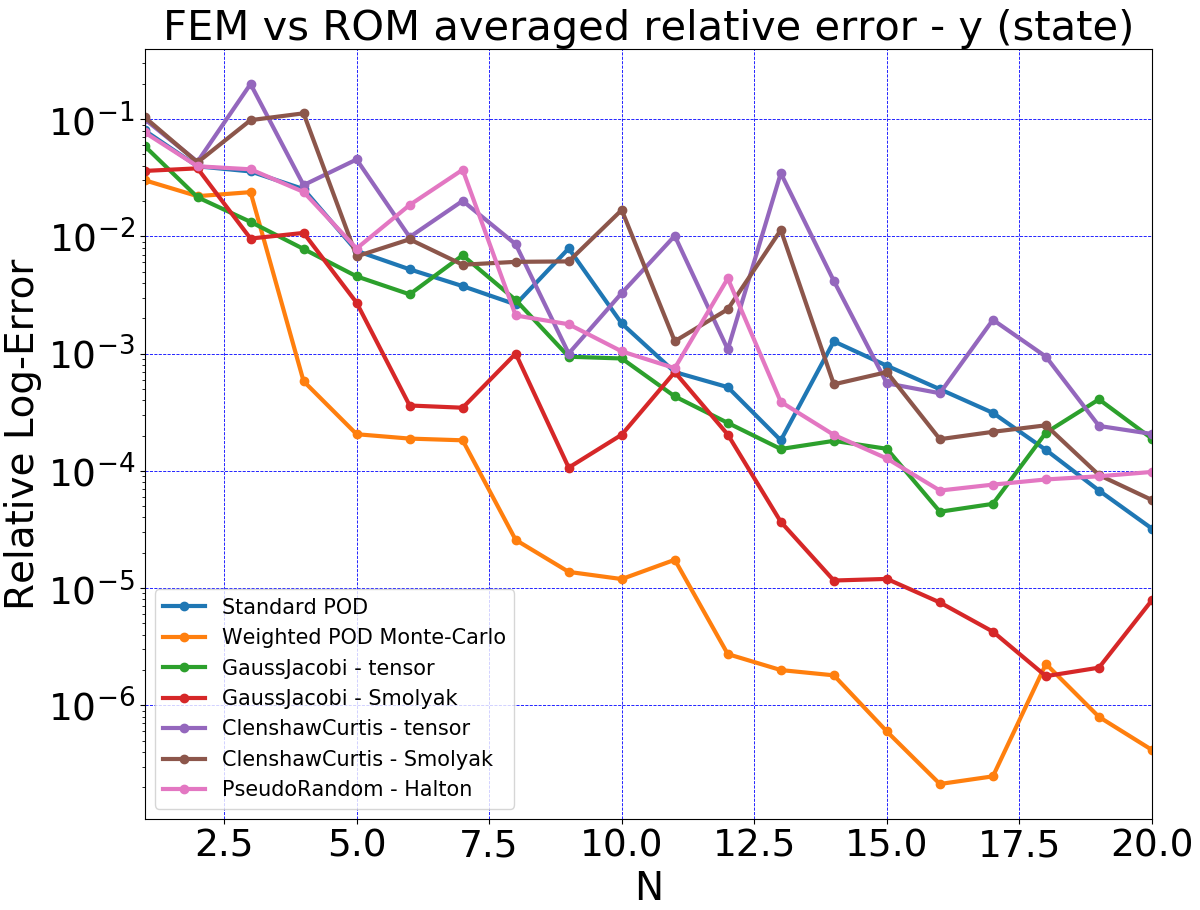}
    \includegraphics[scale=0.12]{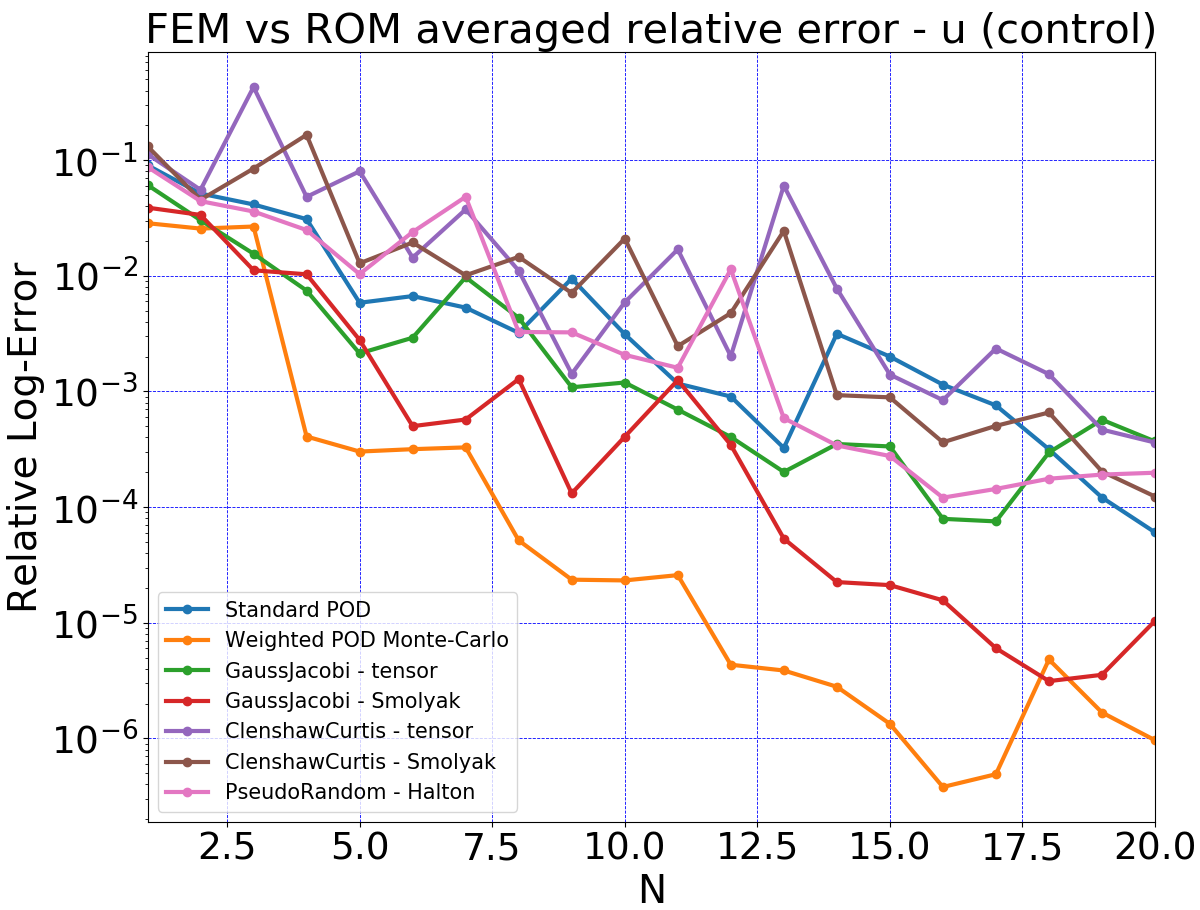}
    \includegraphics[scale=0.12]{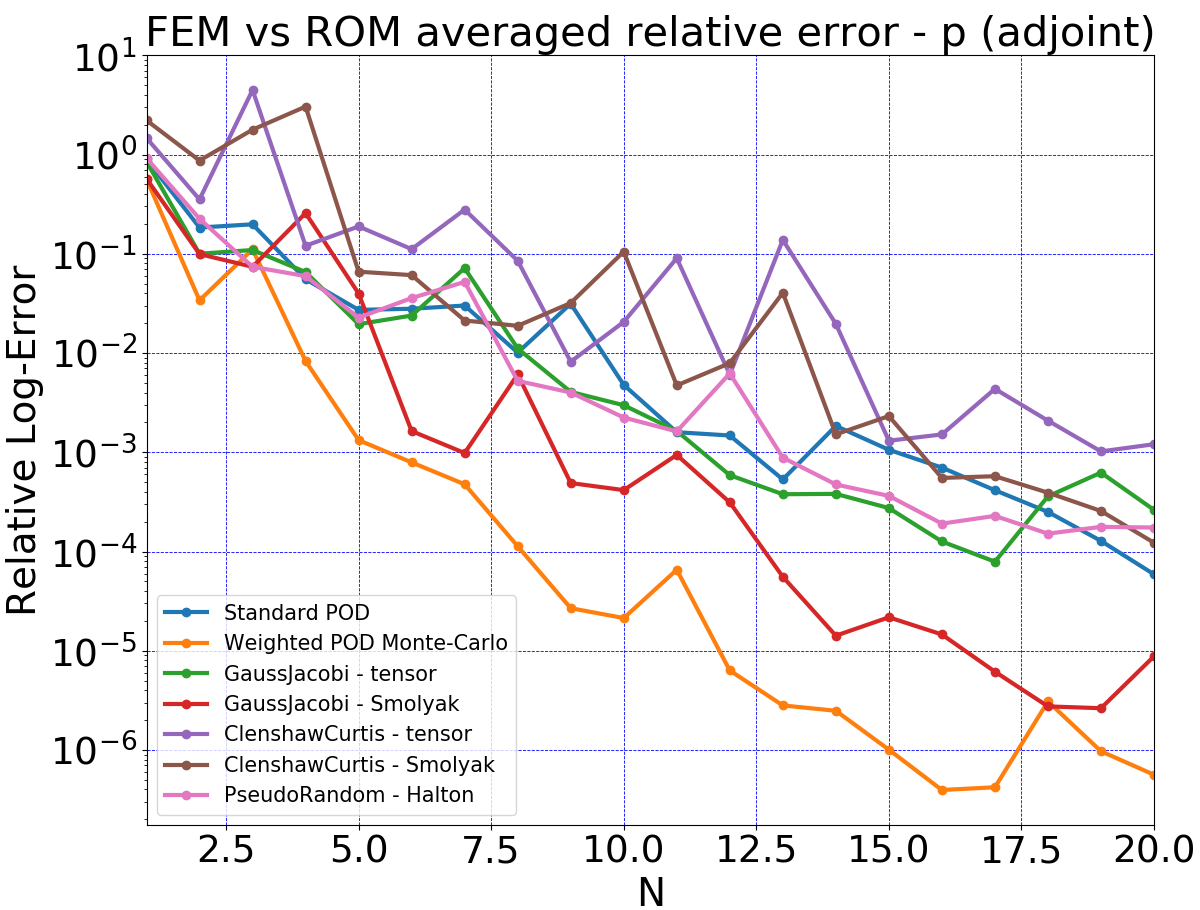}
    \caption{Relative Errors for the Graetz-Poiseuille Problem with $\boldsymbol{\mu}$ following distribution \eqref{beta-graetz} on the parameter space $\mathcal{P}$ - \textit{Offline-Online} Stabilization; State (\underline{left}), Control (\underline{center}), Adjoint (\underline{right}); Standard POD (blue), wPOD Monte-Carlo (orange), Gauss-Jacobi tensor rule (green), Gauss-Jacobi Smolyak grid (red), Clenshaw-Curtis tensor rule (cyan), Clenshaw-Curtis Smolyak grid (dark green), Pseudo-Random based on Halton numbers (pink).}
    \label{fig:plot_graetz_onoffstab}
\end{figure}

\sloppy
In order to have a complete view of the ROM procedure, in Figure \ref{fig:graetz-grids} we first illustrate the grid points for the quadrature procedures explained in Section \ref{wpod}, so that we can see which points are taken during the sampling procedure before being weighted according to the prescribed rule. One can notice that in the tensor rules and in the Smolyak technique with Clenshaw-Curtis procedure several points lie close to the boundary. Therefore, it is expected that those snapshots are pretty much different from the ones of a distribution that is more concentrated inside the domain, for instance, as \eqref{beta-graetz}. In this way, we have a poor representation of the behaviors that are more likely to happen.
\sloppy
Singular value decay of the snapshots matrices is shown in Figure \ref{fig:plot_Graetz Geom-sing values} and later the projection error onto the POD space in Figure \ref{fig:plot_Graetz Geom-proj error} so as to compare it with the ones between ROM and FEM solutions in Figures \ref{fig:plot_Graetz Geom-not stab} and \ref{fig:plot_graetz_onoffstab}, respectively. In Figure \ref{fig:plot_Graetz Geom-sing values} we see that generally, all the procedures show a decreasing tendency over $N$; however, the weighted Monte-Carlo strategy overcomes all the other procedures by at least a factor of $10^4$. This means that weighted Monte-Carlo POD needs less reduced basis to reach good accuracy levels for all three variables. This phenomenon also reflects in the projection errors (and in the Offline-Online ones, too), where the distance between the FEM space and the ROM one is always definitely smaller using the weighted Monte-Carlo POD. This trend of the projection errors indicates that the performance differences for the relative errors shown later will be due to the different approximation quality of weighted POD spaces and not provoked by some effect on the ROMs built on top of these spaces. \\
We remark that a different number of bases for the three variables can be chosen using truncation error tolerances, but for the sake of simplicity (and following the code implementation of the RBniCS library \cite{RBniCS}) we define a priori value of $N$ for all the variables.
\sloppy
\sloppy

\begin{figure}
        \centering
        \includegraphics[scale=0.2]{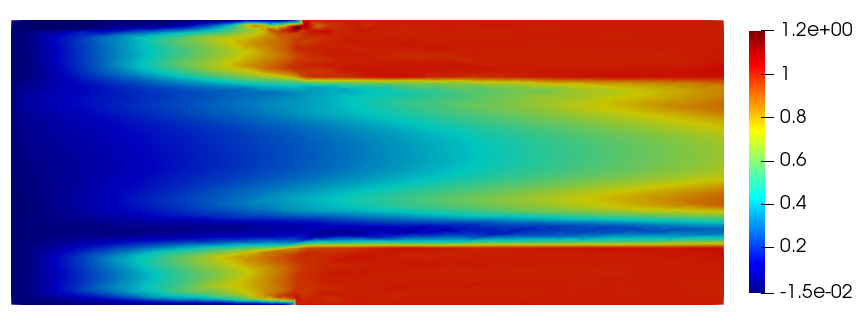}
        \includegraphics[scale=0.2]{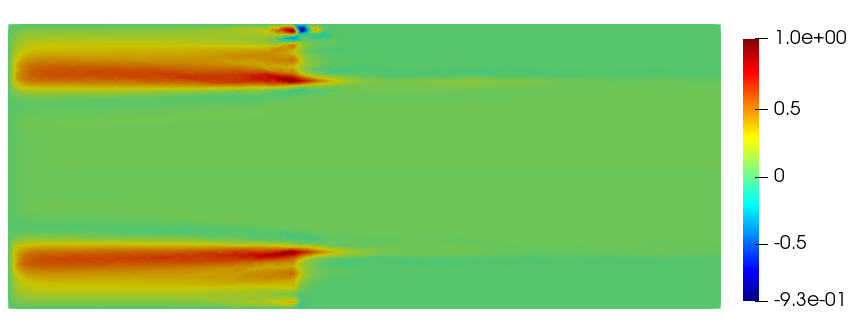} \\
        \includegraphics[scale=0.2]{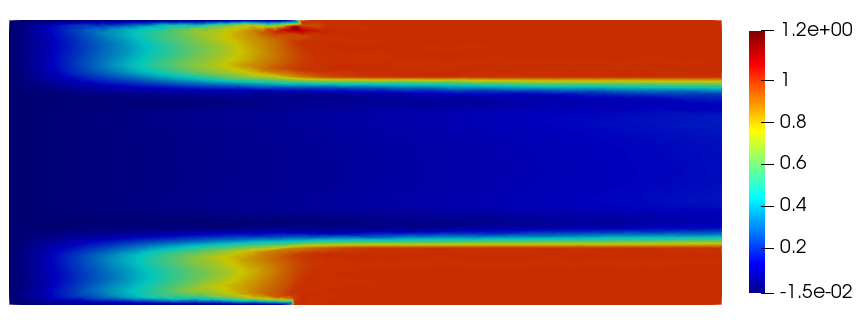}
        \includegraphics[scale=0.2]{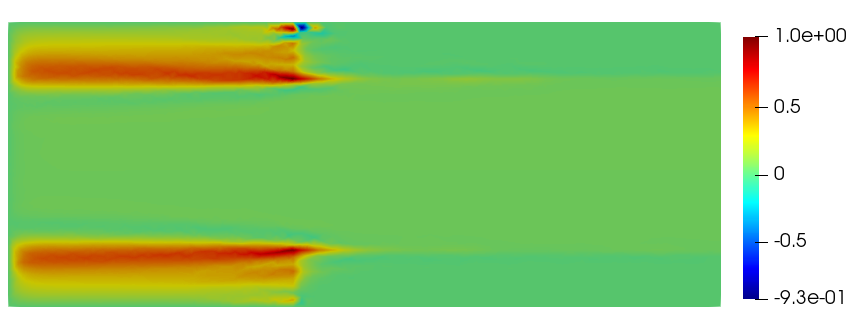} 
        \caption{(\underline{top}) wPOD Monte-Carlo \textit{Offline-Only} stabilized and (\underline{bottom})  \textit{Offline-Online} stabilized solution, $y$ (\underline{right}) and $u$ (\underline{left}), $\boldsymbol{\mu} = (10^{5},1.5)$,  $h=0.034$, $\alpha=0.01$, $N_{\text{train}}=100$, $\delta_K=1.0$, $N=20$.}
        \label{fig:graetz_onoff_stab}
\end{figure}

\sloppy
We firstly exploit the \textit{Offline-Only} stabilization procedure, which results regarding errors are shown in Figure \ref{fig:plot_Graetz Geom-not stab}. The performance is not good for any kind of wPOD. Moreover, the Standard POD does not perform good, either. Relative errors never drop under $10^{-2}$ for any variables. {All of them do not follow the same trends as the corresponding projection errors in Figure \ref{fig:plot_Graetz Geom-proj error}, which means that the Offline-Only Galerkin projector is not closed to the orthogonal one between the ROM and the FEM spaces.} This implies that more stabilization is necessary in this case.
\sloppy
\sloppy
In Figure \ref{fig:plot_graetz_onoffstab}, relative errors of the Offline-Online stabilization procedure are presented. Here the trend seems better than the \textit{Offline-Only} one, because these quantities decay faster along the value of $N$, {following a similar behavior with respect to the projection errors in Figure \ref{fig:plot_Graetz Geom-proj error}, as one would have expected from previous considerations.} The wPOD Monte-Carlo is the best performer for all $y,u,p$ variables, as a matter of fact, it reaches $e_{y, 16} = 2.13 \cdot 10^{-7}$ for the state, for the adjoint $e_{p, 16}= 3.95 \cdot 10^{-7}$ and the control $e_{u, 16}=3.80 \cdot 10^{-7}$. This procedure has a better performance of the Standard POD, which its accuracy is at least $100$ times inferior to the wPOD Monte-Carlo after $N>11$. Concerning other rules, it can be noticed that Smolyak grid techniques perform better than their tensor-rule counterparts, despite having a training set whose cardinality is similar, but less of $100$: $93$ and $91$ for the Clenshaw-Curtis and Gauss-Jacobi sparse grids, respectively.
\sloppy
\sloppy
In Figure \ref{fig:graetz_onoff_stab}, we visually compare the two possibilities of stabilization for the geometrical parametrization of the Graetz-Pouiseuille problem for the wPOD Monte-Carlo.
\sloppy
\sloppy
In Table \ref{speed-g}, we compare the speedup-index for all wPOD algorithms. We see that computational values are all of the same order of magnitude. For the wPOD Monte-Carlo we calculate $87$ reduced solutions in the time of a FEM one.
\sloppy
\begin{table}
\centering

\begin{tabular}{|c|c|c|c|c|c|c|c|}
\hline \multicolumn{8}{|c|}{Speedup-index Graetz-Poiseuille Problem: \textit{Offline-Online} Stabilization - $\mu_1, \mu_2 \sim $ Beta(5,3)}  \\
\hline $N$ & POD & wPOD & Gauss tensor & Gauss Smolyak & CC tensor & CC Smolyak & Ps. Random \\
\hline $4$ & $113.0$ & $108.9$ & $110.1$  & $110.1$ & $106.5$ & $109.4$ & $112.0$ \\
\hline $8$ & $108.4$ & $105.1$ & $104.9$ & $106.1$ & $102.1$ & $105.9$ & $107.4$\\
\hline $12$ & $103.3$ & $100.2$ & $99.9$ & $99.8$ & $99.1$ & $96.9$ & $101.7$\\
\hline $16$ & $97.2$ & $92.5$ & $95.1$ & $94.5$ & $92.6$ & $94.2$ & $96.9$\\
\hline $20$ & $90.5$ & $87.3$ & $87.0$ & $88.0$ & $85.8$ & $86.3$ & $89.7$\\
\hline
\end{tabular}
  \caption{Average Speedup-index of \textit{Offline-Online} Stabilization for the Graetz-Poiseuille Problem under geometrical parametrization. From left to right: Standard POD, wPOD Monte-Carlo, Gauss-Jacobi tensor, Gauss-Jacobi Smolyak grid, Clenshaw-Curtis tensor, Clenshaw-Curtis Smolyak grid, Pseudo-Random based on Halton numbers.}
  \label{speed-g}
\end{table}

\begin{figure}
    \centering
    \includegraphics[scale=0.15]{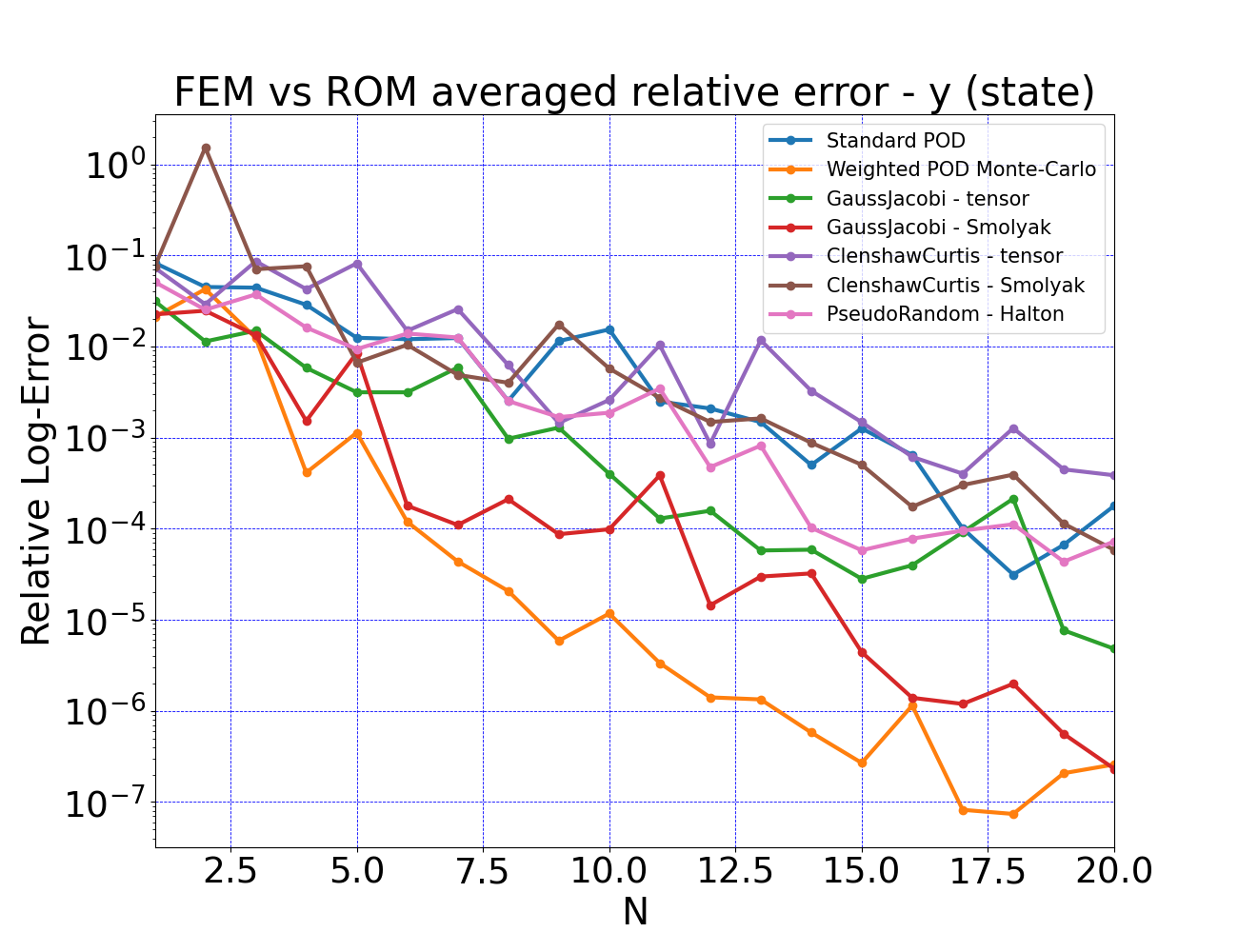}
    \includegraphics[scale=0.15]{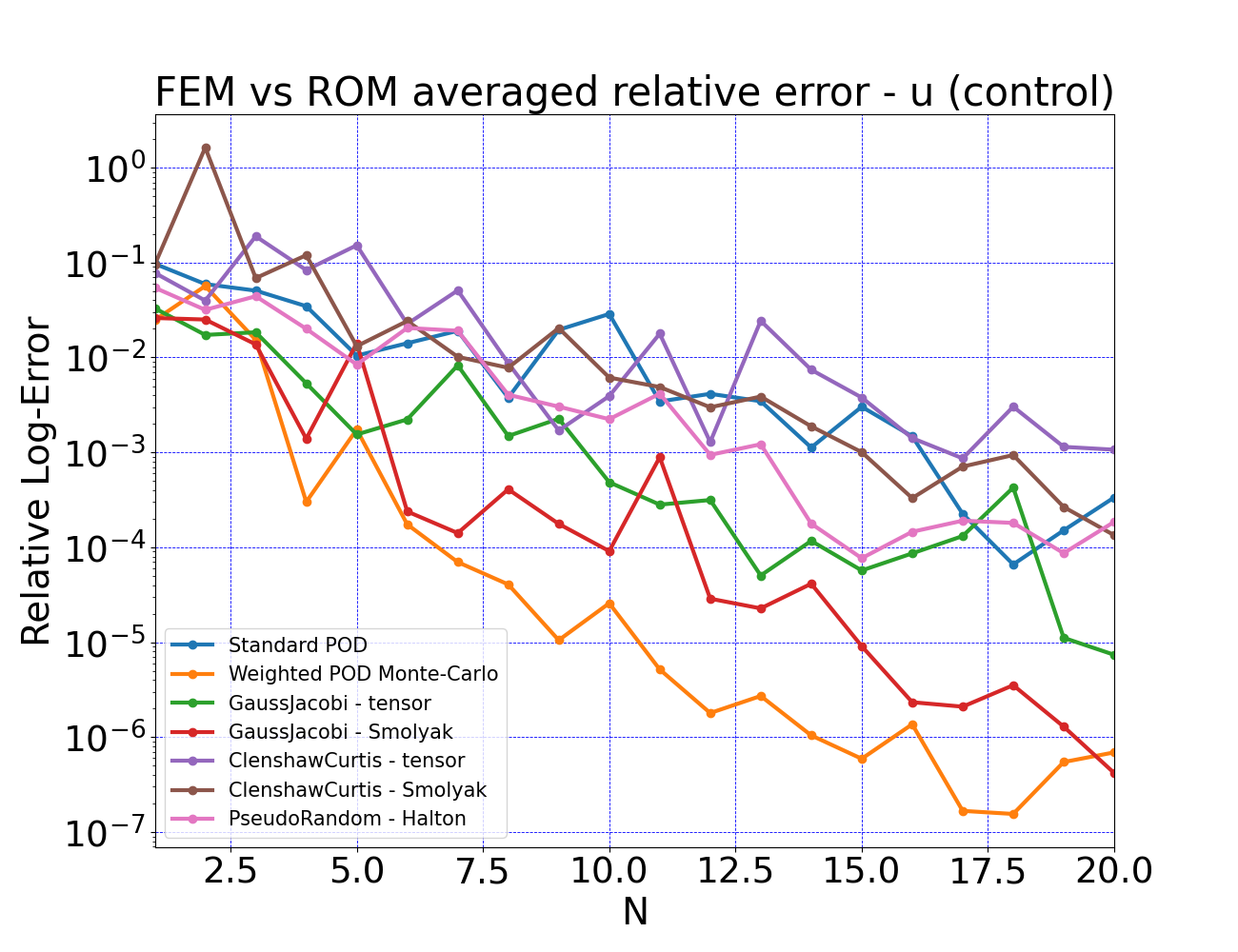}
    \includegraphics[scale=0.15]{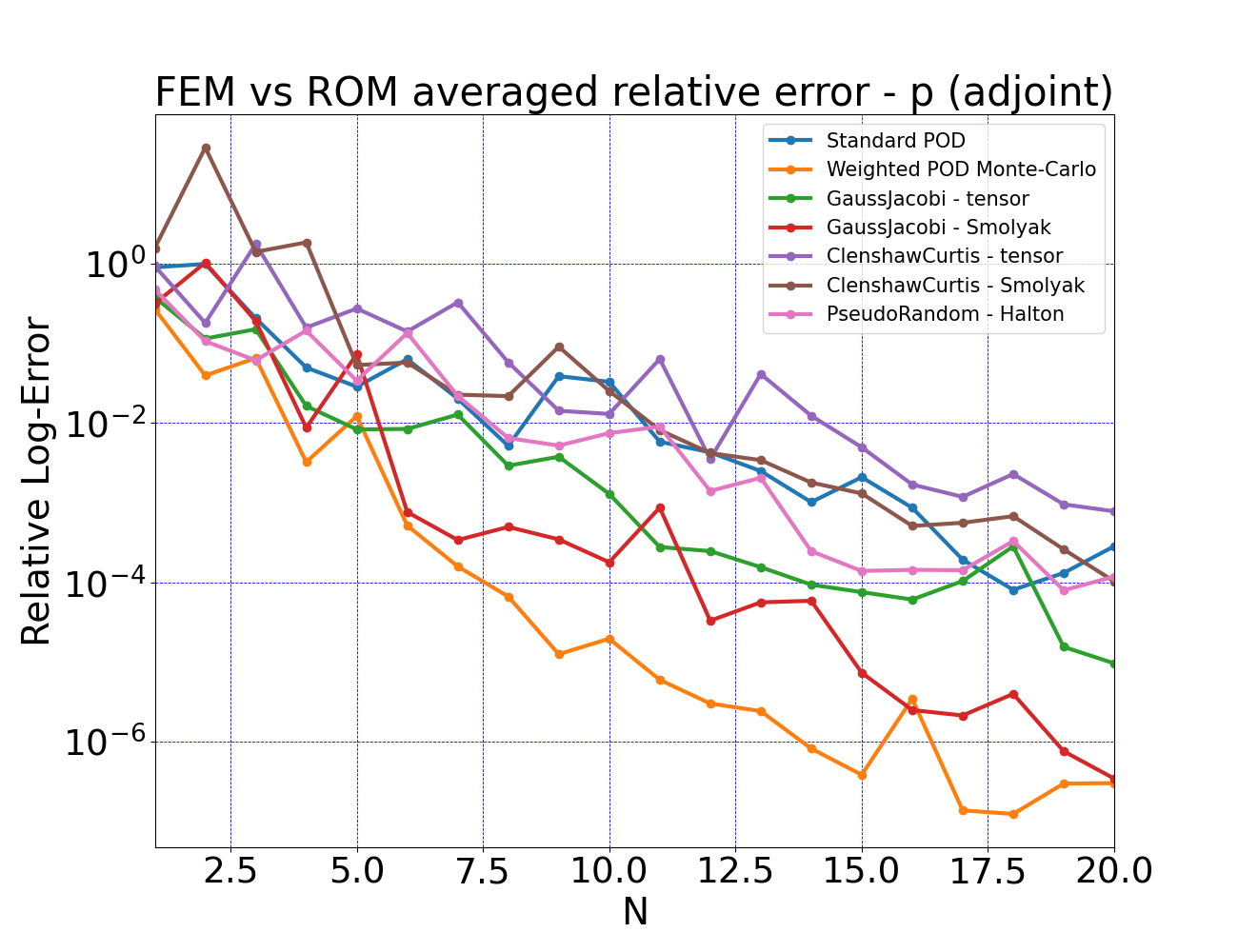}
    \caption{Relative Errors for the Graetz-Poiseuille Problem with $\boldsymbol{\mu}$ following distribution \eqref{beta-graetz-2} on the parameter space $\mathcal{P}$- \textit{Offline-Online} Stabilization; State (\underline{left}), Control (\underline{center}), Adjoint (\underline{right}); Standard POD (blue), wPOD Monte-Carlo (orange), Gauss-Jacobi tensor rule (green), Gauss-Jacobi Smolyak grid (red), Clenshaw-Curtis tensor rule (cyan), Clenshaw-Curtis Smolyak grid (dark green), Pseudo-Random based on Halton numbers (pink).}
    \label{fig:plot_graetz_onoffstab-2}
\end{figure}

\begin{figure}
    \centering
    \includegraphics[scale=0.15]{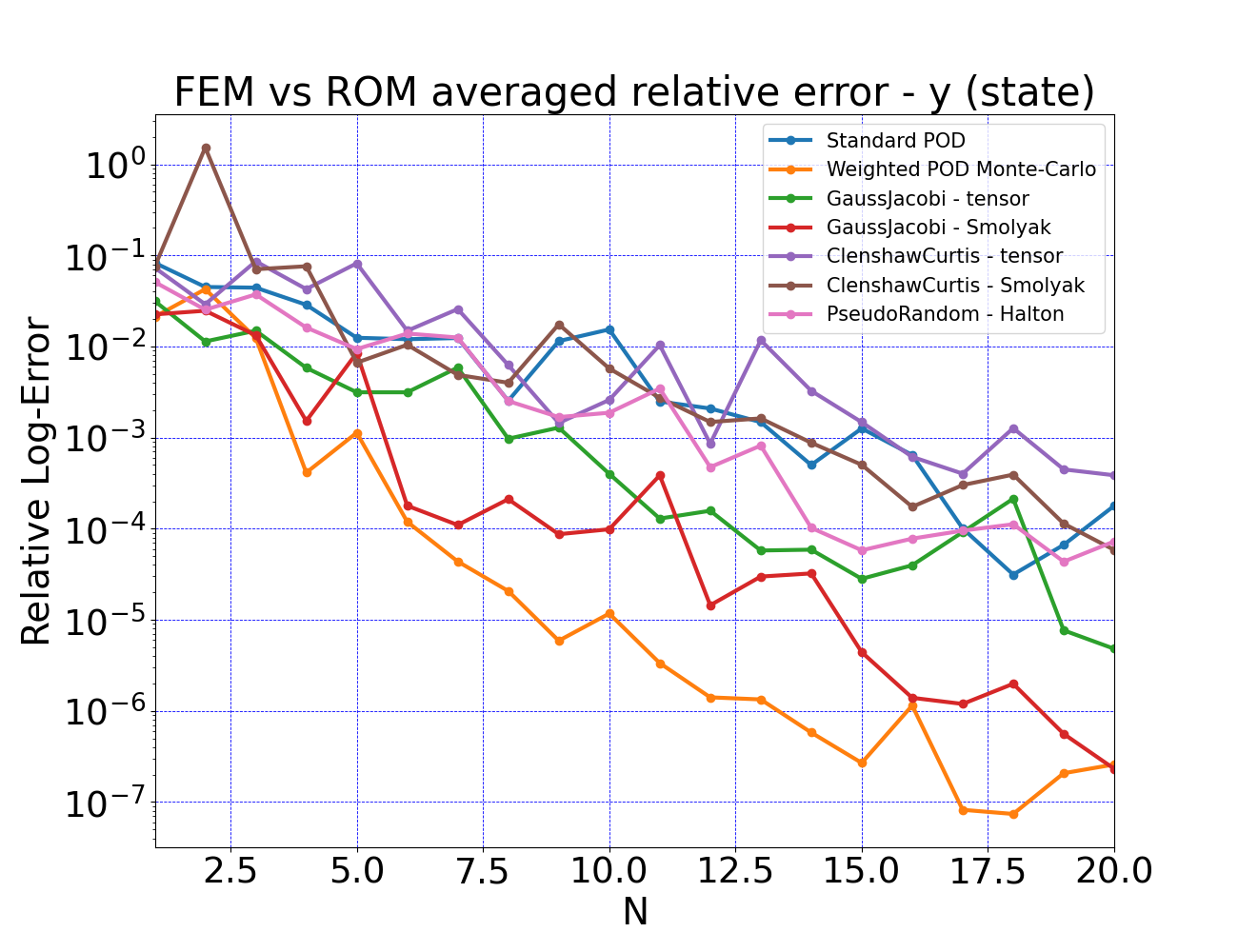}
    \includegraphics[scale=0.15]{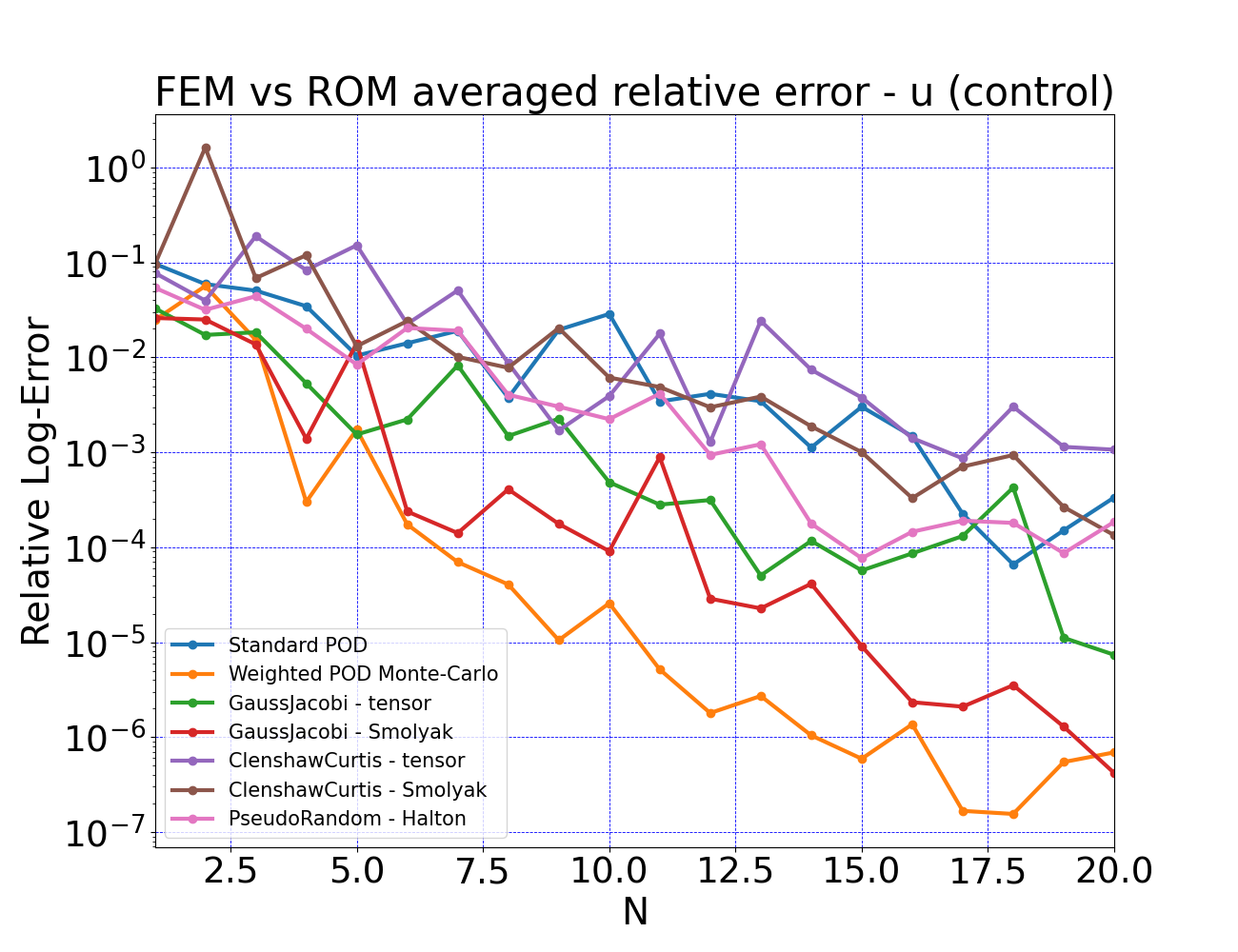}
    \includegraphics[scale=0.15]{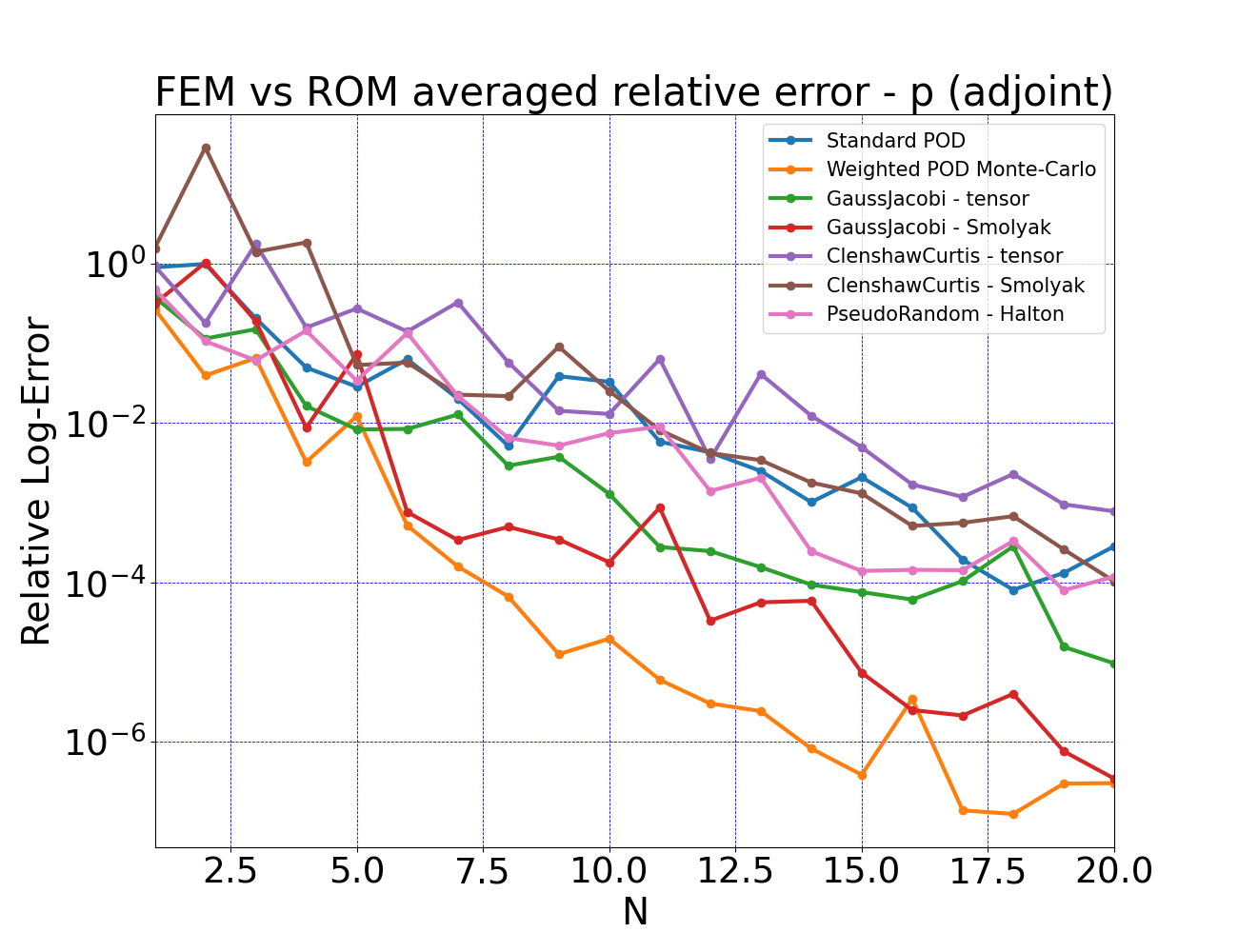}
    \caption{Relative Errors for the Graetz-Poiseuille Problem with $\boldsymbol{\mu}$ following distribution \eqref{beta-graetz-3} on the parameter space $\mathcal{P}$- \textit{Offline-Online} Stabilization; State (\underline{left}), Control (\underline{center}), Adjoint (\underline{right}); Standard POD (blue), wPOD Monte-Carlo (orange), Gauss-Jacobi tensor rule (green), Gauss-Jacobi Smolyak grid (red), Clenshaw-Curtis tensor rule (cyan), Clenshaw-Curtis Smolyak grid (dark green), Pseudo-Random based on Halton numbers (pink).}
    \label{fig:plot_graetz_onoffstab-3}
\end{figure}

We ran other experiments concerning two more distributions regarding the Graetz-Poiseuille Problem with the same parameter space $\mathcal{P}$. For the sake of brevity, we only show the error for the Offline-Online Stabilization procedure. The first distribution that we study is the following one
\begin{equation}\label{beta-graetz-2}
\begin{aligned}
    {\mu_1} \sim 1 + \big(10^5 - 1\big) X_1, \text{ where } X_1 \sim \text{Beta}(10,10), \\
     \mu_2 \sim 0.5 + \big( 1.5 - 0.5\big) X_2, \text{ where } X_2 \sim \text{Beta}(10,10).
\end{aligned}
\end{equation}  
In Figure \ref{fig:plot_graetz_onoffstab-2}, we see again that the Weighted Monte-Carlo method performs very well in terms of accuracy
Finally, we simulate the same problem with the following distribution of the parameter on $\mathcal{P}$
\begin{equation}\label{beta-graetz-3}
\begin{aligned}
    {\mu_1} \sim 1 + \big(10^5 - 1\big) X_1, \text{ where } X_1 \sim \text{Beta}(20,20), \\
     \mu_2 \sim 0.5 + \big( 1.5 - 0.5\big) X_2, \text{ where } X_2 \sim \text{Beta}(20,20),
\end{aligned}
\end{equation}  
where we concentrate even more the information in the middle of the parameter space. As expected, all quadrature rules that have several samples close to the boundary do not perform as well as the Weighted Monte-Carlo. We can see the trend of the errors in Figure \ref{fig:plot_graetz_onoffstab-3}.

\sloppy
Now we want to present the parabolic version of Problem \eqref{graetz-system}. This unsteady problem has been studied without optimization in \cite{pacciarini2014stabilized, torlo2018stabilized} in a deterministic context and in \cite{torlo2018stabilized} in a UQ one. Instead, the deterministic OCP($\boldsymbol{\mu}$) Graetz Problem under boundary control without Advection-dominancy is studied in \cite{strazzullo2020pod,strazzullo2021model} and the deterministic distributed control configuration is analyzed in \cite{zoccolan2024streamline}.
\sloppy
Recalling Figure \ref{fig:Geometry-Graetz}, for a fixed $T>0$ the unsteady Graetz-Poiseuille Problem is posed as follows: find $(y, u) \in \tilde{\mathcal{Y}} \times \mathcal{U}$ which solves
\begin{equation*}
    \min\limits_{(y,u)} \frac{1}{2} \int\limits_{\Omega_{obs}(\boldsymbol{\mu})\times (0,T)}(y(\boldsymbol{\mu})- y_d)^2 \; d \Omega  +
\frac{\alpha}{2} \int\limits_{\Omega(\boldsymbol{\mu}) \times (0,T)}u(\boldsymbol{\mu})^2 \; d \Omega, \ \text{  such that }
\vspace{-2mm}
\end{equation*}
\begin{equation}
   \begin{cases}
      \displaystyle \partial_t y(\boldsymbol{\mu})-\frac{1}{{\mu_{1}}} \Delta y(\boldsymbol{\mu})+4 x_1(1-x_1) \partial_{x_0} y(\boldsymbol{\mu})=u(\boldsymbol{\mu}), & \text { in } \Omega(\boldsymbol{\mu})  \times (0,T), \\
y(\boldsymbol{\mu})=0, & \text { on } \Gamma_{1} \cup \Gamma_{5} \cup \Gamma_{6}  \times (0,T), \\
y(\boldsymbol{\mu})=1, & \text { on } \Gamma_{2}(\boldsymbol{\mu}) \cup \Gamma_{4}(\boldsymbol{\mu}) \times (0,T), \\
\displaystyle \frac{\partial y(\boldsymbol{\mu})}{\partial \nu}=0, & \text { on } \Gamma_{3}(\boldsymbol{\mu})  \times (0,T),
\\
\displaystyle y(\boldsymbol{\mu})(0)=y_0(x), & \text { in } \Omega(\boldsymbol{\mu}).
    \end{cases}
\end{equation}
As made in the steady version, we first consider a lifting procedure.
Simulations are run following the space-time setting proposed in Section \ref{time-dep-discr} for a prearranged number of time-steps $N_t$. 

\begin{figure}
        \centering
        \includegraphics[scale=0.165]{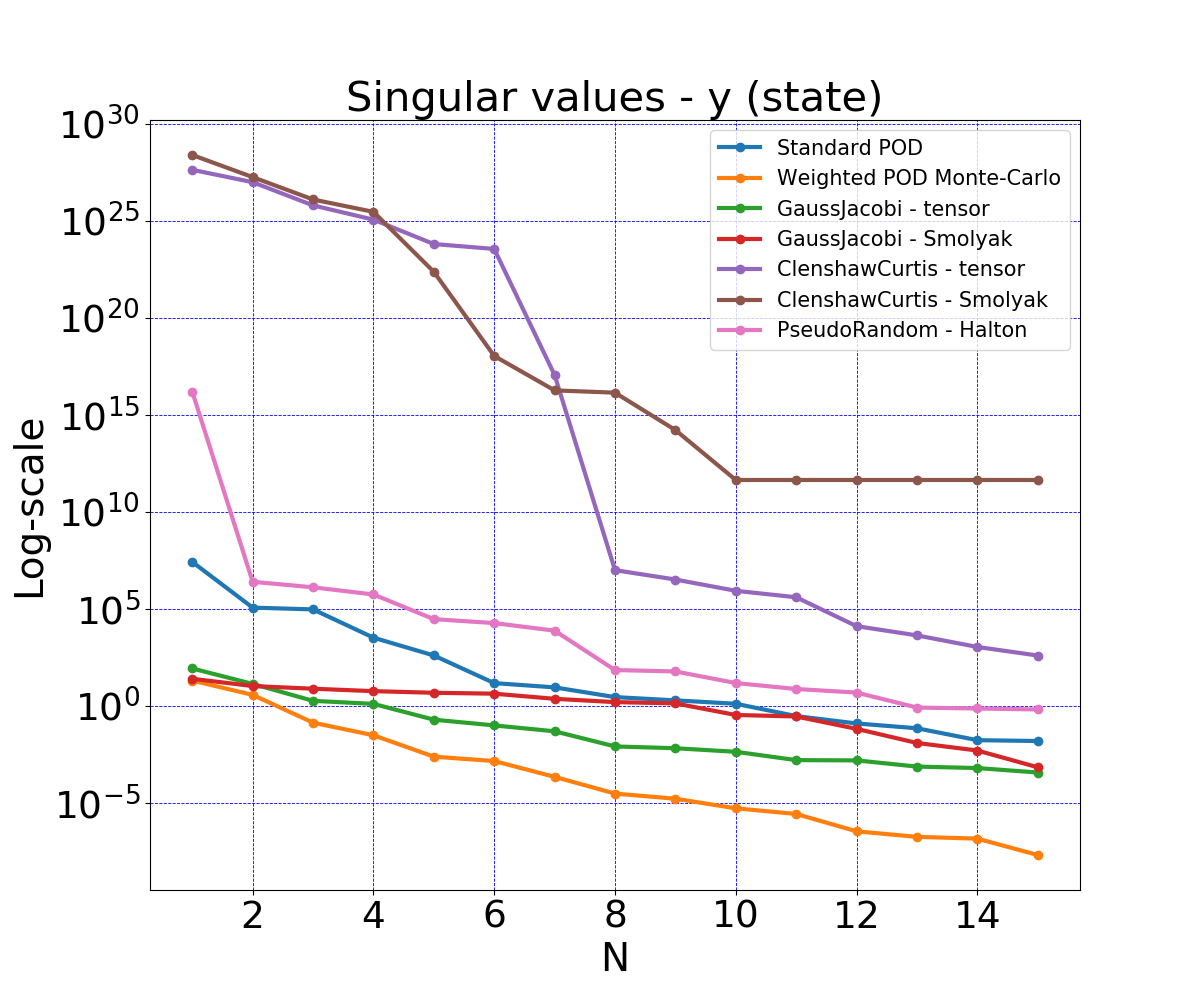} 
        \includegraphics[scale=0.165]{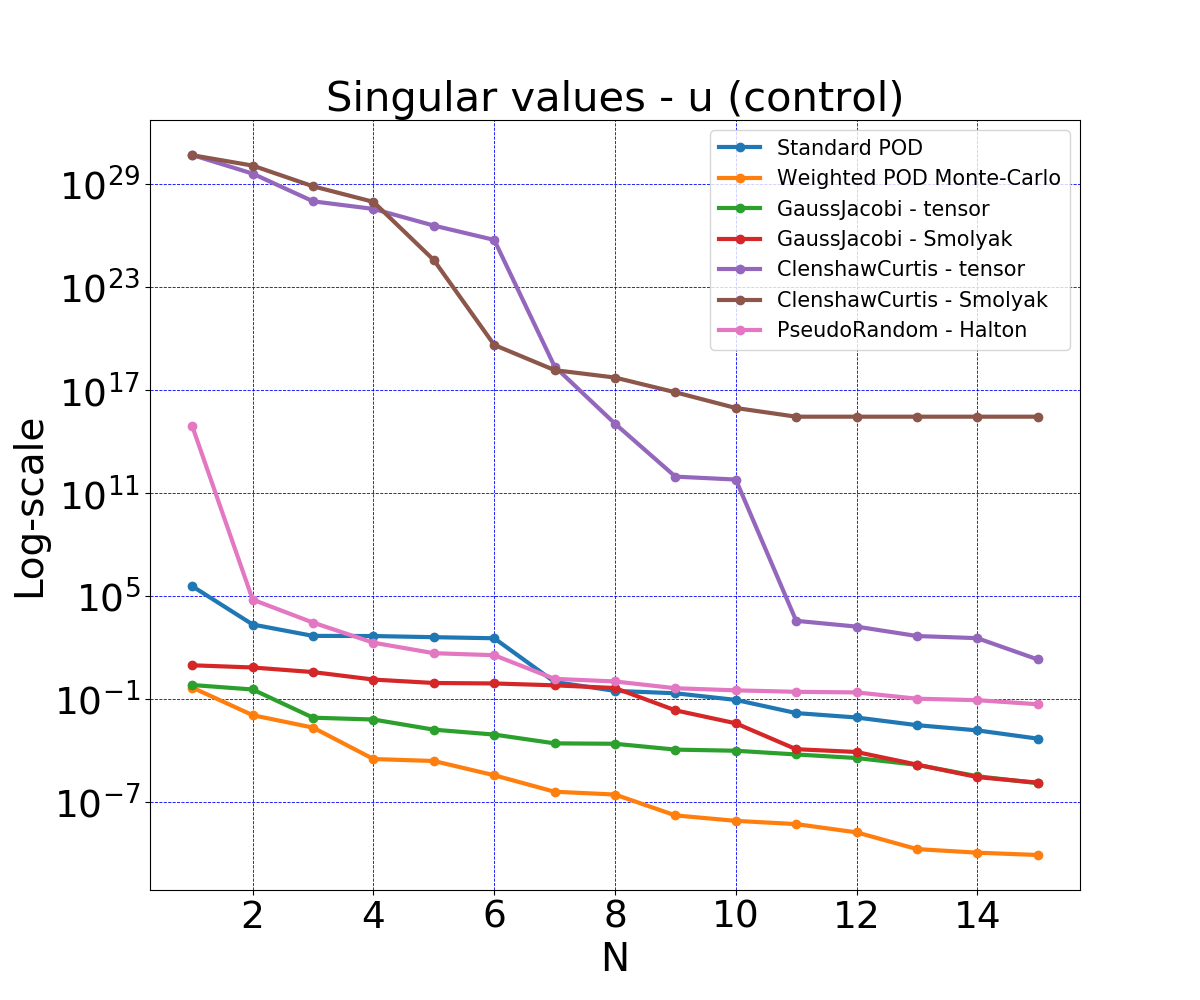} 
        \includegraphics[scale=0.165]{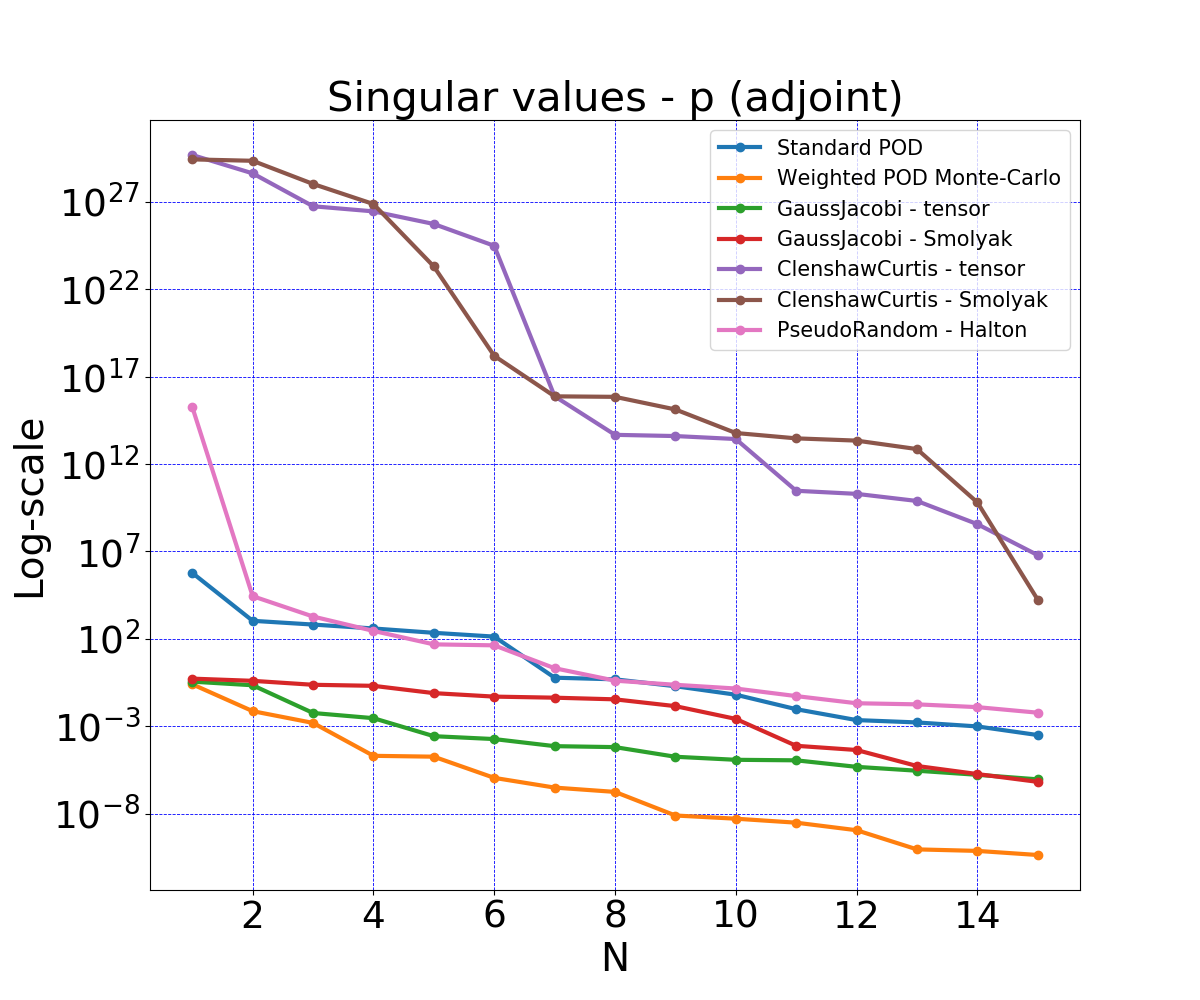} 
        \caption{Singular values decay for the snapshot matrices for the Parabolic Graetz-Poiseuille Problem; State (\underline{left}), Control (\underline{center}), Adjoint (\underline{right}); Standard POD (blue), wPOD Monte-Carlo (orange), Gauss-Jacobi tensor rule (green), Gauss-Jacobi Smolyak grid (red), Clenshaw-Curtis tensor rule (cyan), Clenshaw-Curtis Smolyak grid (dark green), Pseudo-Random based on Halton numbers (pink).}
        \label{fig:plot_Par_graetz-sing values}
\end{figure}

\begin{figure}
        \centering
        \includegraphics[scale=0.165]{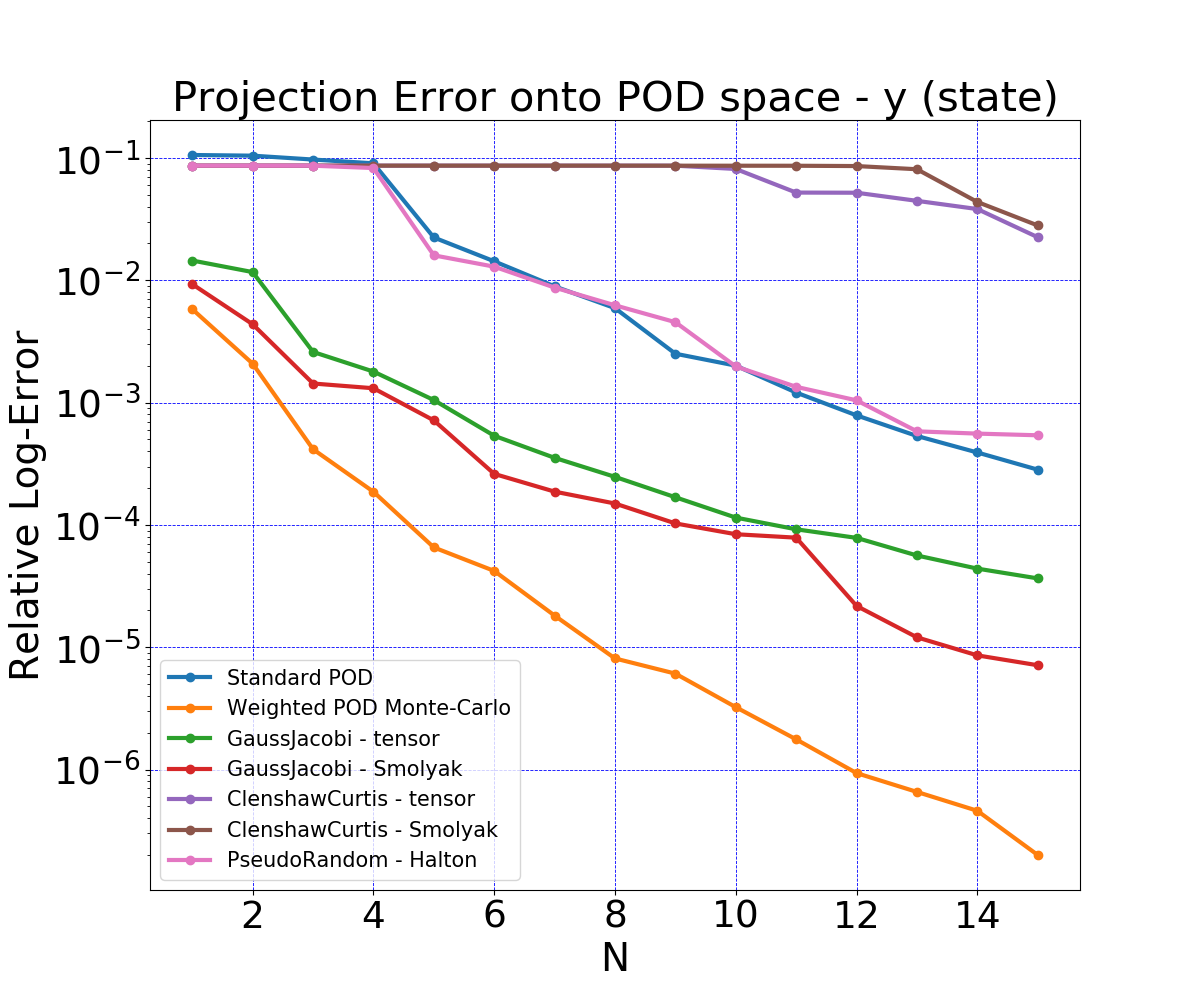} 
        \includegraphics[scale=0.165]{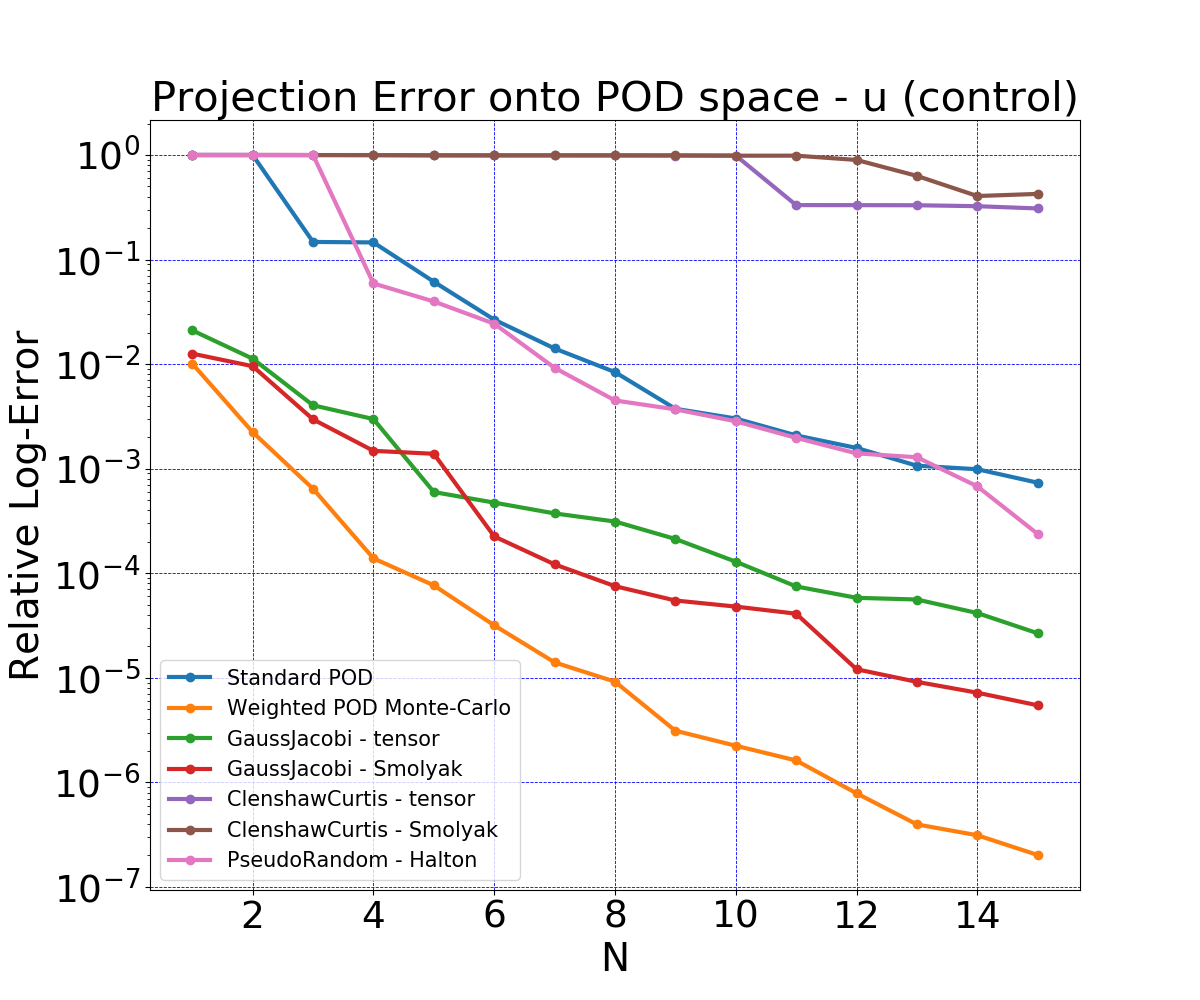} 
        \includegraphics[scale=0.165]{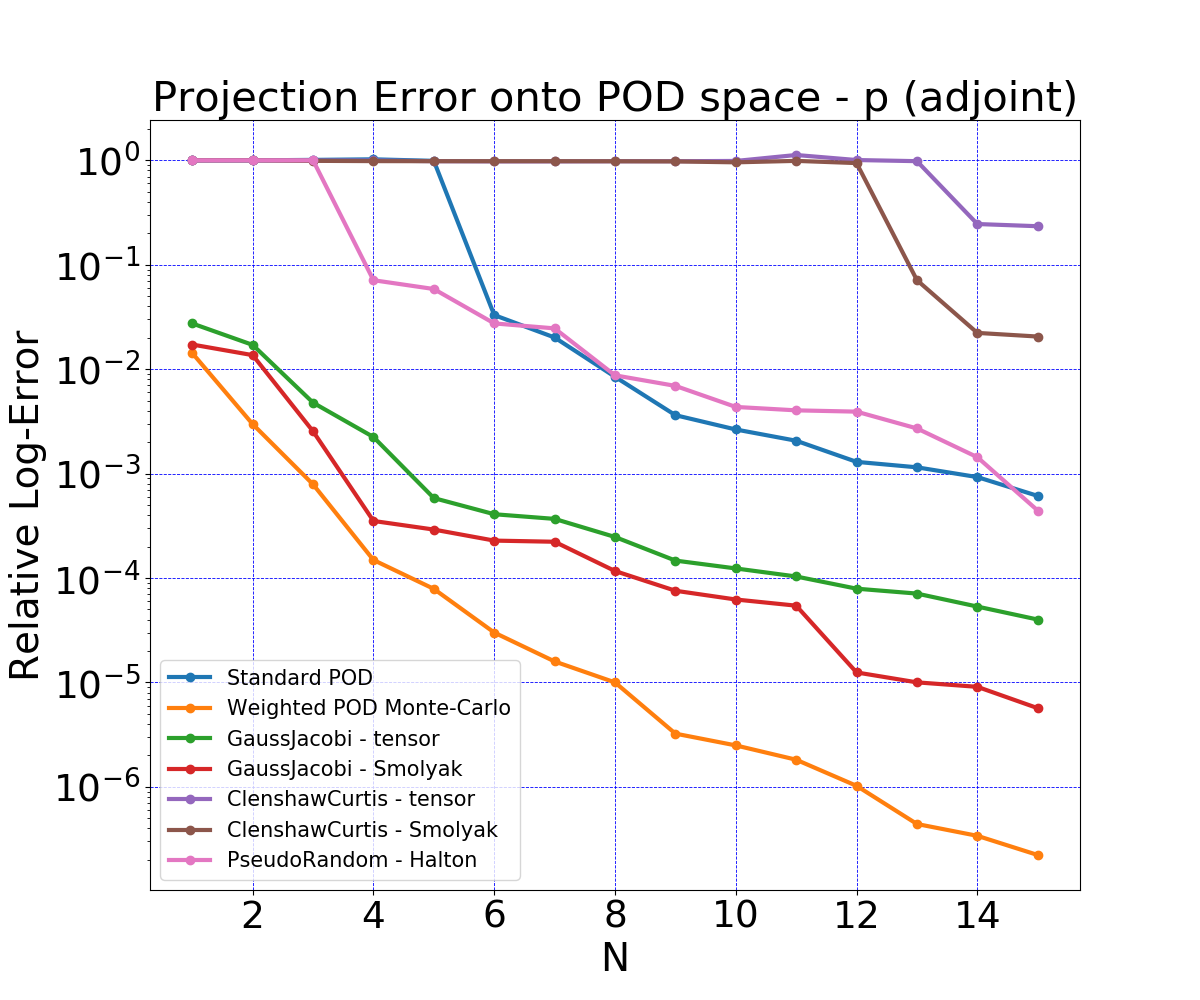} 
        \caption{Projection Errors onto the POD space for the Parabolic Graetz-Poiseuille Problem; State (\underline{left}), Control (\underline{center}), Adjoint (\underline{right}); Standard POD (blue), wPOD Monte-Carlo (orange), Gauss-Jacobi tensor rule (green), Gauss-Jacobi Smolyak grid (red), Clenshaw-Curtis tensor rule (cyan), Clenshaw-Curtis Smolyak grid (dark green), Pseudo-Random based on Halton numbers (pink).}
        \label{fig:plot_Par_graetz-proj error}
\end{figure}

\begin{figure}
    \centering
    \includegraphics[scale=0.127]{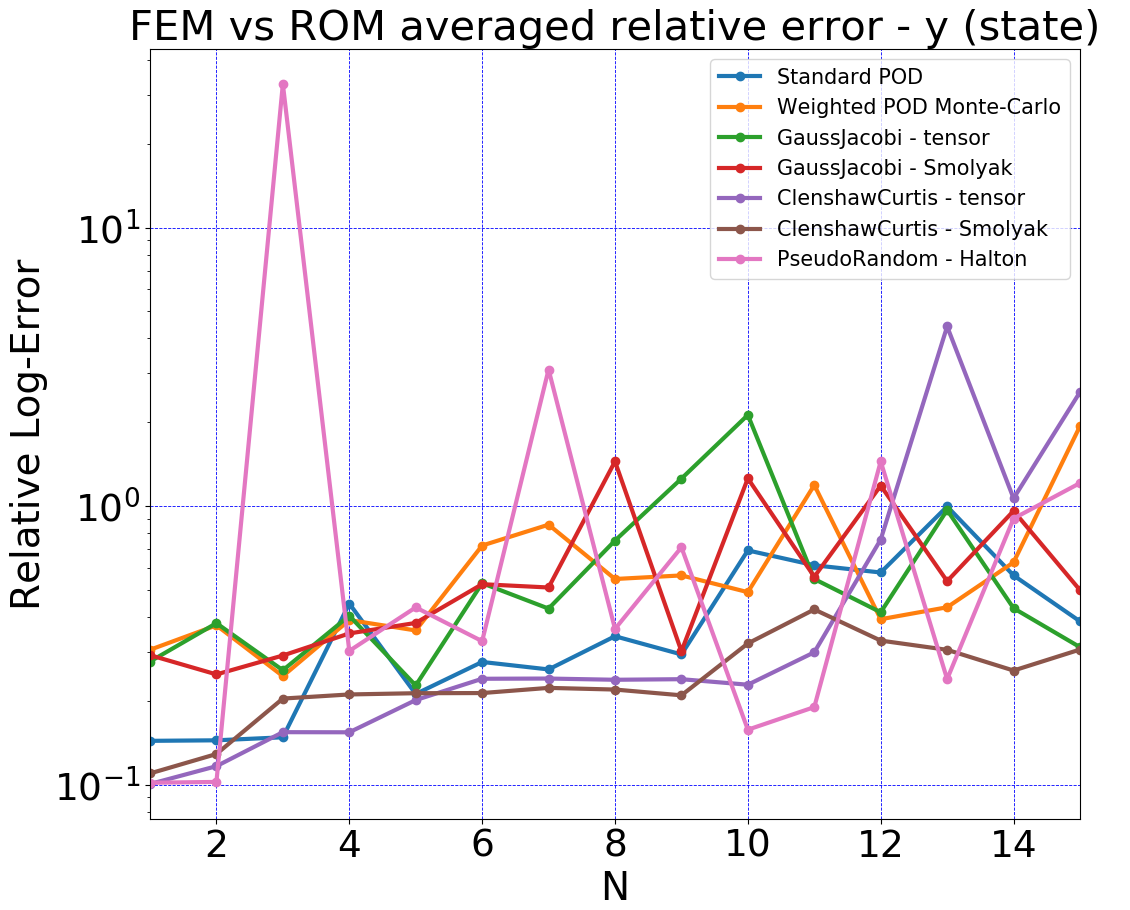}
    \includegraphics[scale=0.127]{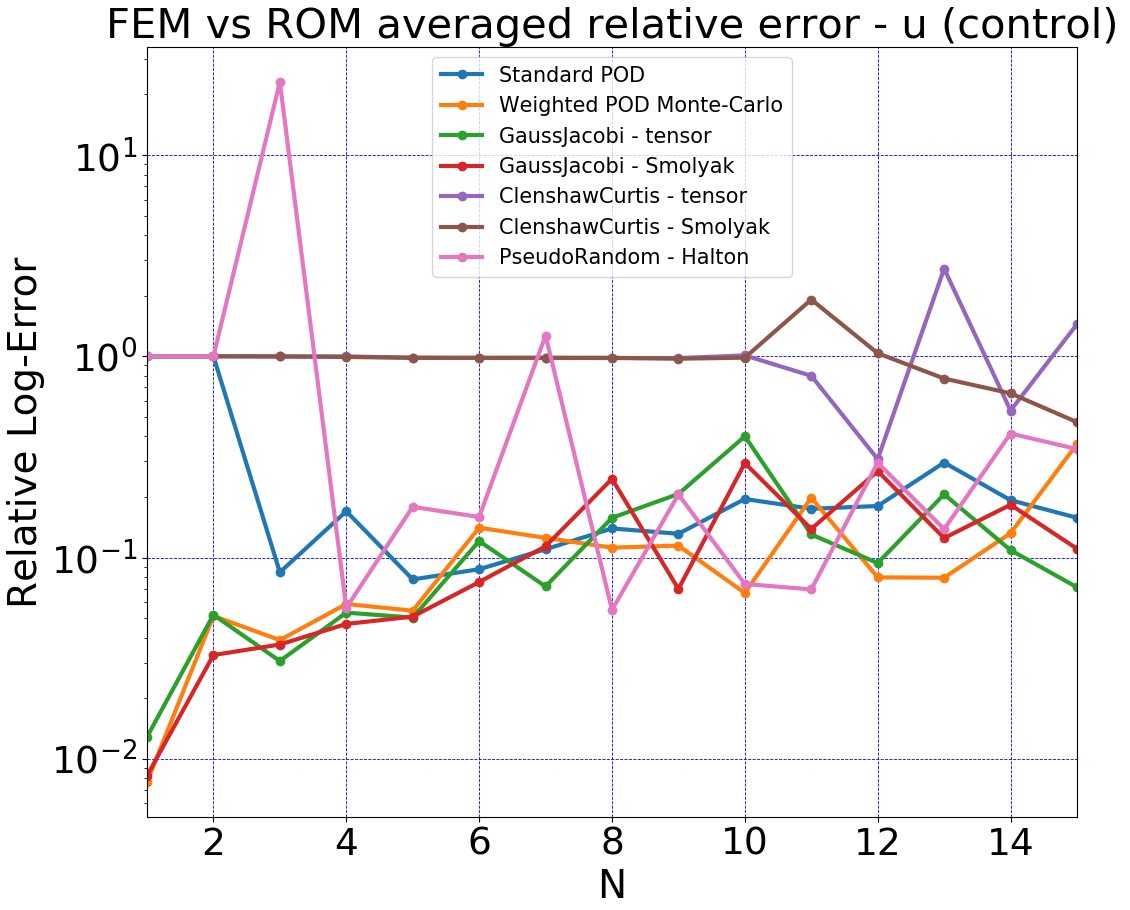}
    \includegraphics[scale=0.127]{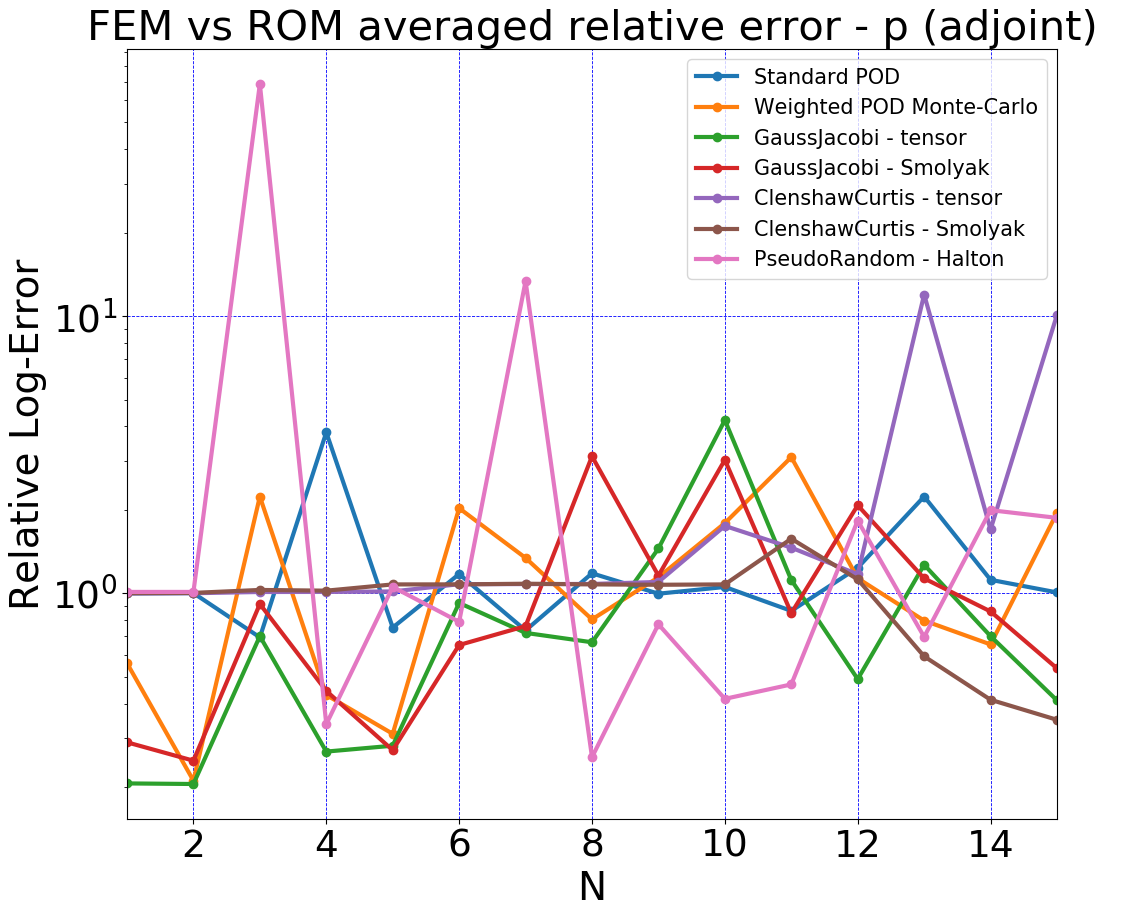}
    \caption{Relative Errors for the Parabolic Graetz-Poiseuille Problem - \textit{Offline-Only} Stabilization; State (\underline{left}), Control (\underline{center}), Adjoint (\underline{right}); Standard POD (blue), wPOD Monte-Carlo (orange), Gauss-Jacobi tensor rule (green), Gauss-Jacobi Smolyak grid (red), Clenshaw-Curtis tensor rule (cyan), Clenshaw-Curtis Smolyak grid (dark green), Pseudo-Random based on Halton numbers (pink).}
    \label{fig:plot_par_graetz_offstab}
\end{figure}

\begin{figure}
    \centering
    \includegraphics[scale=0.127]{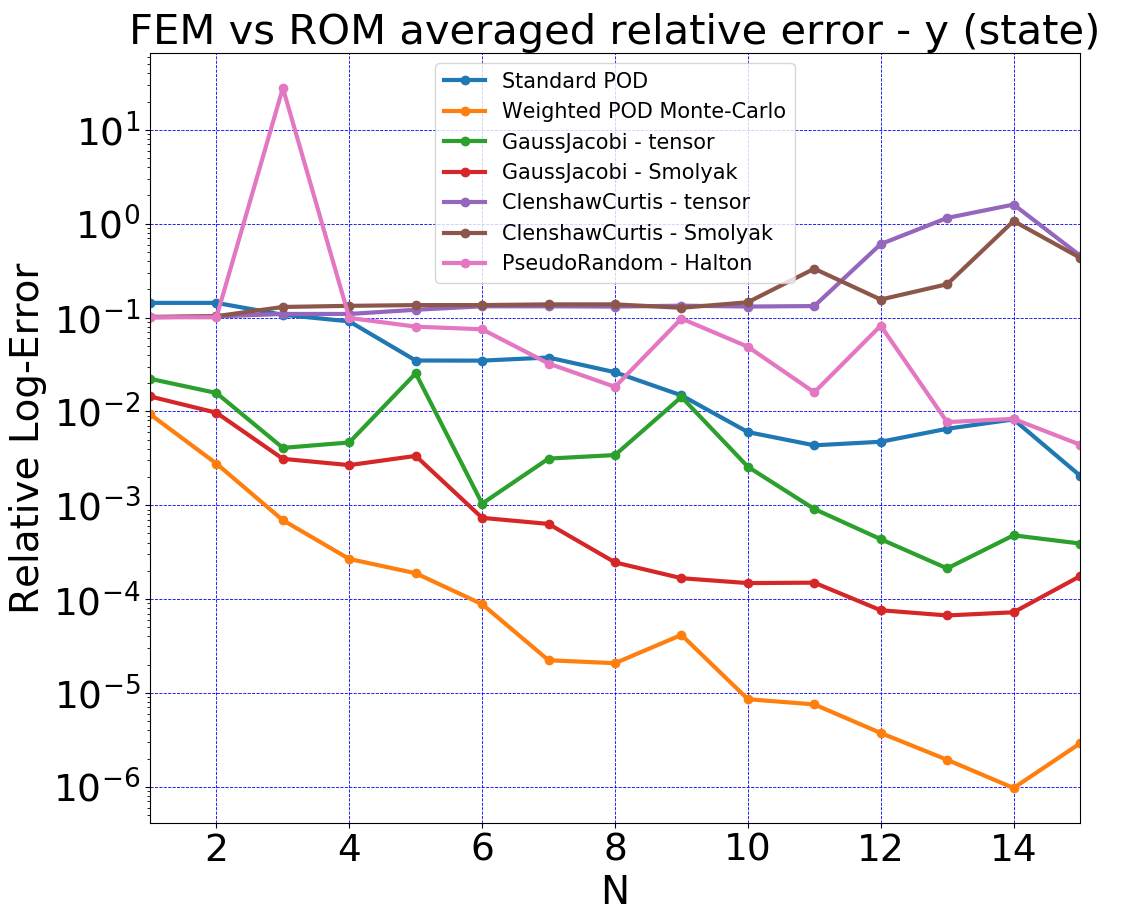}
    \includegraphics[scale=0.127]{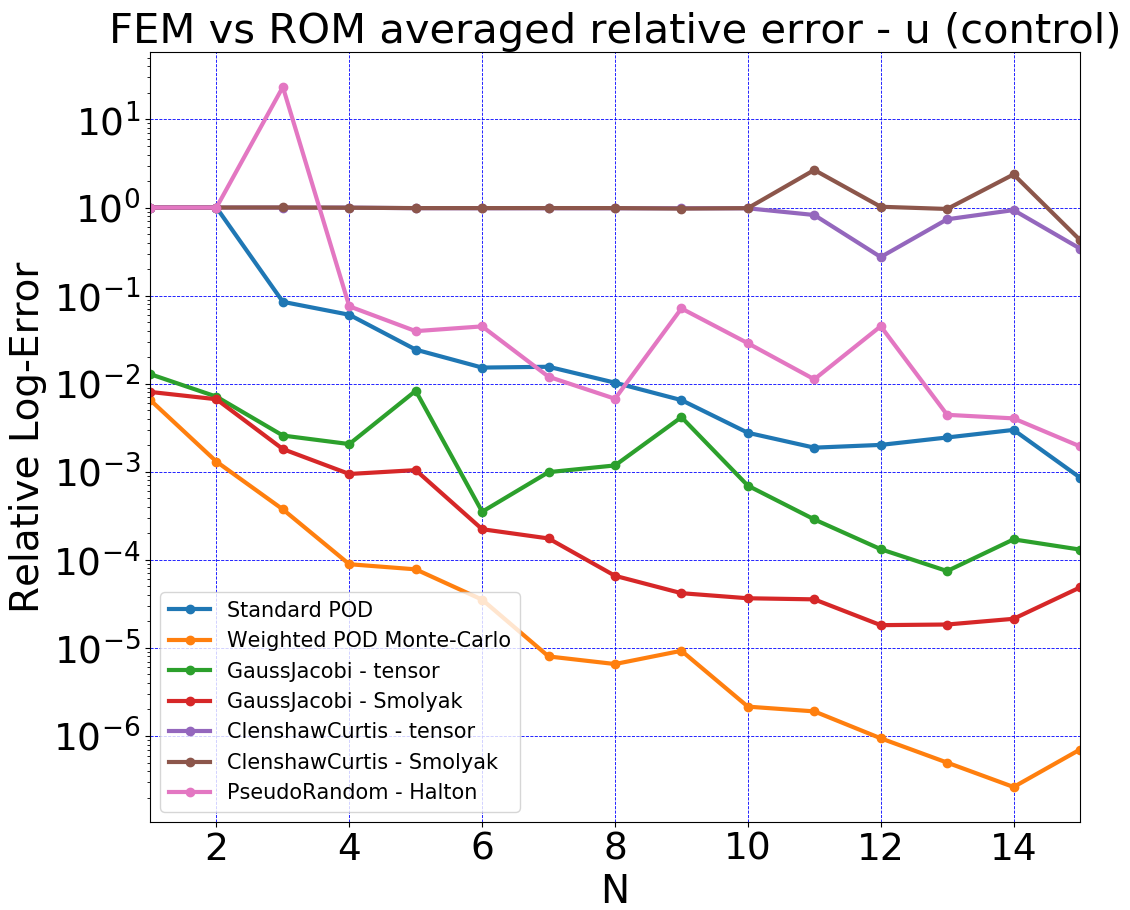}
    \includegraphics[scale=0.127]{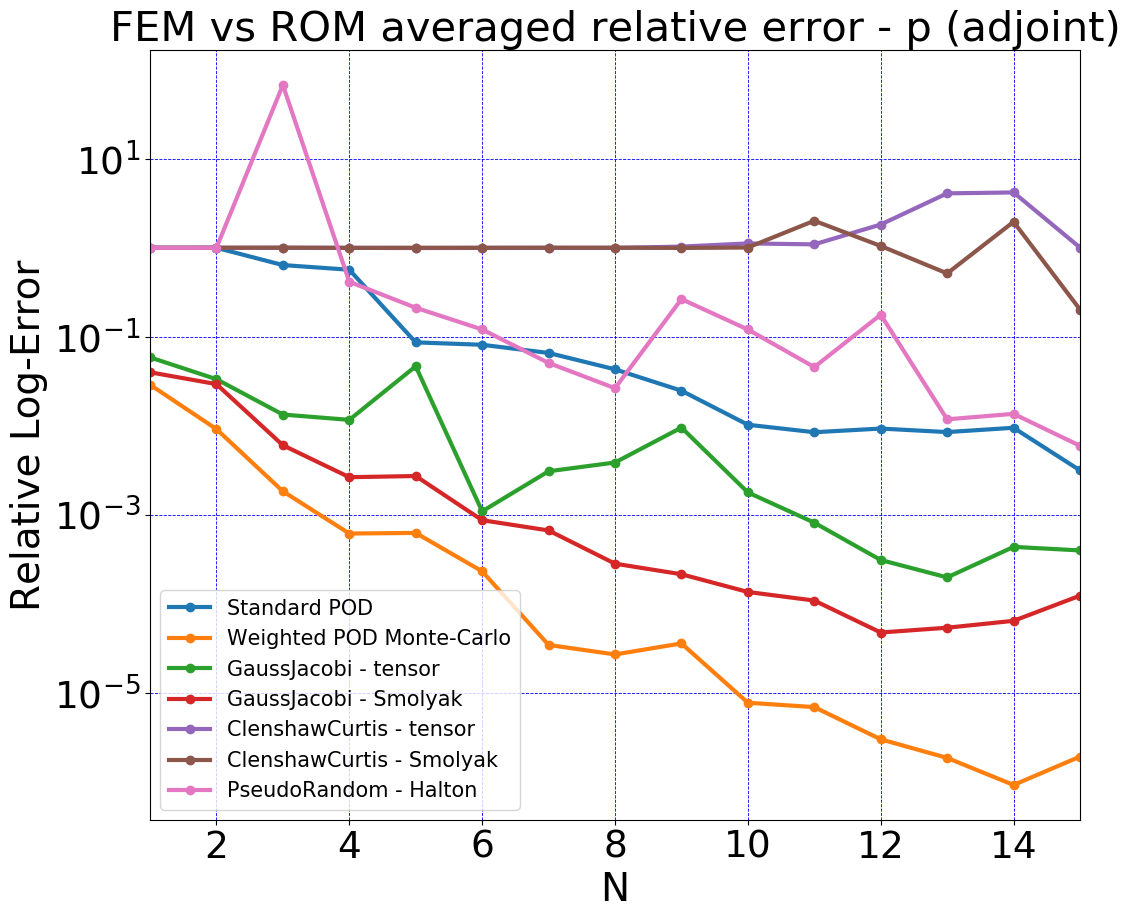}
    \caption{Relative Errors for the Parabolic Graetz-Poiseuille Problem - \textit{Offline-Online} Stabilization; State (\underline{left}), Control (\underline{center}), Adjoint (\underline{right}); Standard POD (blue), wPOD Monte-Carlo (orange), Gauss-Jacobi tensor rule (green), Gauss-Jacobi Smolyak grid (red), Clenshaw-Curtis tensor rule (cyan), Clenshaw-Curtis Smolyak grid (dark green), Pseudo-Random based on Halton numbers (pink).}
    \label{fig:plot_par_graetz_onoffstab}
\end{figure}

The initial condition is $y_0(x)=0$ \emph{for all} $x \in \Omega$ referring to Figure \ref{fig:Geometry-Graetz} and we set $T=3.0$. The penalization parameter is $\alpha=0.01$ and we want the state solution to be similar in the $L^2$-norm to a desired solution profile $y_d(x,t)=1.0$, function defined \emph{for all} $t \in (0,3.0)$ and \emph{for all} $x$ in $\Omega_{obs}$ in Figure \ref{fig:Geometry-Graetz}. Choosing $N_t = 30$, the time step is $\Delta t=0.1$. For the spatial discretization, a quite coarse mesh of $h=0.038$ is implemented: consequently, the total high-fidelity dimension is $\mathcal{N}_{tot}=314820$ and a single FEM space is characterized by $\mathcal{N}=3498$ for a fixed instant $t$ . Again, $\delta_K =1.0$ \emph{for all} $K \in \mathcal{T}_{h}$.
We take $\mathcal{P} := \big[1,10^5\big]\times \big[1, 3.0\big]$ and $\boldsymbol{\mu}$ is determined by the probability distribution
\begin{equation}\label{beta-graetz-p}
\begin{aligned}
    {\mu_1} \sim 1 + \big(10^5 - 1\big) X_1, \text{ where } X_1 \sim \text{Beta}(5,3), \\
     \mu_2 \sim 1.0 + \big( 3.0 - 1.0\big) X_2, \text{ where } X_2 \sim \text{Beta}(5,3).
\end{aligned}
\end{equation}
We choose a training set $\mathcal{P}_{\text{train}}$ of cardinality $N_{\text{train}}=100$ (with the exception of sparse grids, which have similar cardinality) and we performed the wPOD algorithms with $N_{\max}=15$. {In Figure \ref{fig:plot_Par_graetz-sing values} we show the trend of decay of the singular values. As in space-time ROM OCPs, the time instances
are not separated in the POD procedure, any space-time problem is analyzed as a steady one and
each snapshot carries all the time instances describing the whole evolution of the dynamics in $[0, T]$. That is why the values of the singular values are huge, even though always decrease over the size $N$. It is worth noticing that any Clenshaw-Curtis-based procedure does not show a good decay trend, therefore, and one should expect an underperforming outcome for their relative errors as we will show later on. Instead, all the other strategies present better results, especially the weighted Monte-Carlo, which is the method that in principle can reach better accuracy for the three variables. Instead, in Figure \ref{fig:plot_Par_graetz-proj error} the projection errors onto the POD spaces are shown. As the space-time approach is reflected in snapshots that describe the whole time instances, the differences among the studied strategies are more highlighted. The weighted Monte-Carlo and the two Gauss-Jacobi-based methods show the better possibility to reach smaller errors a priori, a possibility that will be actually reflected in the case of the relative errors in Figure \ref{fig:plot_par_graetz_onoffstab} later on, too.}
\sloppy
\sloppy
In Figure \ref{fig:plot_par_graetz_offstab} we present relative errors related to \textit{Offline-Only} stabilization. Also in the parabolic case, this procedure does not perform well. Therefore an online stabilization is needed.
\sloppy
As a matter of fact, one can see in Figure \ref{fig:plot_par_graetz_onoffstab} that the trends for \textit{Offline-Online} stabilization seem a lot better than the previous strategy. Besides the Clenshaw-Curtis quadrature rule, errors decrease along the dimension $N$. For the time-dependent cases we averaged the contribution in the testing set of the errors defined in \eqref{rel-error} where all the time instances are considered as components of the same high fidelity vector. Again, the best performance is given by the wPOD Monte-Carlo, where the following values are reached for $N=14$: $e_{y, 14} = 9.71 \cdot 10^{-7}$,$e_{p, 14}= 9.21 \cdot 10^{-7}$, and $e_{u, 14}=2.64 \cdot 10^{-7}$.
\sloppy
\sloppy
Finally, in Table \ref{speed-par-g} we illustrate the performance of the speedup-index. All weighted algorithms performs similar: we compute an order of magnitude of $10^4$ reduced solution in the time of a FEM one. This efficiency is given by the nature of the space-time procedure, where each snapshot carries {all the time instances}, and the reduction is very effective. 
\sloppy
\begin{table}
\centering

\begin{tabular}{|c|c|c|c|c|c|c|c|}
\hline  \multicolumn{8}{|c|}{Speedup-index Parabolic Graetz-Poiseuille Problem: \textit{Offline-Online} Stab. -  $\mu_1, \mu_2 \sim $ Beta(5,3)}  \\
\hline $N$ & POD & wPOD & Gauss tensor & Gauss Smolyak & CC tensor & CC Smolyak & Ps. Random\\
\hline $3$ & $14299.7$ & $14571.0$ & $13970.7$ & $14013.4$ & $14524.6$ & $14578.1$ & $14106.2$\\
\hline $6$ & $14666.3$ & $15393.5$ & $14621.8$ & $14952.8$ & $15302.6$ & $15117.8$& $14482.0$\\
\hline $9$ & $14245.6$ &  $14803.1$ & $14125.6$ & $14546.5$ & $14756.7$ & $14608.5$& $13986.9$ \\
\hline $12$ & $13693.6$ & $14206.2$  & $13554.3$ & $13935.7$ & $14050.5$ & $14075.4$ & $13453.0$\\
\hline $15$ & $13090.9$ & $13606.4$  & $13055.8$ & $13455.1$ & $13548.6$ & $13544.0$ & $12875.2$\\
\hline
\end{tabular}
  \caption{Average Speedup-index of \textit{Offline-Online} Stabilization for the Parabolic Graetz-Poiseuille Problem under geometrical parametrization. From left to right: Standard POD, wPOD Monte-Carlo, Gauss-Jacobi tensor, Gauss-Jacobi Smolyak grid, Clenshaw-Curtis tensor, Clenshaw-Curtis Smolyak grid, Pseudo-Random based on Halton numbers.}
  \label{speed-par-g}
\end{table}

\subsection{Numerical Tests for Propagating Front in a Square Problem} \label{sec:c square}
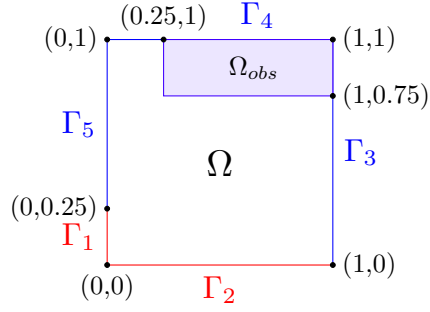
\begin{figure}[h!]
\vspace{-.3cm}
   \centering
        \begin{tikzpicture}[scale=3.0]
                  \draw[red] (0,0) -- (0,0.25) node[midway, left, scale=1.2]{$\Gamma_1$};
                 \draw[red] (1,0) -- (0,0) node[midway, below, scale=1.2]{$\Gamma_2$};
                 \draw[blue] (1,1) -- (1,0) node[midway, right, scale=1.2]{$\Gamma_3$};
                 \draw[blue] (0,1) -- (1,1) node[midway, xshift=0.5cm, above, scale=1.2]{$\Gamma_4$};
                 \draw[blue] (0,0.25) -- (0,1) node[midway, left, scale=1.2]{$\Gamma_5$};
                 \draw[black] (0.5,0.45) node[scale=1.4]{$\Omega$};
		         \draw[color=DEblue!100, fill=DEblue!10] (.25, 1) -- (0.25, 0.75) -- (1,0.75) -- (1,1) -- (.25, 1)node[midway, left, scale=1.2]{};
                 \draw[black] (0.65,0.875) node[scale=1.]{{$\Omega_{obs}$}};
                 \draw[black] (0.5,0.5) node[scale=1.4]{};
                 \filldraw[black] (0,0.25) circle (0.3pt) node[above, left]{(0,0.25)};
                 \filldraw[black] (0,1) circle (0.3pt) node[left]{(0,1)};
                 \filldraw[black] (1,0.75) circle (0.3pt) node[right]{(1,0.75)};
                 \filldraw[black] (1,1) circle (0.3pt) node[right]{(1,1)}; \filldraw[black] (0.25,1) circle (0.3pt) node[above]{(0.25,1)};
                 \filldraw[black] (1,0) circle (0.3pt) node[below, right]{(1,0)};
                 \filldraw[black] (0,0) circle (0.3pt) node[below]{(0,0)};
        \end{tikzpicture} 
        \caption{Geometry of the Propagating Front in a Square Problem}
        \label{fig:geometry-square}
\end{figure}

Here we analyze an Advection-Dominated PDE problem illustrated without control in a deterministic and in a stochastic context in \cite{pacciarini2014stabilized,torlo2018stabilized} and in \cite{torlo2018stabilized}, respectively.
A distributed control is applied all over the domain $\Omega$, which is the square   $(0,1) \times (0,1)$, as shown under Cartesian coordinates $(x_0,x_1)$ in Figure \ref{fig:geometry-square}. The boundary is composed as follows: $\Gamma_1 := \{0\} \times [0,0.25]$, $\Gamma_2 := [0,1] \times \{0\}$, $\Gamma_3 := \{1\} \times [0,1]$, $\Gamma_4 := [0,1] \times \{1\}$, $\Gamma_5 := \{0\} \times [0.25,1]$; instead $\Omega_{obs}:= [0.25,1]\times [0.75,1]$.
Given $\boldsymbol{\mu} = (\mu_1, \mu_2)$, our aim is to solve the following OCP($\boldsymbol{\mu}$) problem: find $(y, u) \in \tilde{Y} \times U$ which solves
\begin{equation*}
\min_{(y,u)} \displaystyle \frac{1}{2} \int_{\Omega_{obs}}(y(\boldsymbol{\mu})- y_d)^2 \; d \Omega \vspace{.5cm} +
\frac{\alpha}{2} \int_{\Omega}u(\boldsymbol{\mu})^2 \; d \Omega, \ \text{  such that }
\vspace{-.6cm}
\end{equation*}
\begin{equation}\label{square-system}
\begin{cases}
\displaystyle
-\frac{1}{\mu_1} \Delta y(\boldsymbol{\mu})+[\cos{\mu_2},\sin{\mu_2}] \cdot \nabla y(\boldsymbol{\mu})=u(\boldsymbol{\mu}), & \text { in } \Omega, \\
\displaystyle
y(\boldsymbol{\mu})=1, & \text { on } \Gamma_{1} \cup \Gamma_{2}, \\
\displaystyle
y(\boldsymbol{\mu})=0, & \text { on } \Gamma_{3} \cup \Gamma_{4} \cup \Gamma_{5}.
\end{cases}
\end{equation}
In this case, we have that {the domain of definition of our state $y$ is $$\tilde{Y}:= \big\{v \in H^{1} \big(\Omega\big) \text{ s.t. } \mathrm{BC} \text{ in }  (\ref{square-system}) \big\}.$$ 
Again, we define a lifting function $R_y \in H^1\big(\Omega\big)$ such that satisfies $\mathrm{BC} \text{ in}$ (\ref{square-system}), applying a lifting procedure before the Lagrangian approach. We define $\bar{y} := y - R_y$, with $\bar{y} \in Y$ and $Y:= H^1_{0}(\Omega)$, }$U = \mathrm{L}^{2}(\Omega)$ and $Q := Y^{*}$, with $p=0$ on $\partial \Omega$.

\begin{figure}
        \centering
        \includegraphics[scale=0.2]{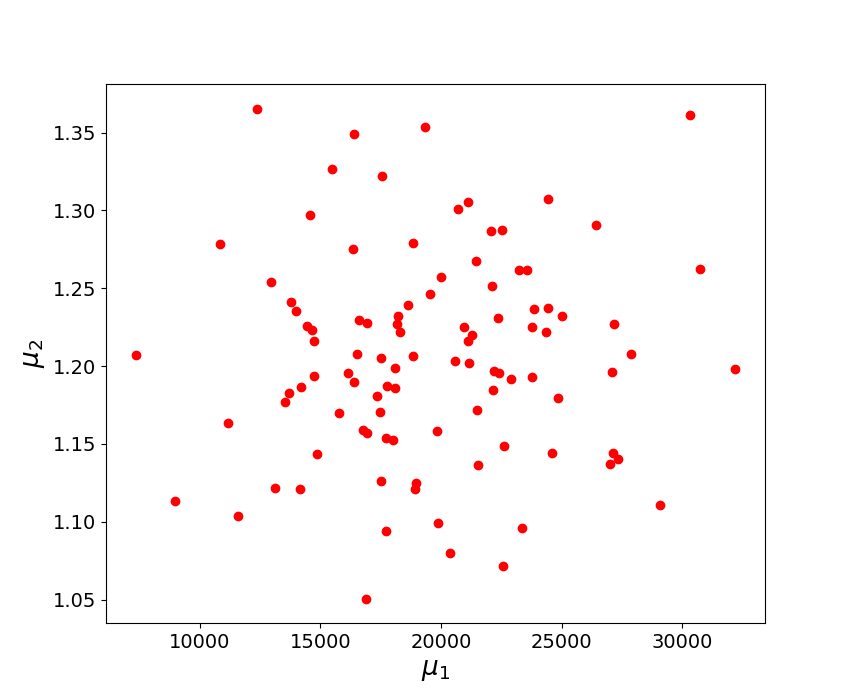}  
        \includegraphics[scale=0.2]{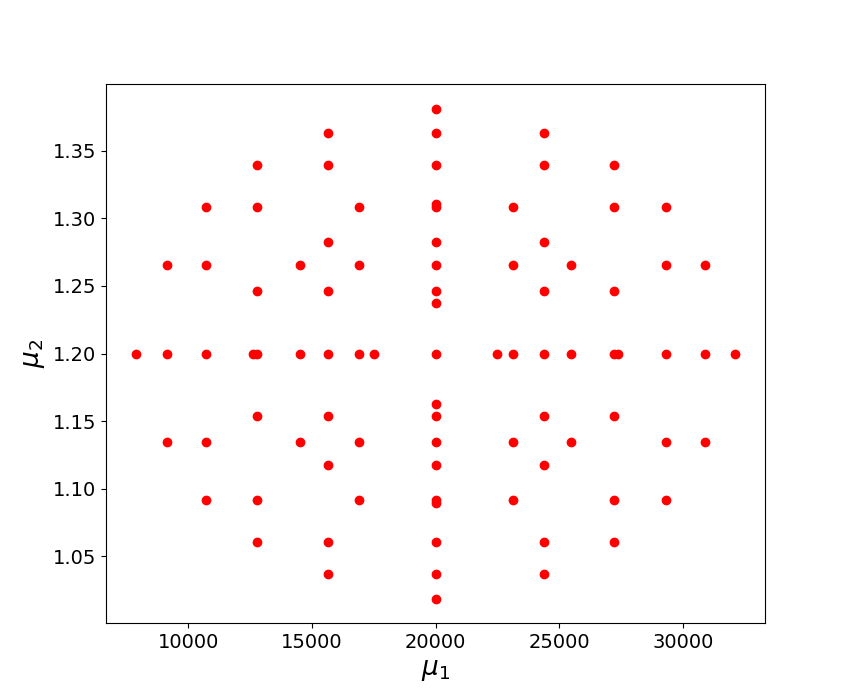} 
         \includegraphics[scale=0.2]{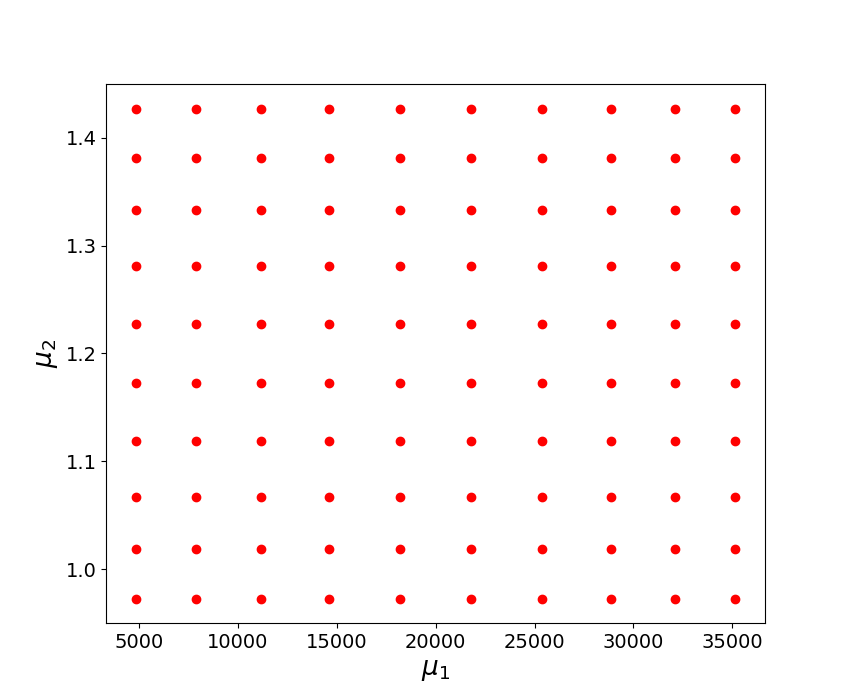} \\
        \includegraphics[scale=0.2]{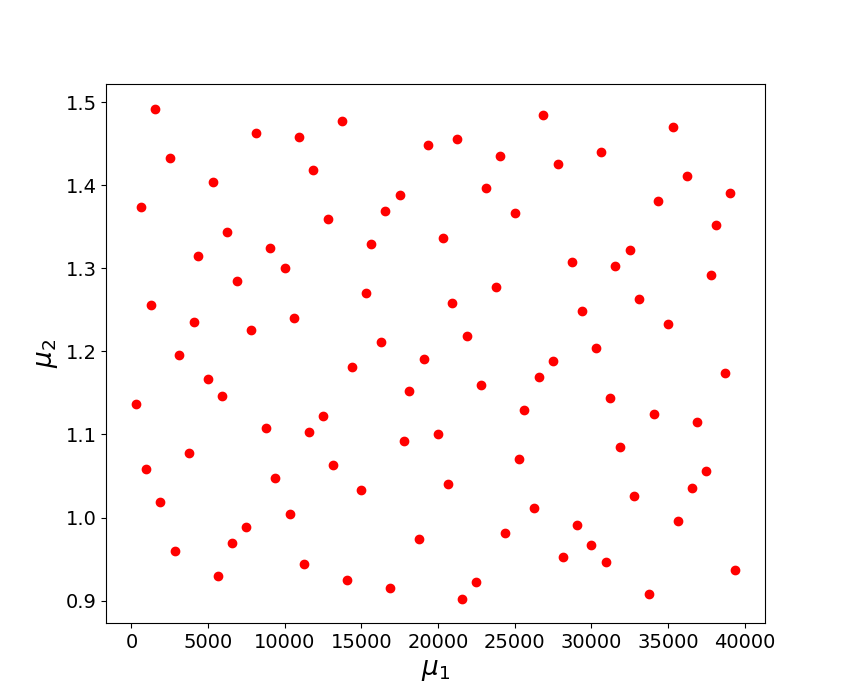} 
        \includegraphics[scale=0.2]{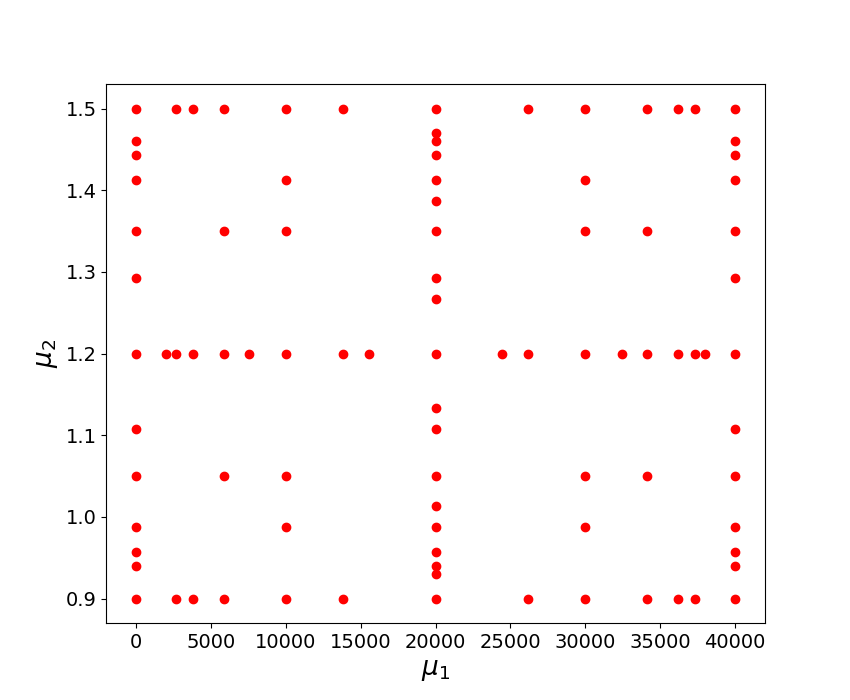} 
        \includegraphics[scale=0.2]{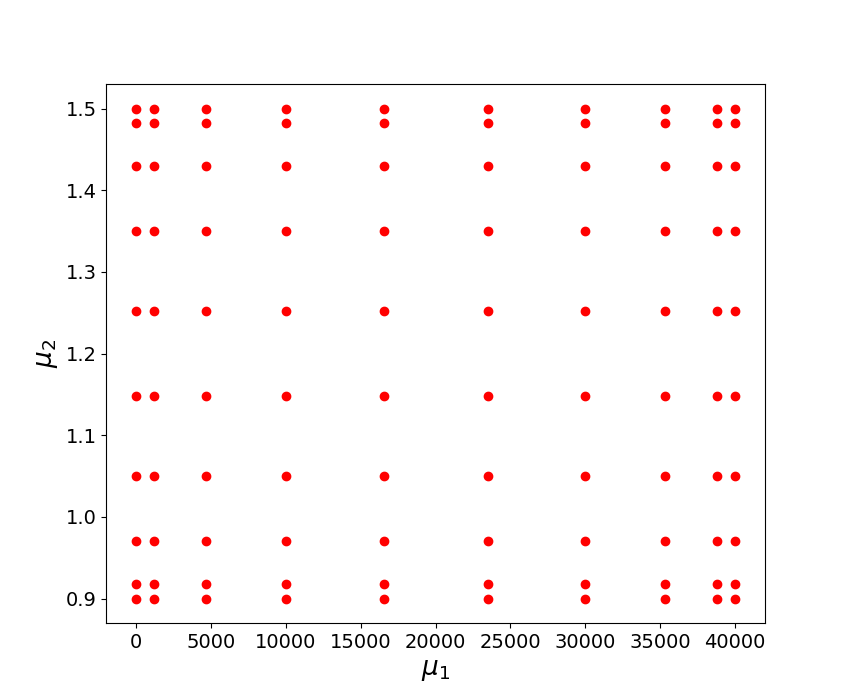} 
        \caption{Grid points for the quadrature formulae of the weighted POD regarding the Propagating Front in a Square Problem; (\underline{top}) Monte-Carlo method with $\boldsymbol{\mu}$ following distribution \eqref{beta-square} on the parameter space $\mathcal{P}$ (\underline{left}), Smolyak grid based on a Gauss-Jacobi rule (\underline{center}), Tensor grid based on a Gauss-Jacobi rule (\underline{right}); (\underline{bottom}) Pseudo-Random method based on a Halton sequence (\underline{left}), Smolyak grid based on a Clenshaw-Curtis rule (\underline{center}), Tensor grid based on a Clenshaw-Curtis rule (\underline{right}).}
        \label{fig:square-grids}
\end{figure}

The mesh size $h$ is equal to $0.025$, which entails an overall dimension of the truth approximation of $12087$. Consequently, we have $\mathcal{N}=4029$ for state, control and adjoint spaces. Concerning stabilization, $\delta_K =1.0$ \emph{for all} $K \in \mathcal{T}_{h}$. The penalization parameter is $\alpha=0.01$ and we pursue {the state solution to be similar in the $L^2$-norm to} $y_d(x)=0.5$, defined \emph{for all} $x$ in $\Omega_{obs}$ of Figure \ref{fig:geometry-square}. 
In our test cases, $\mathcal{P}:=\big[1,4\cdot 10^4\big] \times \big[0.9,1.5\big]$ and $\boldsymbol{\mu}$ follow the subsequent probability distribution:
\begin{equation}\label{beta-square}
\begin{aligned}
    {\mu_1} \sim 1 + \big(4 \cdot 10^4 - 1\big) X_1, \text{ where } X_1 \sim \text{Beta}(10,10), \\
     \mu_2 \sim 0.9 + \big( 1.5 - 0.9\big) X_2, \text{ where } X_2 \sim \text{Beta}(10,10),
\end{aligned}
\end{equation}
where $\mu_1$ and $\mu_2$ are independent random variables.
The training set $\mathcal{P}_{\text{train}}$ and the testing set $\mathcal{P}_{\text{test}}$ have both cardinality equal to $n_{\text{train}}=100$, with exception of sparse grid samplings, whose cardinality is similar to $100$. In Figure \ref{fig:square-grids} the grid points for different sampling strategies of the weighted POD are illustrated. Similarly to the other example, by using tensor rules or the Smolyak technique with Clenshaw-Curtis quadrature, several points are taken close to the boundary. Thus, one would expect a lower performance of those rules according to the distribution \eqref{beta-square}.
\sloppy

\begin{figure}
        \centering
        \includegraphics[scale=0.16]{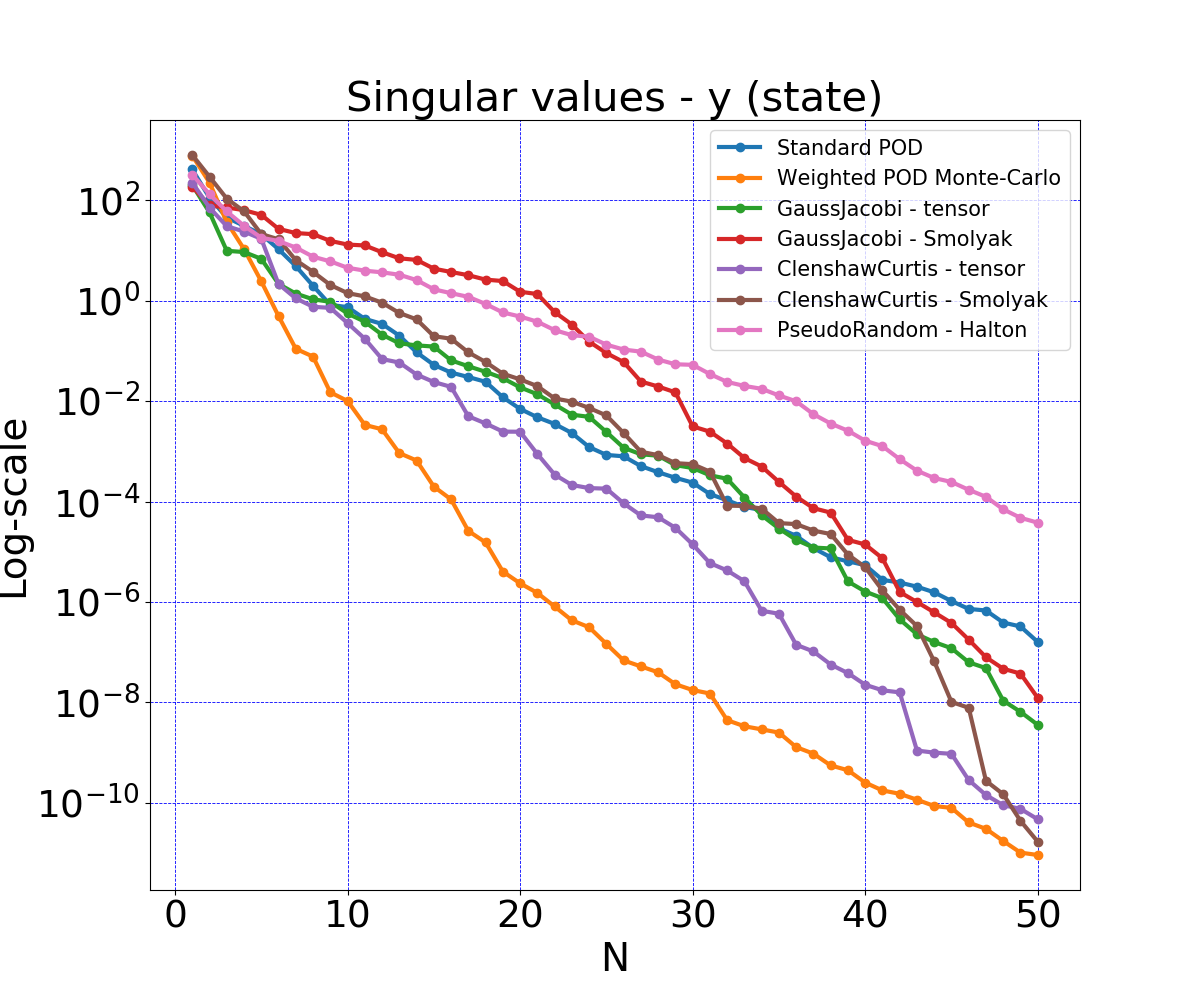} 
        \includegraphics[scale=0.16]{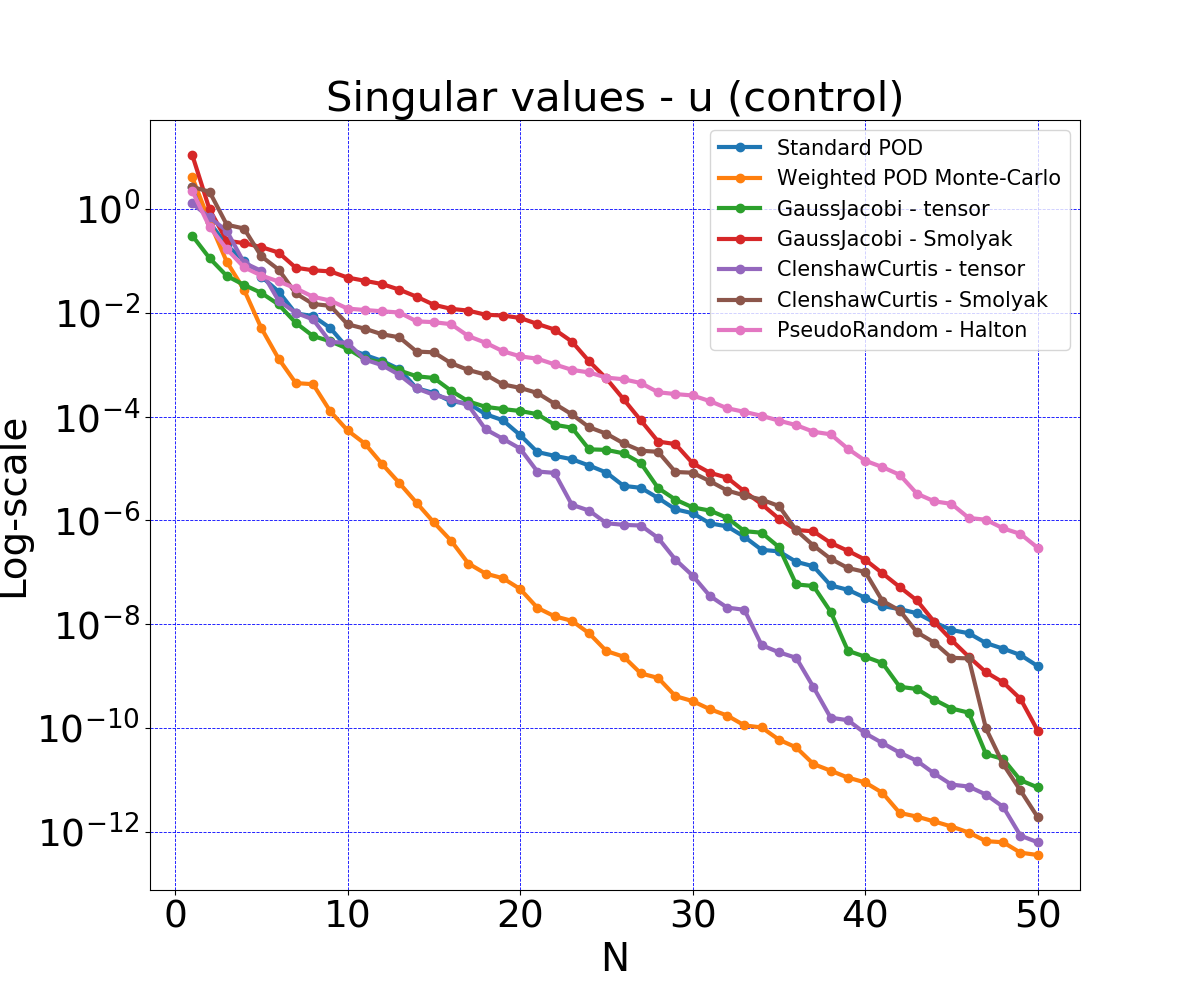} 
        \includegraphics[scale=0.16]{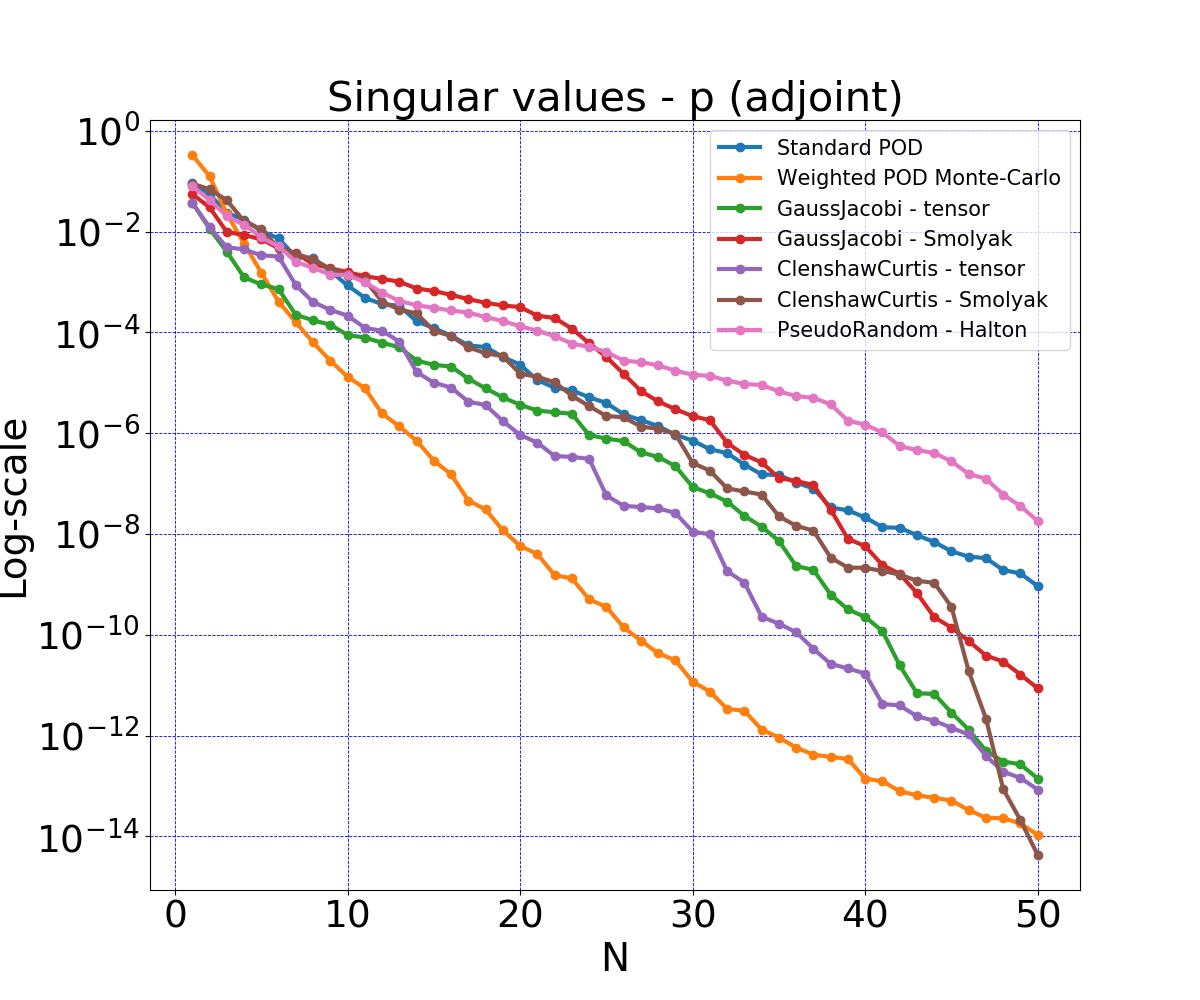} 
        \caption{Singular values decay for the snapshot matrices for the Propagating Front in a Problem with $\boldsymbol{\mu}$ following distribution \eqref{beta-square} on the parameter space $\mathcal{P}$; State (\underline{left}), Control (\underline{center}), Adjoint (\underline{right}); Standard POD (blue), wPOD Monte-Carlo (orange), Gauss-Jacobi tensor rule (green), Gauss-Jacobi Smolyak grid (red), Clenshaw-Curtis tensor rule (cyan), Clenshaw-Curtis Smolyak grid (dark green), Pseudo-Random based on Halton numbers (pink).}
        \label{fig:plot_Square-sing values}
\end{figure}

\begin{figure}
        \centering
        \includegraphics[scale=0.16]{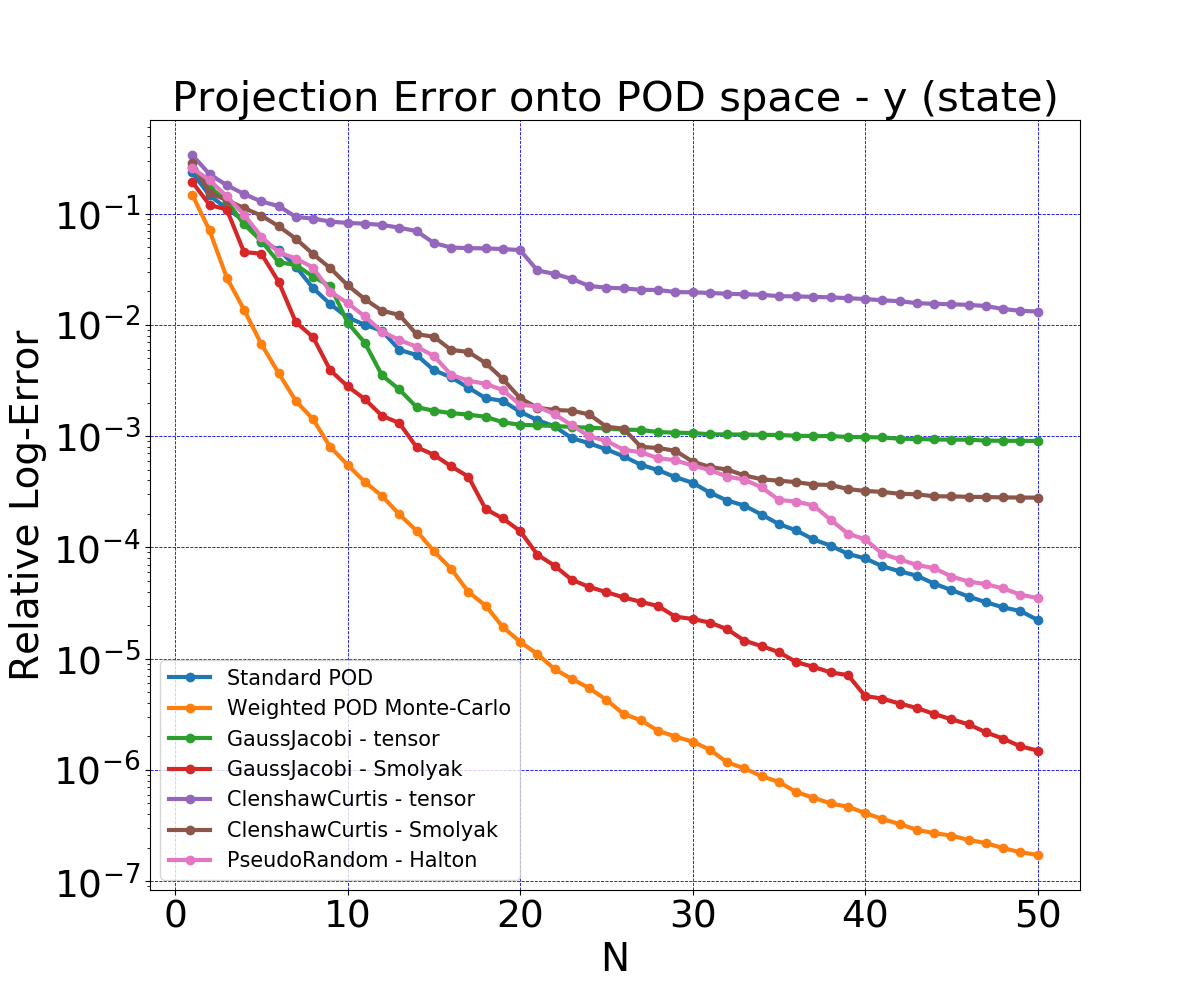} 
        \includegraphics[scale=0.16]{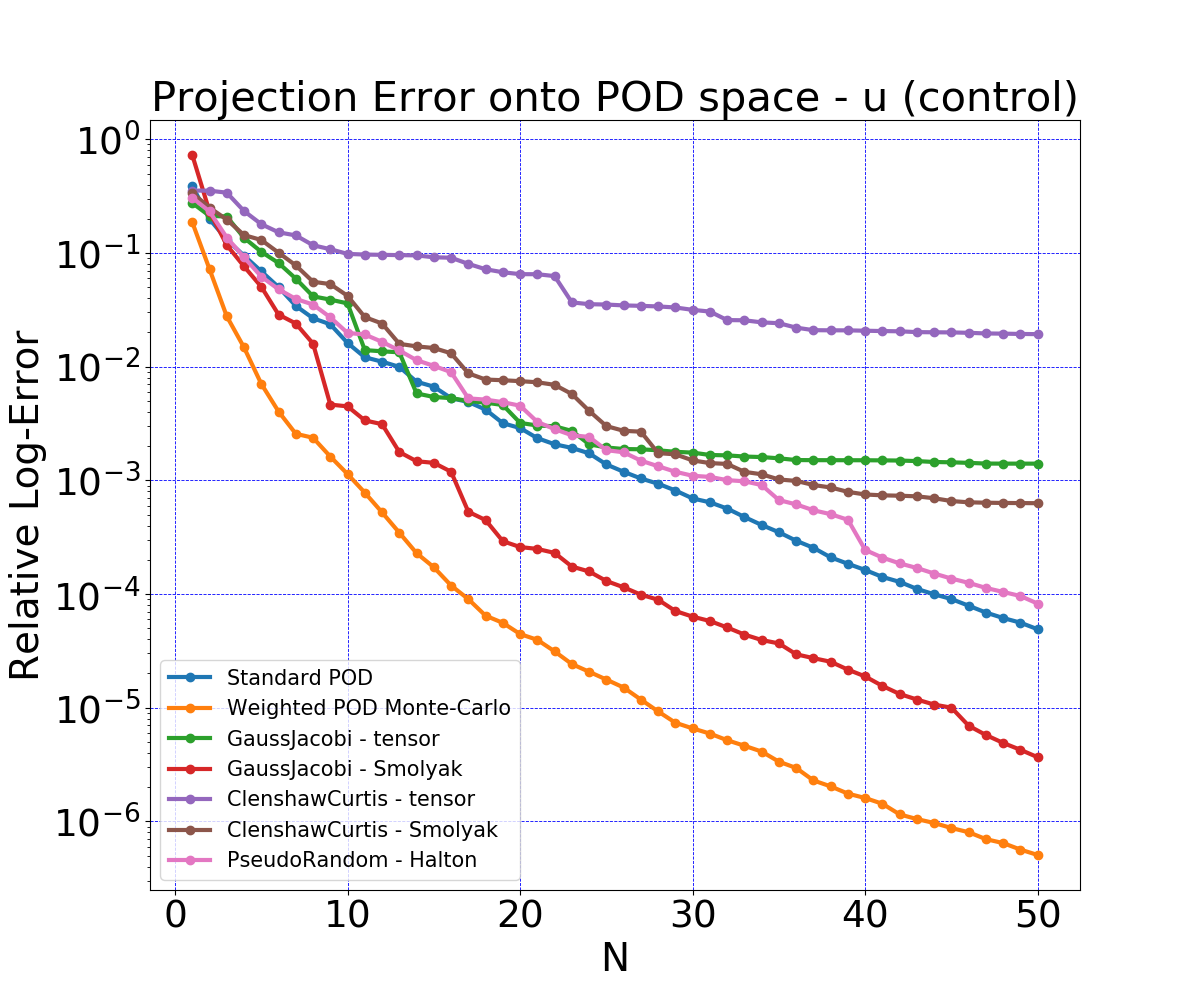} 
        \includegraphics[scale=0.16]{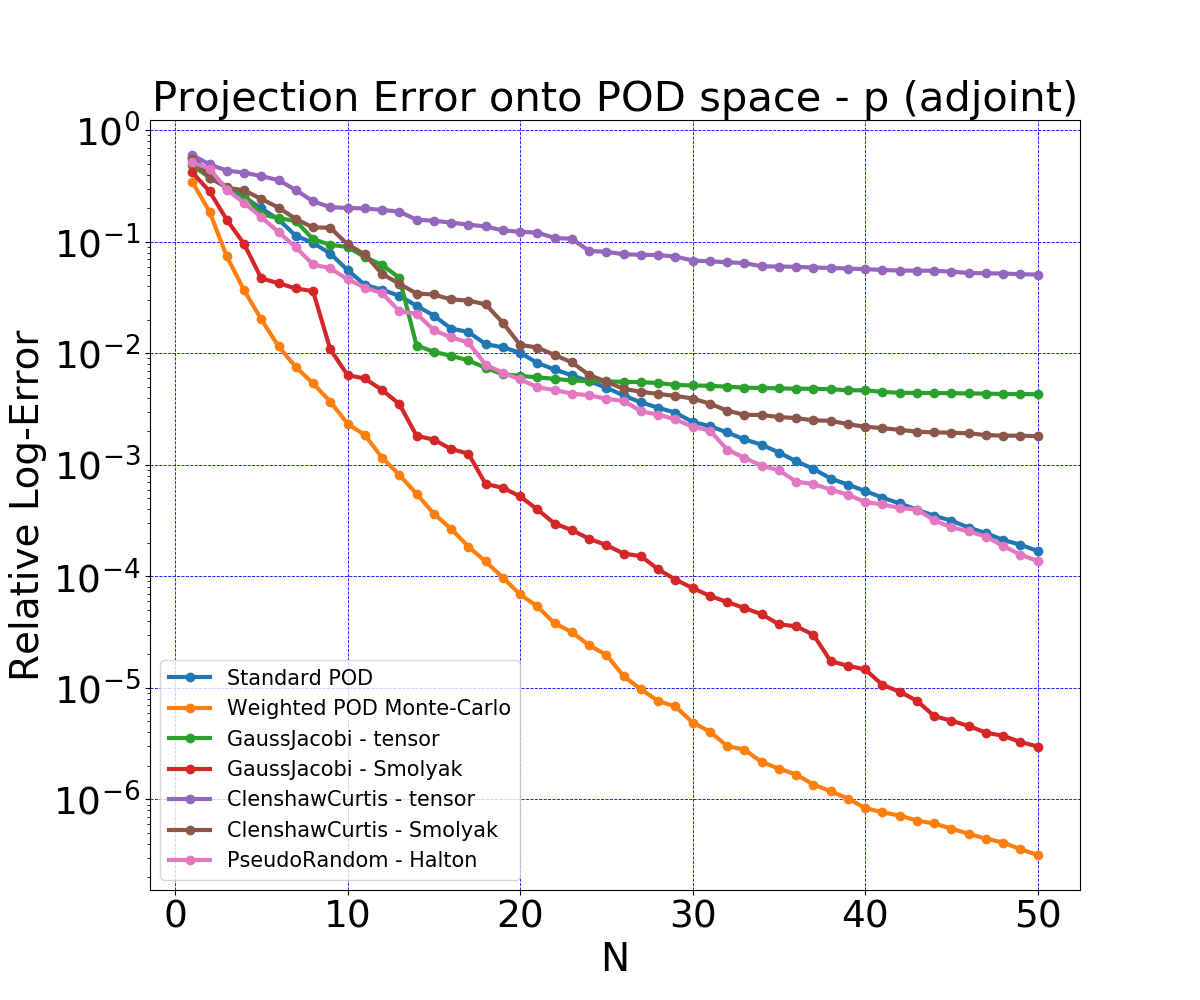} 
        \caption{Projection Errors onto the POD space for the Propagating Front in a Problem with $\boldsymbol{\mu}$ following distribution \eqref{beta-square} on the parameter space $\mathcal{P}$; State (\underline{left}), Control (\underline{center}), Adjoint (\underline{right}); Standard POD (blue), wPOD Monte-Carlo (orange), Gauss-Jacobi tensor rule (green), Gauss-Jacobi Smolyak grid (red), Clenshaw-Curtis tensor rule (cyan), Clenshaw-Curtis Smolyak grid (dark green), Pseudo-Random based on Halton numbers (pink).}
        \label{fig:plot_Square-proj error}
\end{figure}

\begin{figure}
        \centering
        \includegraphics[scale=0.129]{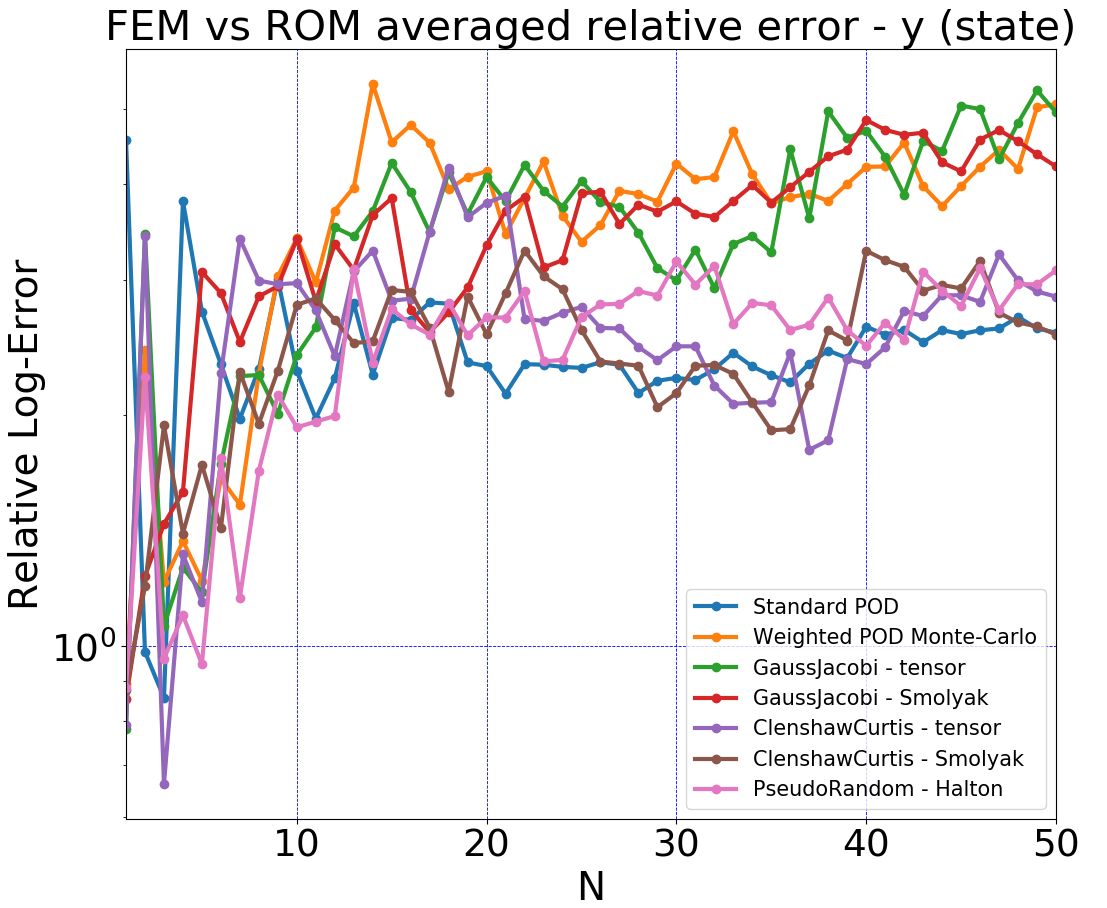}
        \includegraphics[scale=0.129]{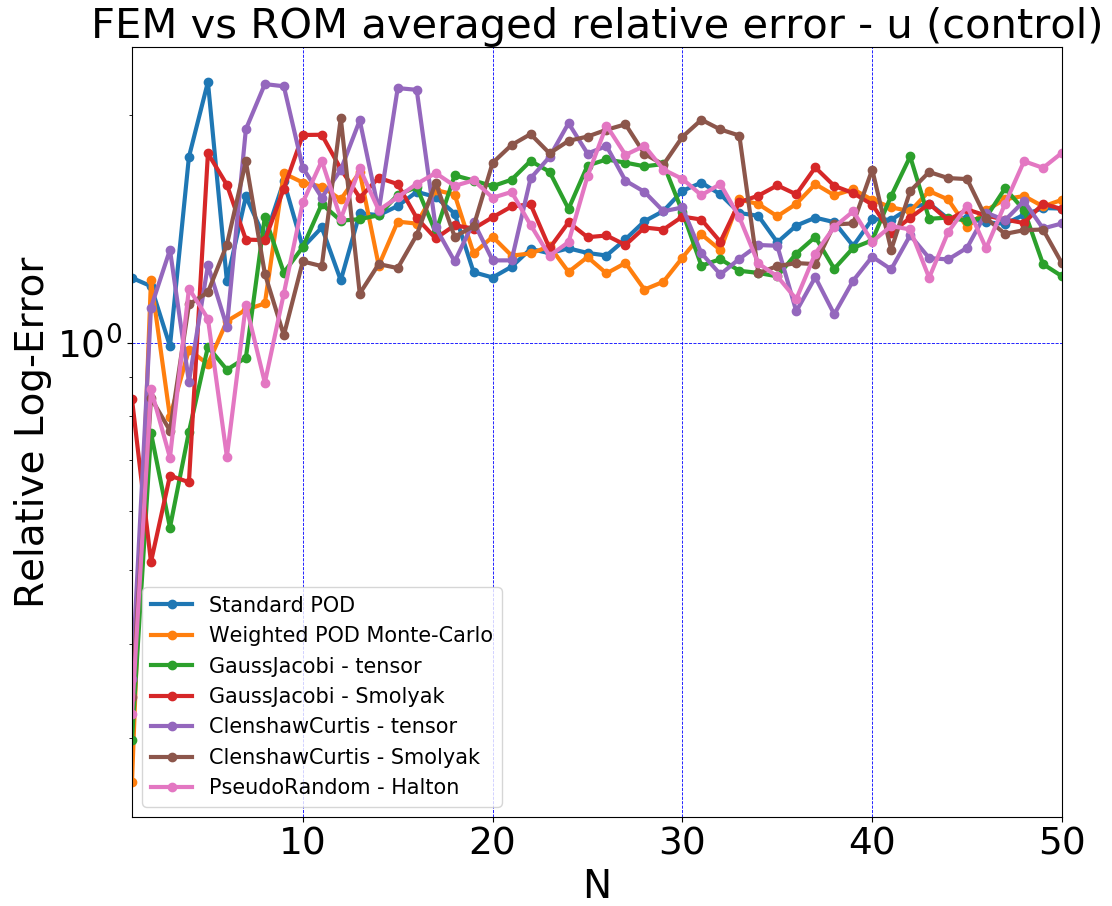}
        \includegraphics[scale=0.129]{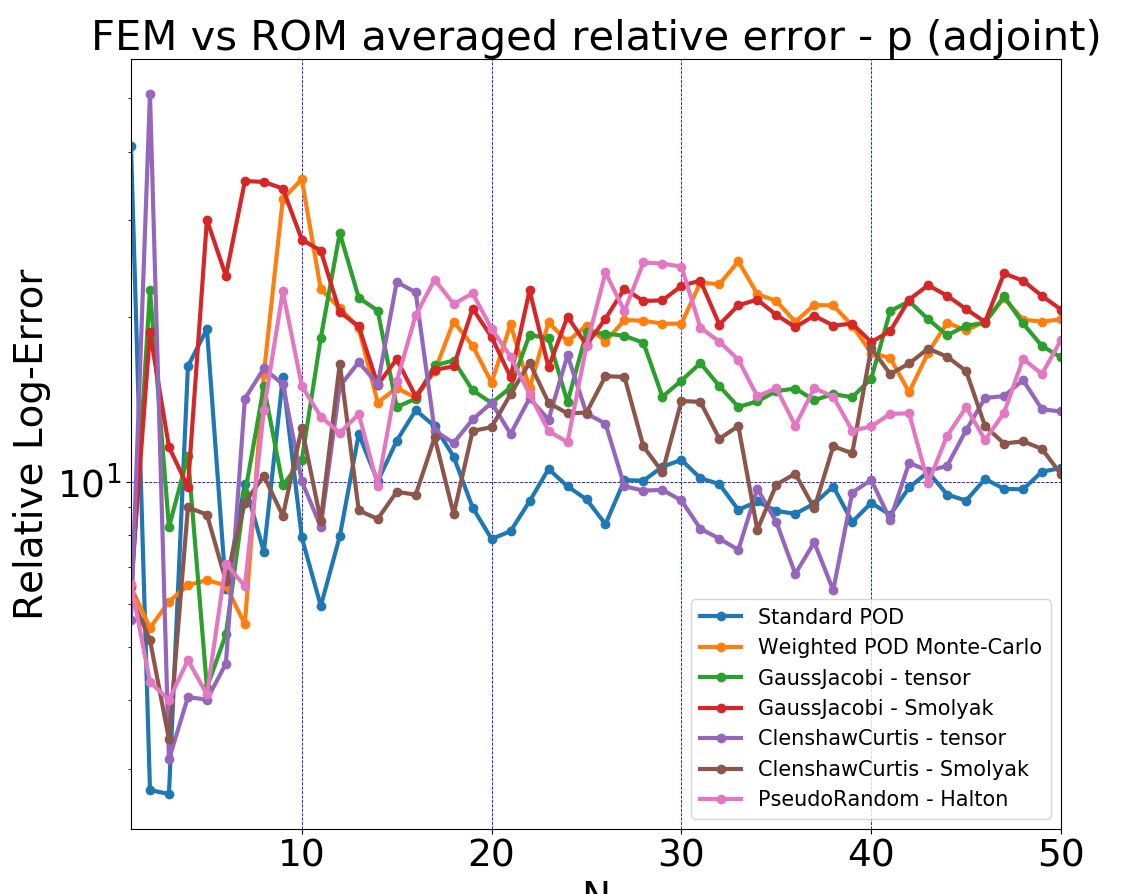}
        \caption{Relative Errors for the Propagating Front in a Problem with $\boldsymbol{\mu}$ following distribution \eqref{beta-square} on the parameter space $\mathcal{P}$- \textit{Offline-Only} Stabilization; State (\underline{left}), Control (\underline{center}), Adjoint (\underline{right}); Standard POD (blue), wPOD Monte-Carlo (orange), Gauss-Jacobi tensor rule (green), Gauss-Jacobi Smolyak grid (red), Clenshaw-Curtis tensor rule (cyan), Clenshaw-Curtis Smolyak grid (dark green), Pseudo-Random based on Halton numbers (pink).}
        \label{fig:square_off_stab}
\end{figure}

\begin{figure}
        \centering
        \includegraphics[scale=0.127]{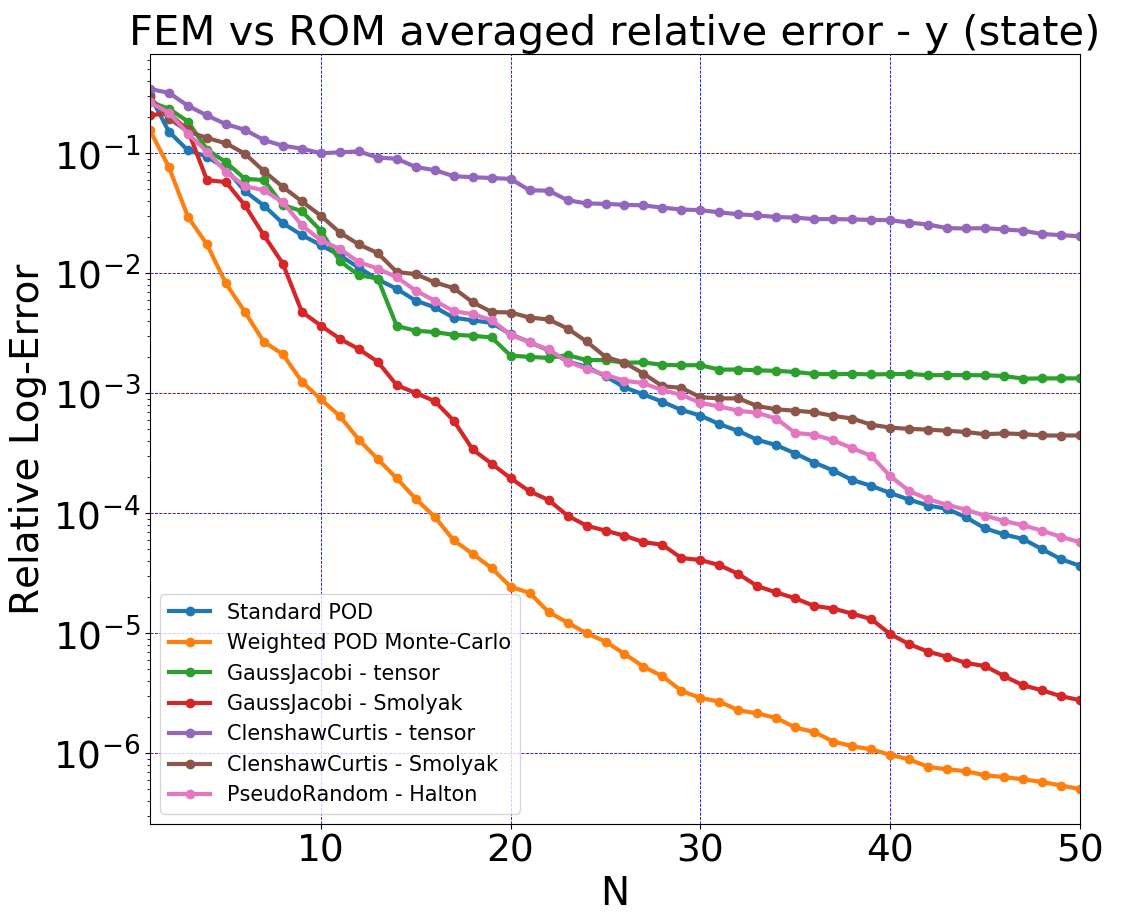}
        \includegraphics[scale=0.127]{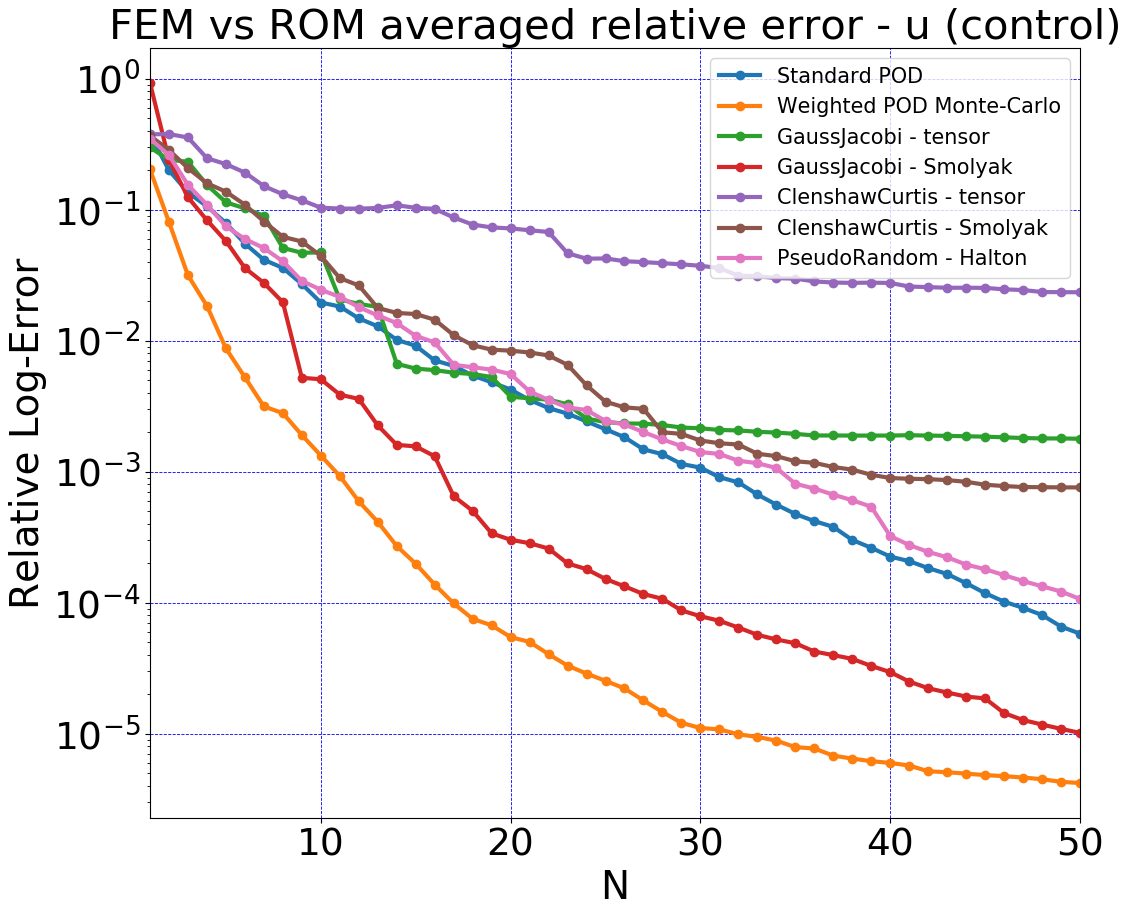}
        \includegraphics[scale=0.127]{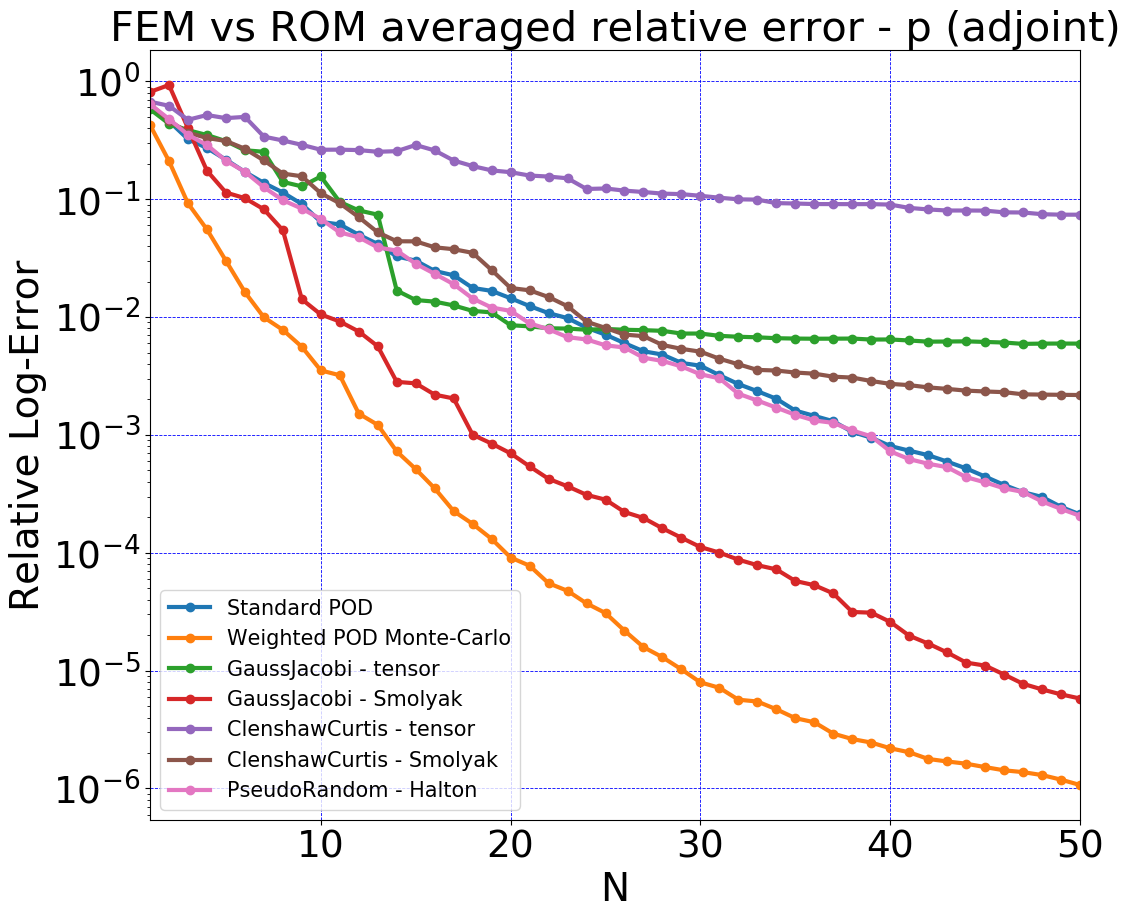}
        \caption{Relative Errors for the Propagating Front in a Problem with $\boldsymbol{\mu}$ following distribution \eqref{beta-square} on the parameter space $\mathcal{P}$- \textit{Offline-Online} Stabilization; State (\underline{left}), Control (\underline{center}), Adjoint (\underline{right}); Standard POD (blue), wPOD Monte-Carlo (orange), Gauss-Jacobi tensor rule (green), Gauss-Jacobi Smolyak grid (red), Clenshaw-Curtis tensor rule (cyan), Clenshaw-Curtis Smolyak grid (dark green), Pseudo-Random based on Halton numbers (pink).}
        \label{fig:square_onoff_stab}
\end{figure}

\begin{figure}
        \centering
        \includegraphics[scale=0.14]{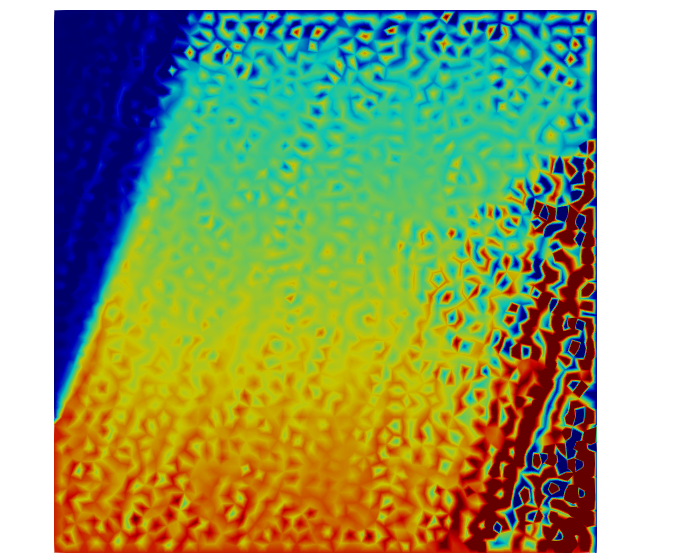}
        \includegraphics[scale=0.14]{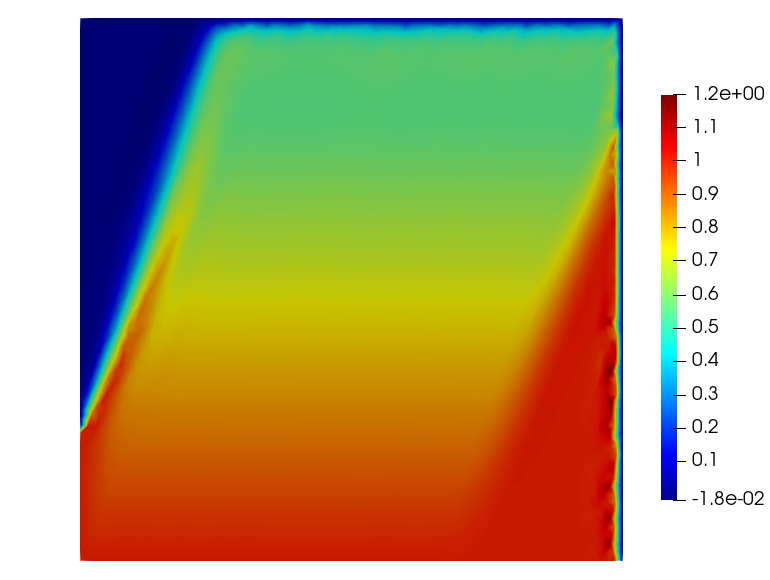}
        \includegraphics[scale=0.14]{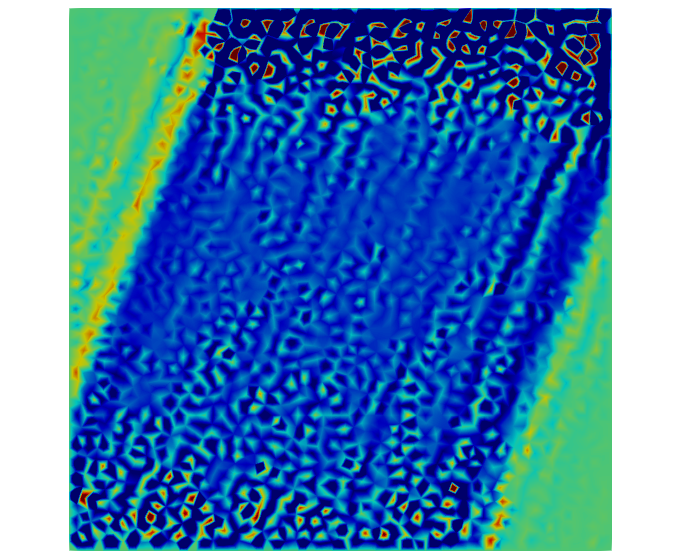}
        \includegraphics[scale=0.14]{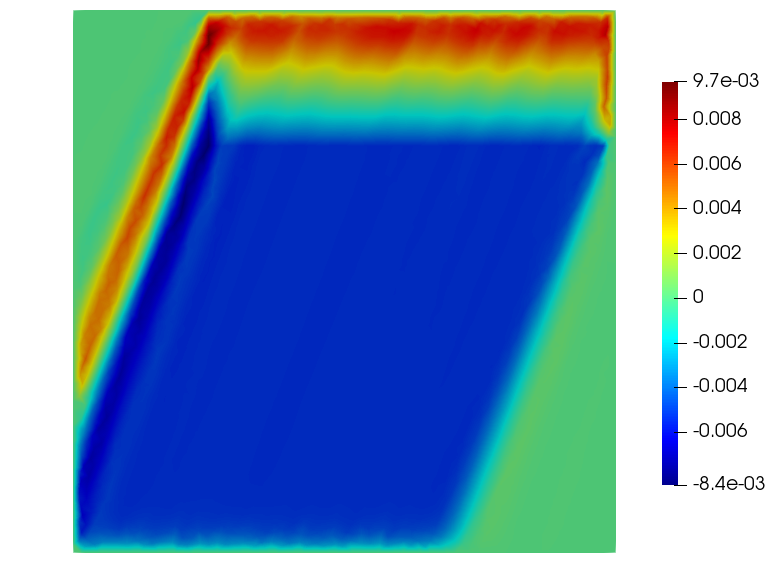}
        \caption{FEM not stabilized and wPOD Monte-Carlo \textit{Offline-Online} stabilized solution for $y$ (\underline{left}) and for $p$ (\underline{right}), $\boldsymbol{\mu} = (2 \cdot 10^4, 1.2)$, $h=0.025$ $\alpha=0.01$, $N_{\text{train}}=100$, $\delta_K=1.0$, $N=50$.}
        \label{fig:square_figures}
\end{figure}

We apply a wPOD procedure for a $N_{\operatorname{max}}=50$ dimension. {In Figures \ref{fig:plot_Square-sing values} and \ref{fig:plot_Square-proj error}, we plot the singular value decay for the snapshot matrices and the projection errors for the three variables, respectively. Concerning the singular values, we see that the best performance still belongs to the weighted Monte-Carlo. It is the case for the projection errors, too, where we also notice that some procedures (Clenshaw-Curtis and Gauss-Jacobi tensor rule) fill their accuracy before reaching $N=50$. This implies that a priori, those methods are not the best strategies to build a weighted POD model for this kind of problem. }
\sloppy
In Figure \ref{fig:square_off_stab}, we show the performance of relative errors for the \textit{Offline-Only} stabilization procedure. As in the Graetz-Poiseuille Problem, these trends are not acceptable, as no quantity drops under $10^{-1}$ for all state, control and adjoint variables. {All these quantities do not follow the same behaviors as the projection errors in Figure \ref{fig:plot_Square-proj error}.} Therefore, a stabilization applied in the Online Phase is needed, too.
\sloppy
\sloppy
In Figure \ref{fig:square_onoff_stab}, relative errors for \textit{Offline-Online} Stabilization procedure are shown. {In this case, the relative errors decrease in the same way as the projection ones in Figure \ref{fig:plot_Square-proj error}.} Again, wPOD Monte-Carlo presents the best behaviour: in this case, it reaches $e_{y, 50} = 5.03 \cdot 10^{-7}$ for the state, for the adjoint $e_{p, 50}= 1.07 \cdot 10^{-6}$, and the control $e_{u, 50}=4.21 \cdot 10^{-6}$. Moreover, the wPOD Monte-Carlo has an accuracy of nearly a factor of $100$ better than a Standard POD in a deterministic context for $N>20$. Also here, Smolyak grids perform better than their tensor counterpart: for instance, we obtain in this case it reaches $e_{y, 50} = 2.77 \cdot 10^{-6}$ for the state, for the adjoint $e_{p, 50}= 5.80 \cdot 10^{-6}$, and the control $e_{u, 50}=1.02 \cdot 10^{-5}$ for Gauss-Jacobi. Concerning the training set, we have $N_{\text{train}}=89$ and $N_{\text{train}}=93$ for the Gauss-Jacobi and the Clenshaw-Curtis ones, respectively. In Figure \ref{fig:square_figures} we see a comparison between the FEM solution for the state and the adjoint without stabilization and the \textit{Offline-Online} Stabilized wPOD Monte-Carlo reduced solution for these variables with $\boldsymbol{\mu} = (2 \cdot 10^4, 1.2)$.
\sloppy
\sloppy
\sloppy
The values of the speedup-index for the \textit{Offline-Online} stabilization for each type of wPOD are reported in Table \ref{speed-s}. For $N=50$ the wPOD Monte-Carlo is the best choice again with a computation of $50$ reduced solutions in the time of a FEM one. All the other possibilities perform a little bit lower for $N=50$; however, all weighted algorithms have similar performances concerning the speedup-index: an order of magnitude of $10^2$ for the first $50$ reduced basis.
\sloppy
\begin{table}
\centering

\begin{tabular}{|c|c|c|c|c|c|c|c|}
\hline  \multicolumn{8}{|c|}{Speedup-index Propagating front in a Square Problem: \textit{Offline-Online} Stab. - $\mu_1, \mu_2 \sim $ Beta(10,10)}  \\
\hline $N$ & POD & wPOD & Gauss tensor & Gauss Smolyak & CC tensor & CC Smolyak & Ps. Random \\
\hline $10$ & $151.3$ & $179.2$  &  $175.0$ & $178.7$ & $181.5$ & $176.4$ & $173.9$\\
\hline $20$ & $123.3$ &  $140.4$ & $139.8$ & $141.0$ & $140.9$ & $140.5$ & $143.6$\\
\hline $30$ & $88.5$ & $103.3$  & $102.6$ & $102.8$ & $100.6$ & $102.6$ & $104.3$\\
\hline $40$ & $61.6$ & $73.7$  & $73.2$ & $69.9$ & $68.6$ & $70.4$ & $70.2$ \\
\hline $50$ & $43.4$ & $50.2$  & $49.0$ & $47.6$ & $46.8$ & $49.2$ & $48.2$ \\
\hline
\end{tabular}
  \caption{Average Speedup-index of \textit{Offline-Online} Stabilization for the Propagating Front in a Square Problem. From left to right: Standard POD, wPOD Monte-Carlo, Gauss-Jacobi tensor, Gauss-Jacobi Smolyak grid, Clenshaw-Curtis tensor, Clenshaw-Curtis Smolyak grid, Pseudo-Random based on Halton numbers.}
  \label{speed-s}
\end{table}

\begin{figure}
    \centering
    \includegraphics[scale=0.15]{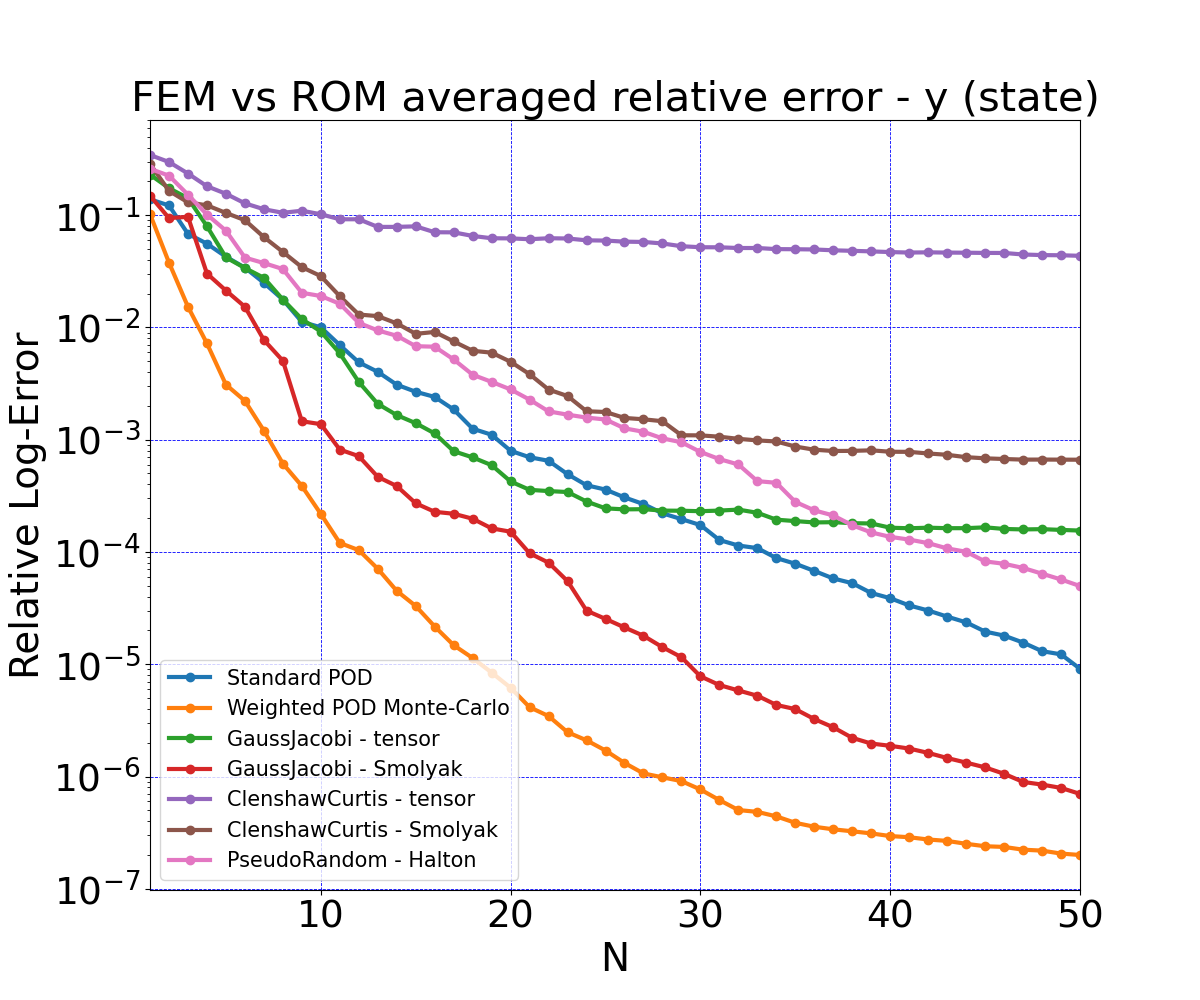}
    \includegraphics[scale=0.15]{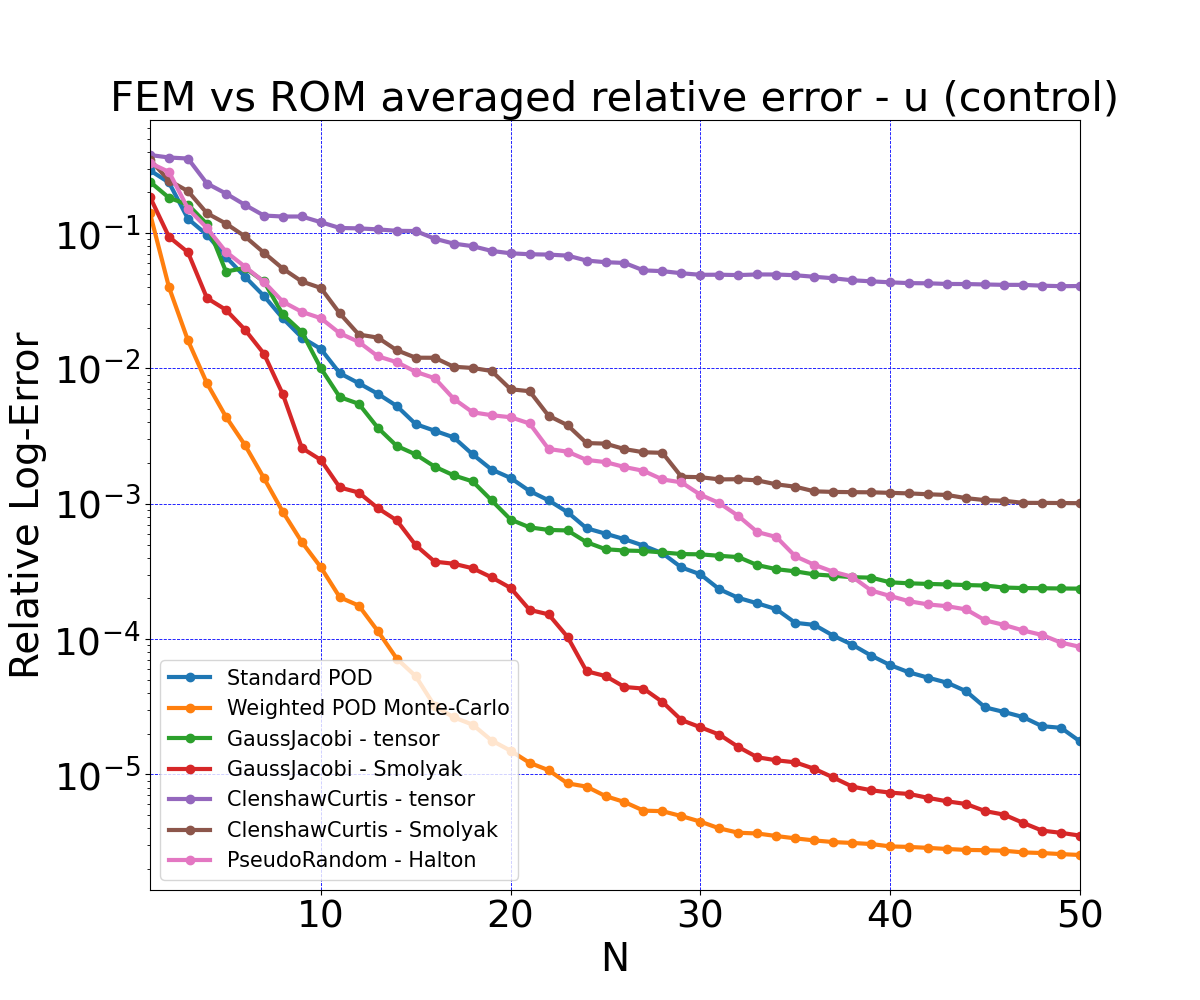}
    \includegraphics[scale=0.15]{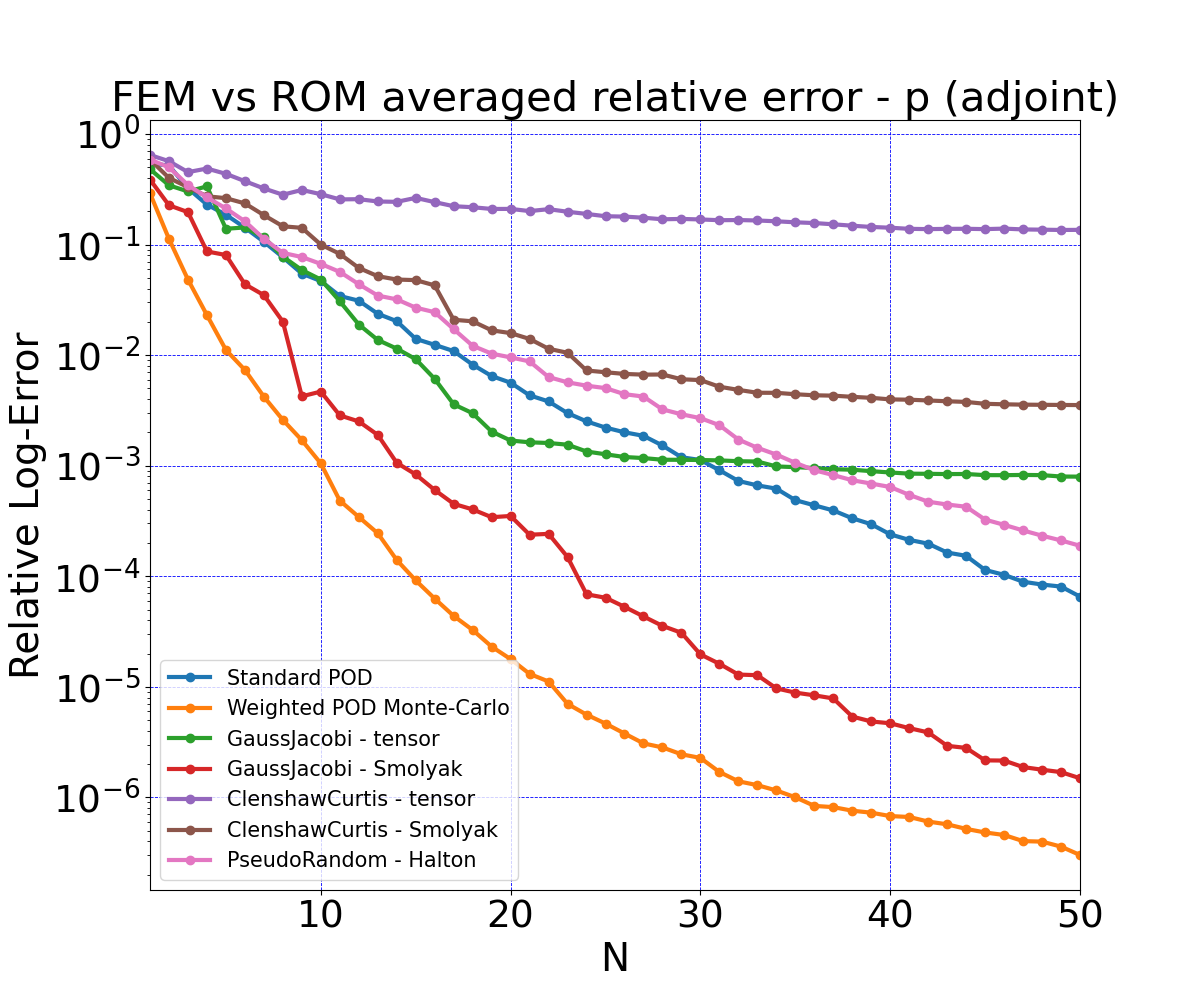}
    \caption{Relative Errors for the Propagating Front in a Square Problem with $\boldsymbol{\mu}$ following distribution \eqref{beta-square-2} on the parameter space $\mathcal{P}$ - \textit{Offline-Online} Stabilization; State (\underline{left}), Control (\underline{center}), Adjoint (\underline{right}); Standard POD (blue), wPOD Monte-Carlo (orange), Gauss-Jacobi tensor rule (green), Gauss-Jacobi Smolyak grid (red), Clenshaw-Curtis tensor rule (cyan), Clenshaw-Curtis Smolyak grid (dark green), Pseudo-Random based on Halton numbers (pink).}
    \label{fig:plot_square_onoffstab-2}
\end{figure}

\begin{figure}
     \centering
    \includegraphics[scale=0.15]{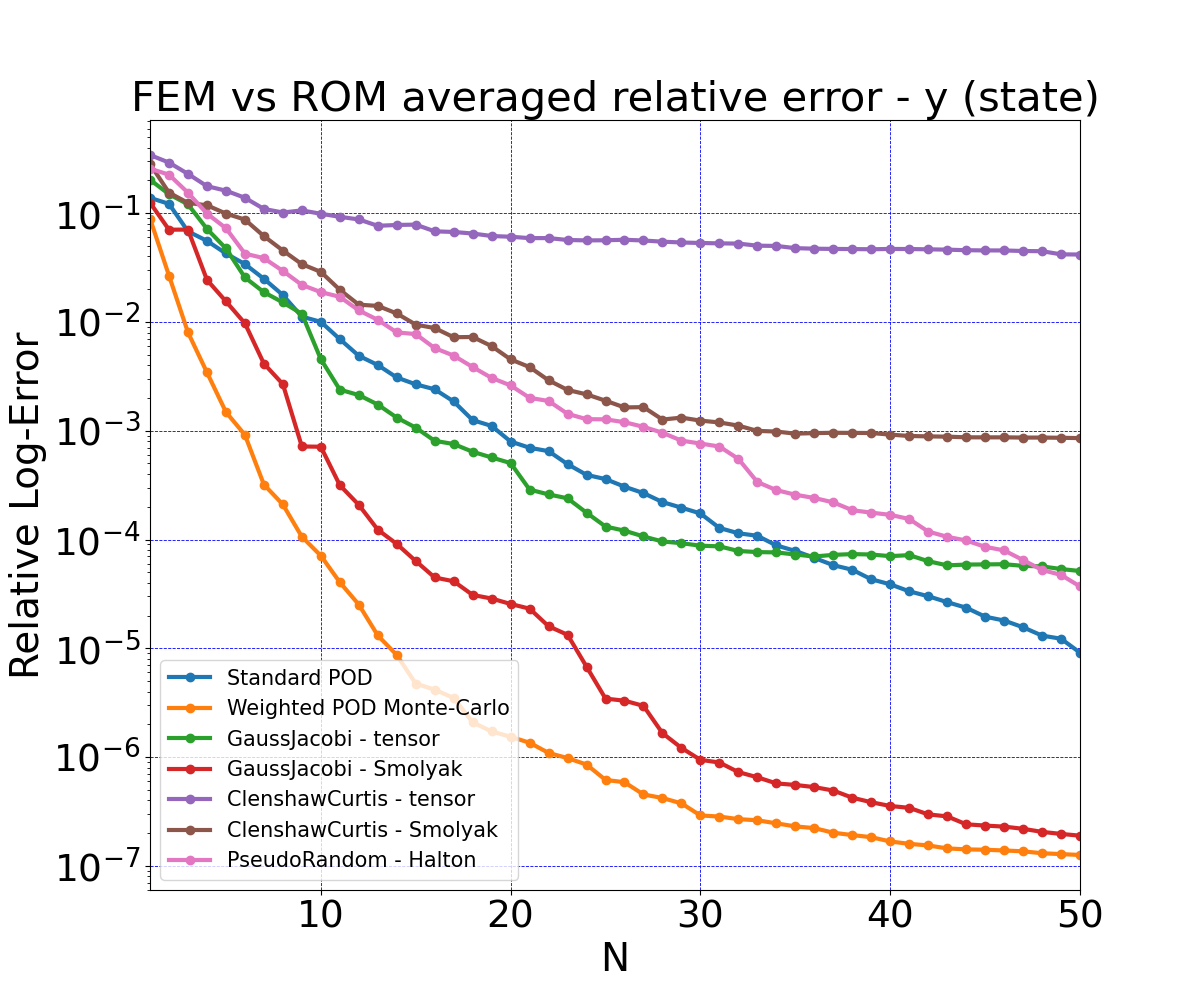}
    \includegraphics[scale=0.15]{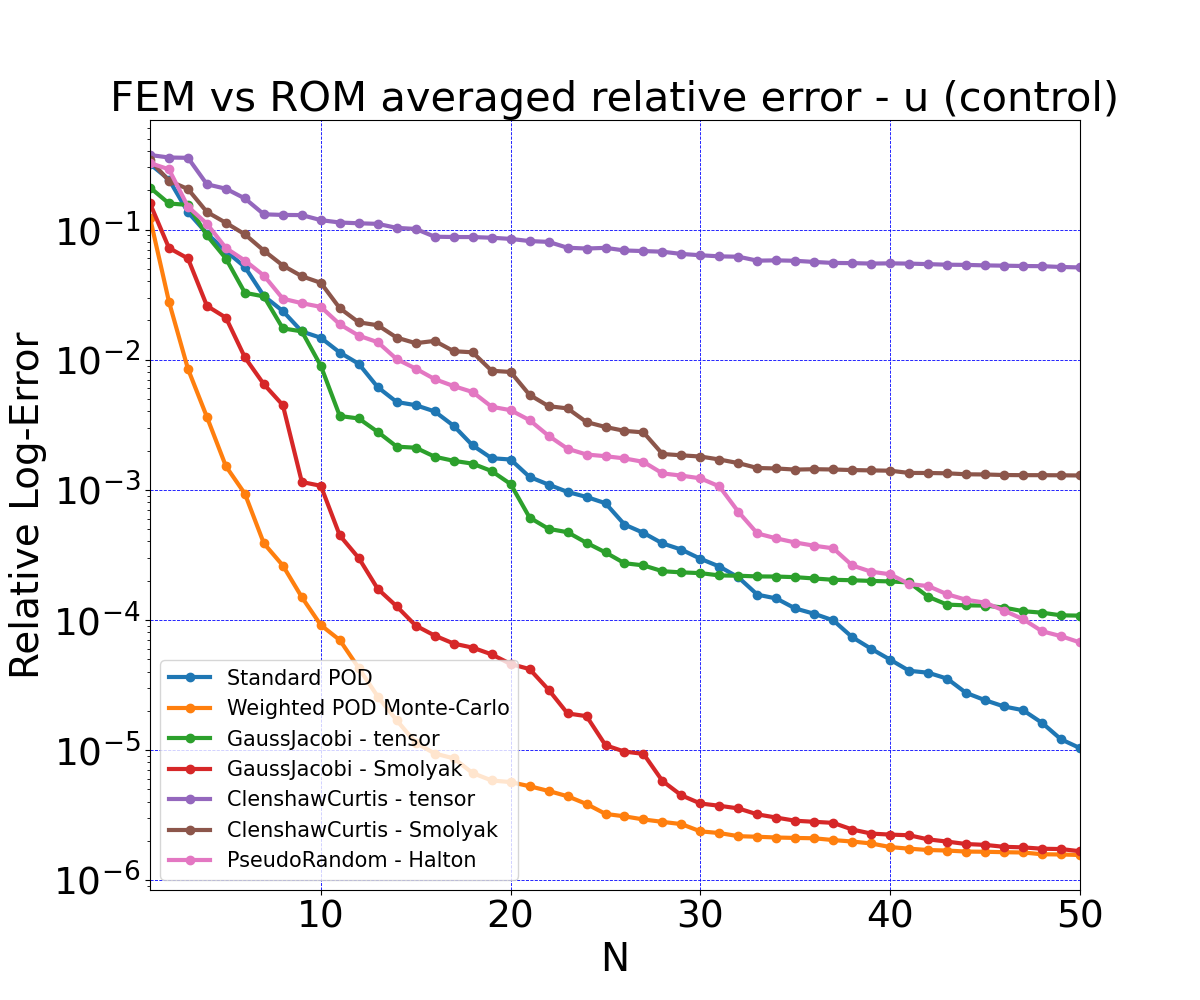}
    \includegraphics[scale=0.15]{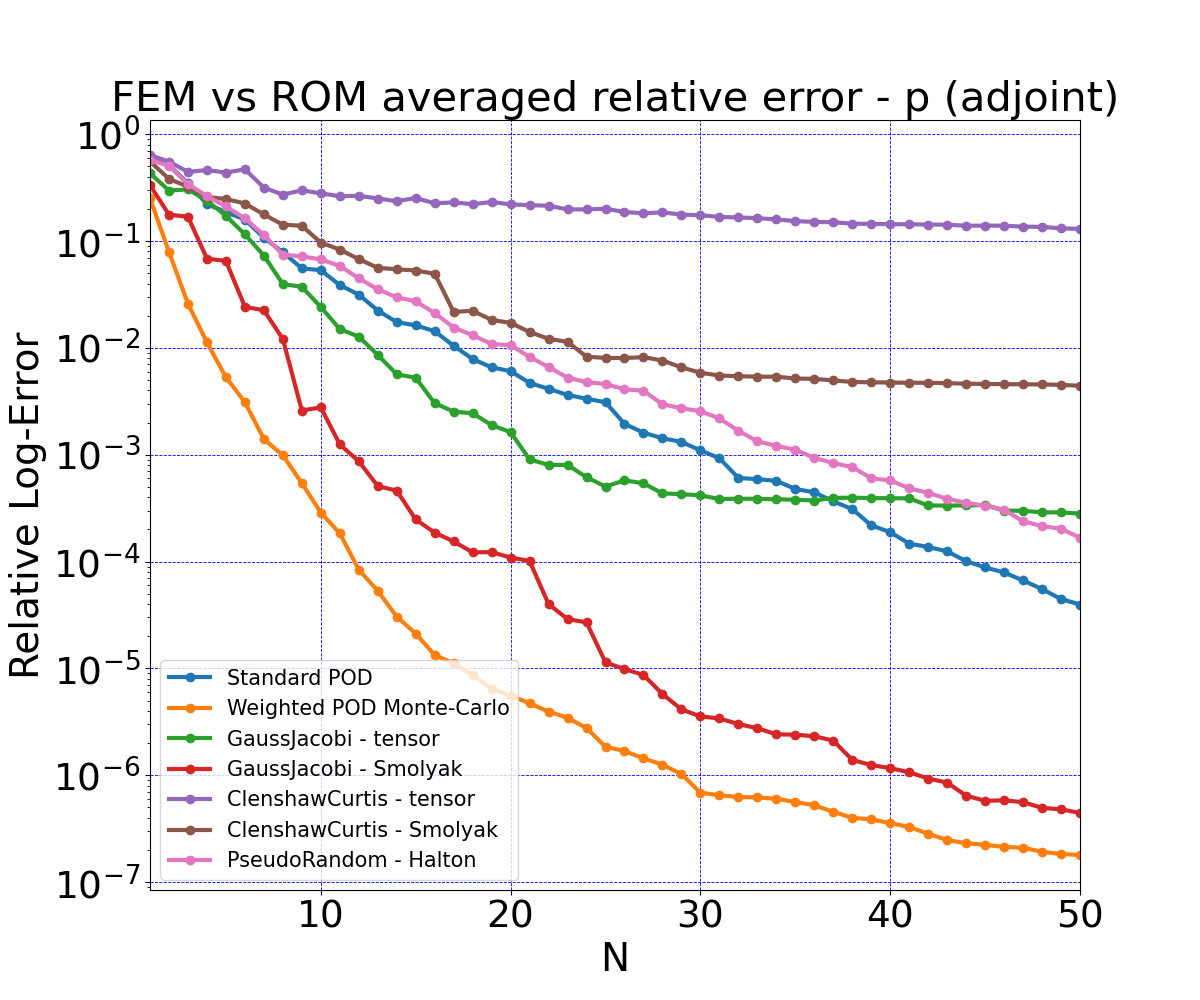}
    \caption{Relative Errors for the Propagating Front in a Square Problem with $\boldsymbol{\mu}$ following distribution \eqref{beta-square-3} on the parameter space $\mathcal{P}$ - \textit{Offline-Online} Stabilization; State (\underline{left}), Control (\underline{center}), Adjoint (\underline{right}); Standard POD (blue), wPOD Monte-Carlo (orange), Gauss-Jacobi tensor rule (green), Gauss-Jacobi Smolyak grid (red), Clenshaw-Curtis tensor rule (cyan), Clenshaw-Curtis Smolyak grid (dark green), Pseudo-Random based on Halton numbers (pink).}
    \label{fig:plot_square_onoffstab-3}
\end{figure}

For the sake of completeness, we simulate other numerical experiments with two more distributions with the same parameter space $\mathcal{P}$, as done in Section \ref{sec:c graetz}. We will only illustrate the Offline-Online errors. In the first experiment, the parameter $\mu$ follows the subsequent distribution
\begin{equation}\label{beta-square-2}
\begin{aligned}
    {\mu_1} \sim 1 + \big(4 \cdot 10^4 - 1\big) X_1, \text{ where } X_1 \sim \text{Beta}(20,20), \\
     \mu_2 \sim 0.9 + \big( 1.5 - 0.9\big) X_2, \text{ where } X_2 \sim \text{Beta}(20,20),
\end{aligned}
\end{equation}
where also in this case most of the parameters are concentrated in the middle of $\mathcal{P}$. Similar to the Graetz-Poiseuille Problem, tensor grids and the Smolyak sampling based on Clenshaw-Curtis quadrature do not perform excellently, whereas Weighted Monte-Carlo outperforms. We show the accuracy in Figure \ref{fig:plot_square_onoffstab-2}. 
Moreover, we simulate the same problem with the following distribution of the parameter on $\mathcal{P}$, which restricts even more the relevant information in the middle of $\mathcal{P}$
\begin{equation}\label{beta-square-3}
\begin{aligned}
    {\mu_1} \sim 1 + \big(4 \cdot 10^4 - 1\big) X_1, \text{ where } X_1 \sim \text{Beta}(30,30), \\
     \mu_2 \sim 0.9 + \big( 1.5 - 0.9\big) X_2, \text{ where } X_2 \sim \text{Beta}(30,30),
\end{aligned}
\end{equation} 
Again, the weighted Monte-Carlo methods performs very well, whereas all quadrature rules that have lots of samples near the boundary have poor results (Figure \ref{fig:plot_square_onoffstab-3}).

Numerical tests of the parabolic version of the Propagating Front in a Square Problem are here illustrated. For a fix $T>0$ and a given $\boldsymbol{\mu} \in \mathcal{P}$ we have to find the pair $(y, u) \in \tilde{\mathcal{Y}} \times \mathcal{U}$ which solves
\begin{equation*}
    \min\limits_{(y,u)} \frac{1}{2} \int\limits_{\Omega_{obs}\times (0,T)}(y(\boldsymbol{\mu})- y_d)^2 \; d \Omega  +
\frac{\alpha}{2} \int\limits_{\Omega \times (0,T)}u(\boldsymbol{\mu})^2 \; d \Omega, \ \text{ such that }
\end{equation*}
\vspace{-2mm}
\begin{equation}\label{par-square-problem}
\begin{cases}
\displaystyle \partial_t y(\boldsymbol{\mu})-\frac{1}{\mu_1} \Delta y(\boldsymbol{\mu})+[\cos{\mu_2},\sin{\mu_2}] \cdot \nabla y(\boldsymbol{\mu})=u(\boldsymbol{\mu}), & \text { in } \Omega  \times (0,T), \\
\displaystyle y(\boldsymbol{\mu})=1, & \text { on } \Gamma_{1} \cup \Gamma_{2}  \times (0,T), \\
\displaystyle y(\boldsymbol{\mu})=0, & \text { on } \Gamma_{3} \cup \Gamma_{4} \cup \Gamma_{5}  \times (0,T), \\
\displaystyle y(\boldsymbol{\mu})(0)=y_0(x), & \text { in } \Omega,
\end{cases}
\end{equation}
where
$y_0(x)=0$ \emph{for all} $x \in {\Omega}$ in Figure \ref{fig:geometry-square}. A final time $T=3.0$ is set. Considering the time discretization, we chose a number of time steps equal to $N_t=30$, then we have $\Delta t=0.1$. Instead, for the spatial approximation, the mesh size is set to $h=0.036$, which implies an overall dimension of the space-time setting equal to $N_{tot}=174780$. For a fixed instant $t$, a single FEM space is characterized by $\mathcal{N}=1942$. For the SUPG procedure, we impose $\delta_K =1.0$ \emph{for all} $K \in \mathcal{T}_{h}$. Setting a penalization parameter $\alpha=0.01$, we try to achieve in a $L^2$-mean a desired solution profile $y_d(x,t)=0.5$, defined \emph{for all} $t \in (0,3)$ and $x$ in $\Omega_{obs}$ of Figure \ref{fig:geometry-square}. 

\begin{figure}
        \centering
        \includegraphics[scale=0.16]{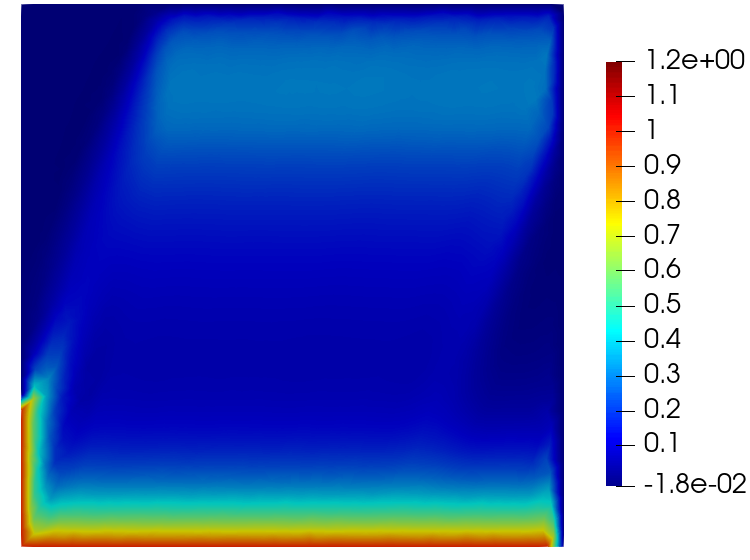}
        \includegraphics[scale=0.16]{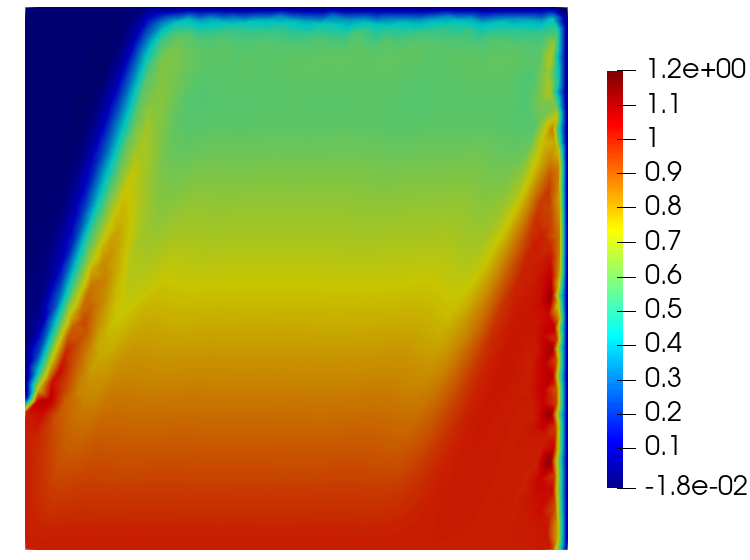}
        \includegraphics[scale=0.16]{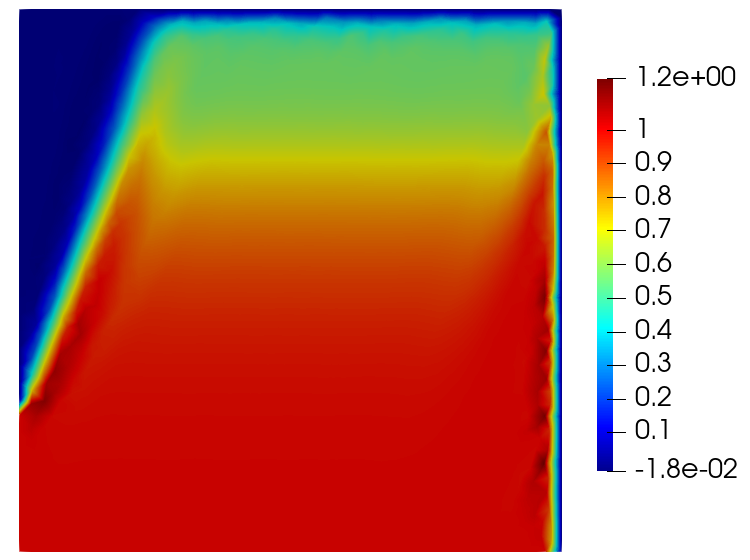}
        \caption{wPOD Monte-Carlo \textit{Offline-Online} stabilized reduced solution of $y$, for $t=0.1$, $t=1.5$, $t=3.0$, $\boldsymbol{\mu} = (2 \cdot 10^4, 1.2)$, $h=0.036$, $\alpha=0.01$, $N_{\text{train}}=100$, $\delta_K=1.0$, $N=30$.}
        \label{fig:par-onff-square-y}
\end{figure}

\begin{figure}
        \centering
        \includegraphics[scale=0.16]{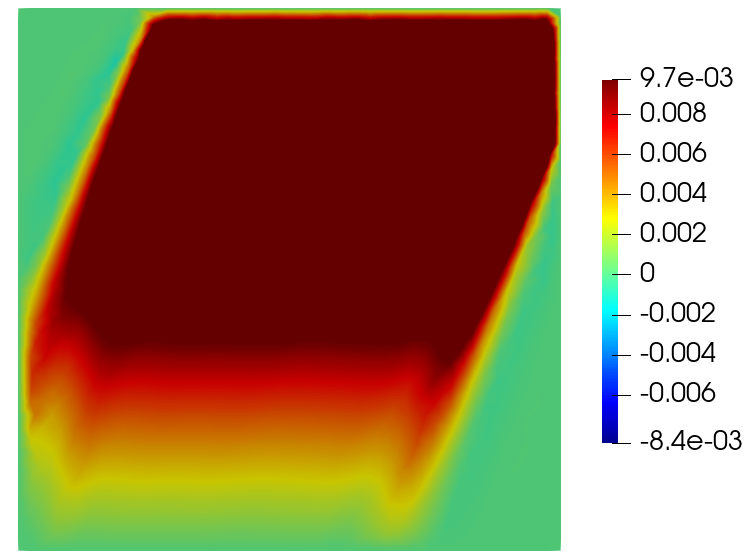}
        \includegraphics[scale=0.16]{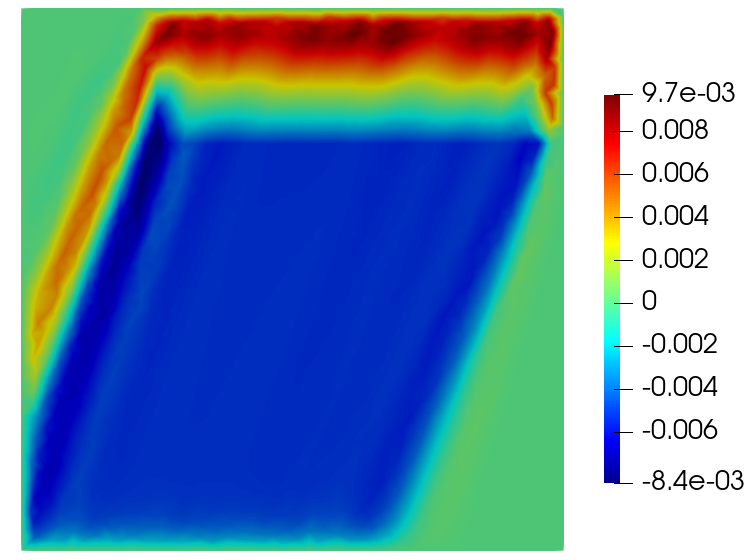}
        \includegraphics[scale=0.16]{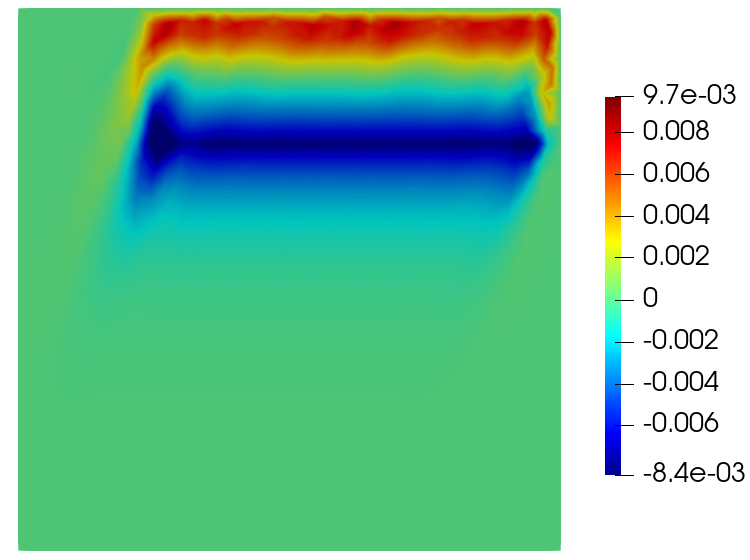}
        \caption{wPOD Monte-Carlo \textit{Offline-Online} stabilized reduced solution of $p$, for $t=0.1$, $t=1.5$, $t=3.0$, $\boldsymbol{\mu} = (2 \cdot 10^4, 1.2)$, $h=0.036$, $\alpha=0.01$, $N_{\text{train}}=100$, $\delta_K=1.0$, $N=30$.}
        \label{fig:par-onff-square-p}
\end{figure}

$\mathcal{P} := \big[1,4 \cdot 10^4\big] \times \big[0.9,1.5\big]$, as in the steady version. We suppose that $\boldsymbol{\mu}$ follows the probability distribution \eqref{beta-square}. Our training set has cardinality $N_{\text{train}}=100$, with the exception of Gauss-Jacobi and Clenshaw-Curtis Smolyak grids with $N_{\text{train}}=89$ and $N_{\text{train}}=93$, respectively, which are the number of nodes nearest to $100$ for this kind of procedure. In Figure \ref{fig:par-onff-square-y} and \ref{fig:par-onff-square-p}, we show a representative stabilized FEM solution for $\boldsymbol{\mu}=(2\cdot 10^{4},1.2)$ for some instants of time of the state $y$ and the adjoint $p$, respectively. We choose to perform all wPOD procedures with $N_{max}=30$.
\sloppy
\sloppy

\begin{figure}
        \centering
        \includegraphics[scale=0.16]{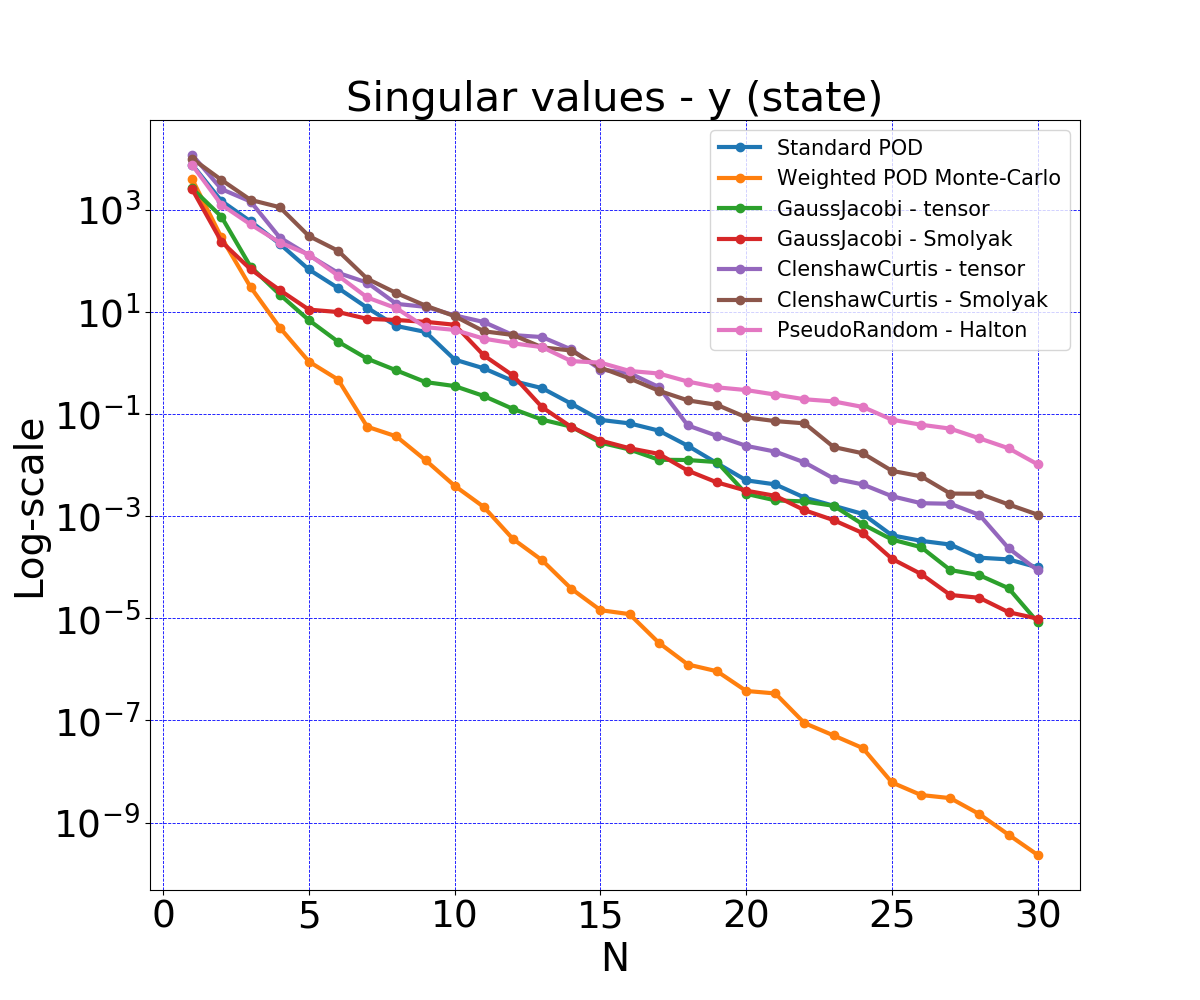} 
        \includegraphics[scale=0.16]{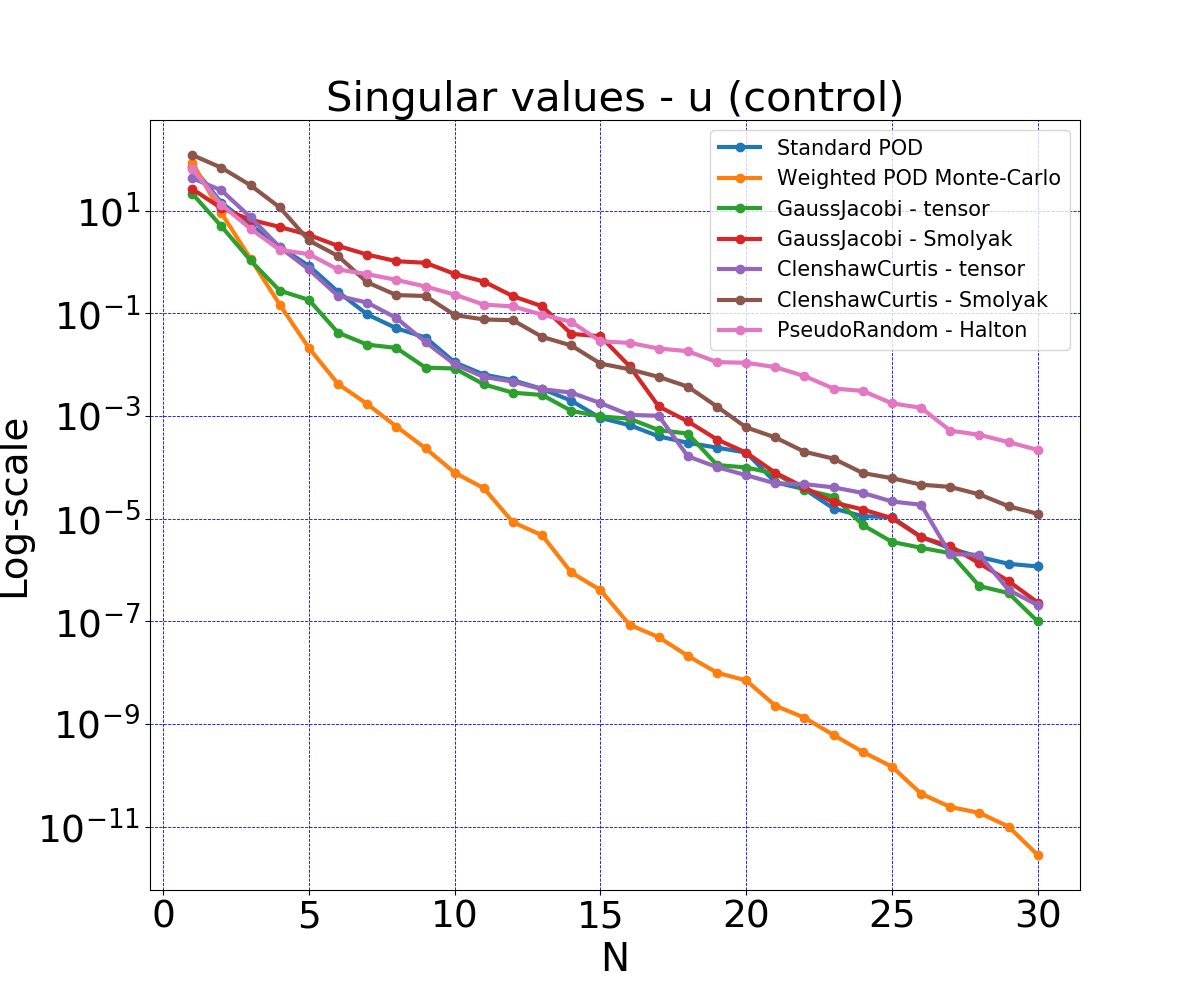} 
        \includegraphics[scale=0.16]{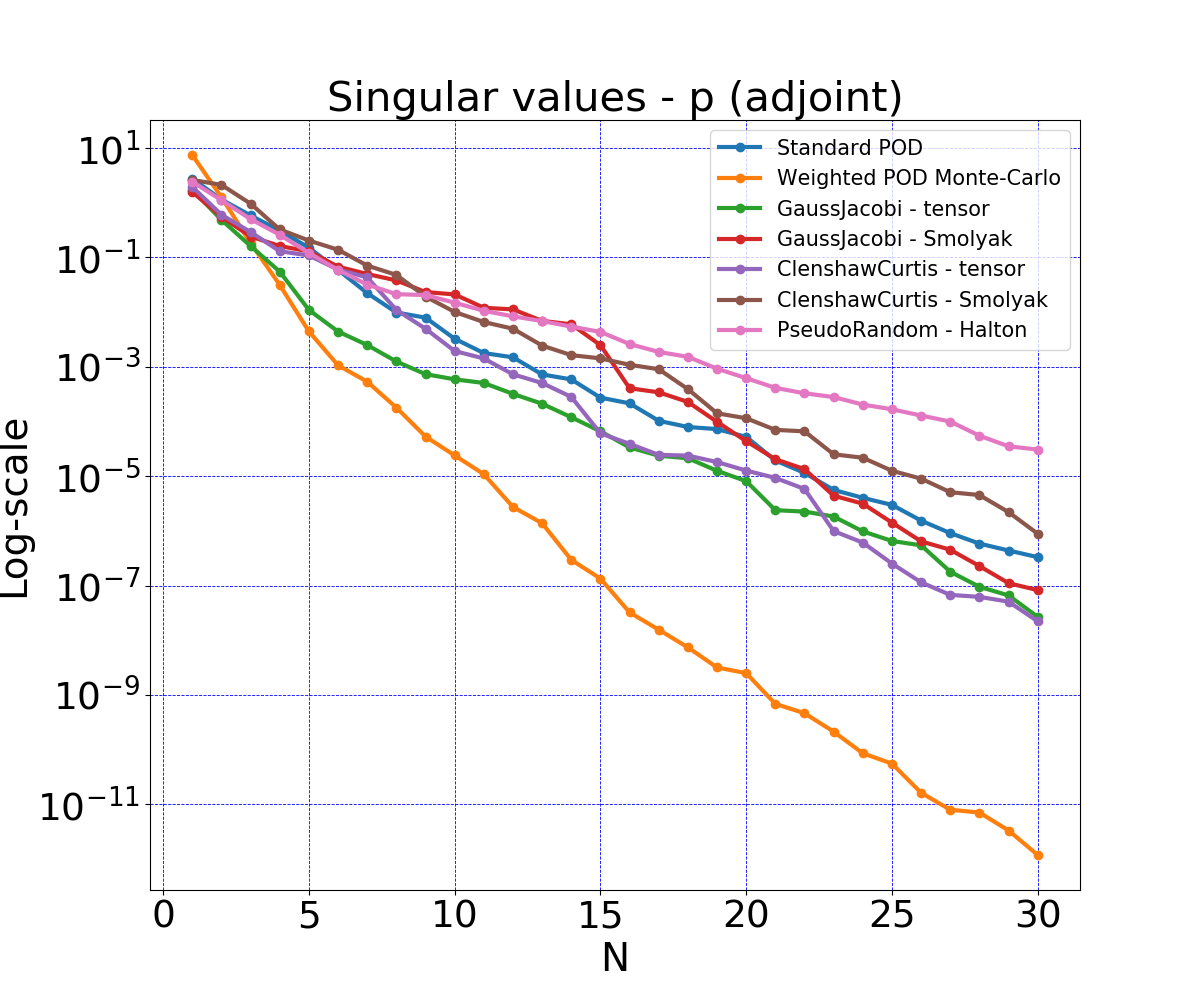} 
        \caption{Singular values decay for the snapshot matrices for the Parabolic Propagating Front in a Problem; State (\underline{left}), Control (\underline{center}), Adjoint (\underline{right}); Standard POD (blue), wPOD Monte-Carlo (orange), Gauss-Jacobi tensor rule (green), Gauss-Jacobi Smolyak grid (red), Clenshaw-Curtis tensor rule (cyan), Clenshaw-Curtis Smolyak grid (dark green), Pseudo-Random based on Halton numbers (pink).}
        \label{fig:plot_Par_Square-sing values}
\end{figure}

\sloppy
Let us move to the error analysis. {In Figure \ref{fig:plot_Par_Square-sing values} we plot the singular value decay for the snapshots matrix for the state, the control, and the adjoint. As motivated in the Graetz-Poiseuille problem, the size of the singular value is huge due to the evolutionary nature hidden in the snapshots. However, these values are in a different scale with respect to the Graetz-Poiseuille ones, because here the problem is not characterized by a parametric geometrical dimension, but the parameters $\mu_1$ and $\mu_2$ directly enter the physical properties of the case of study. Also in this case, the weighted Monte-Carlo POD has the better performance by far. Moreover, in order to compare the real error between the reduced solution and the FEM one, we present also the projection errors between the FEM and the POD space for the three components $y,u,$ and $p$. } 
\sloppy
\sloppy
\begin{figure}
        \centering
        \includegraphics[scale=0.16]{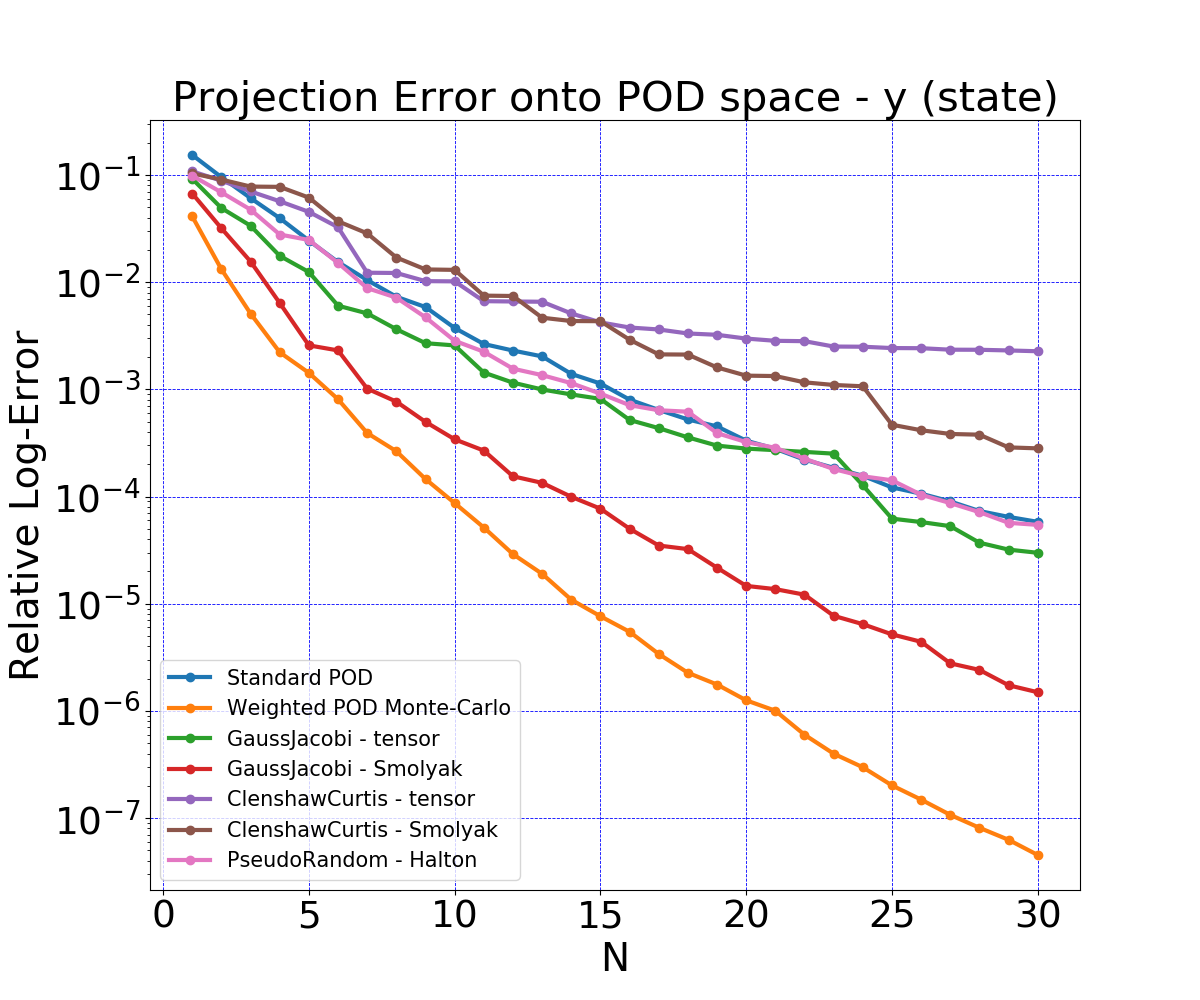} 
        \includegraphics[scale=0.16]{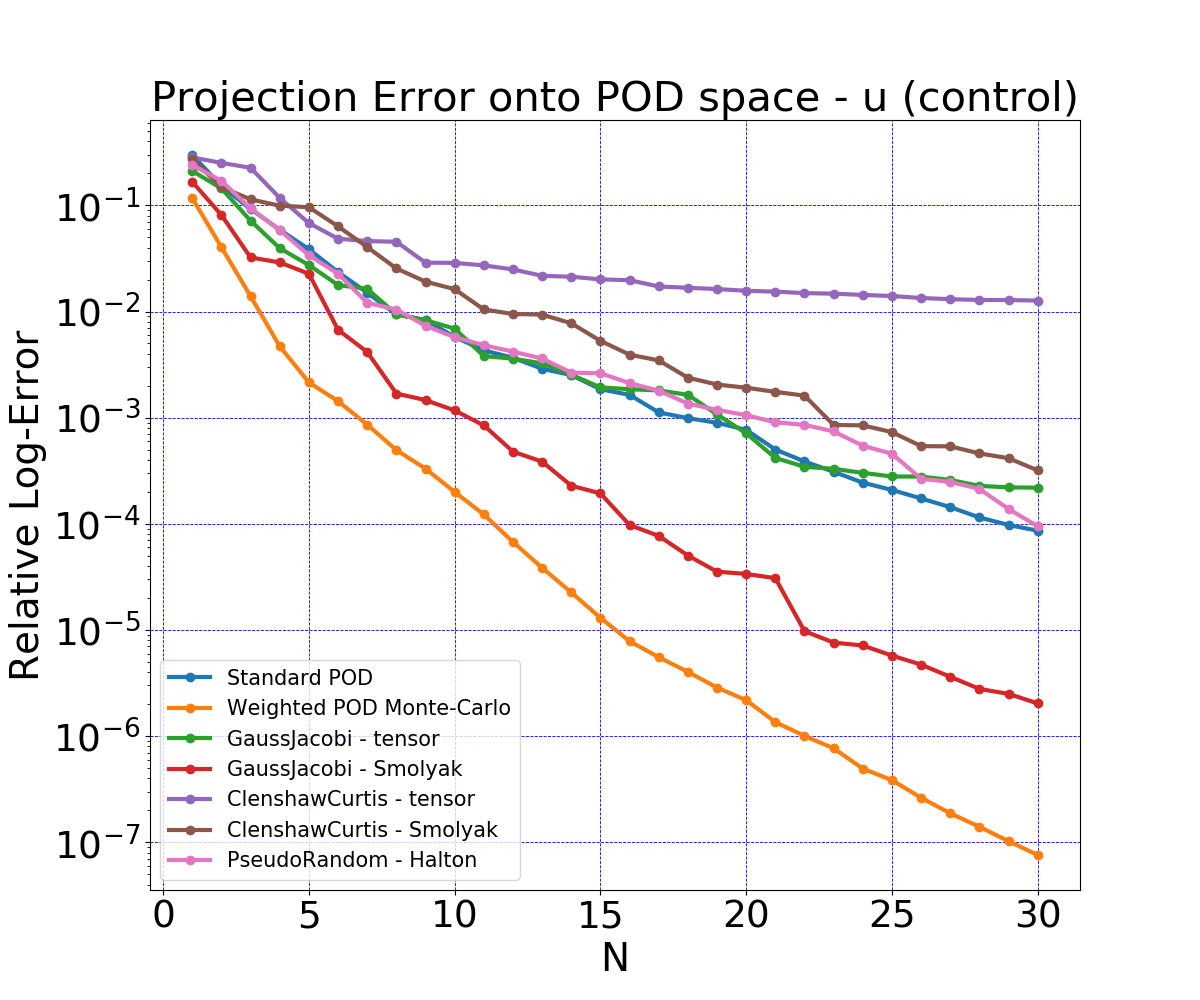} 
        \includegraphics[scale=0.16]{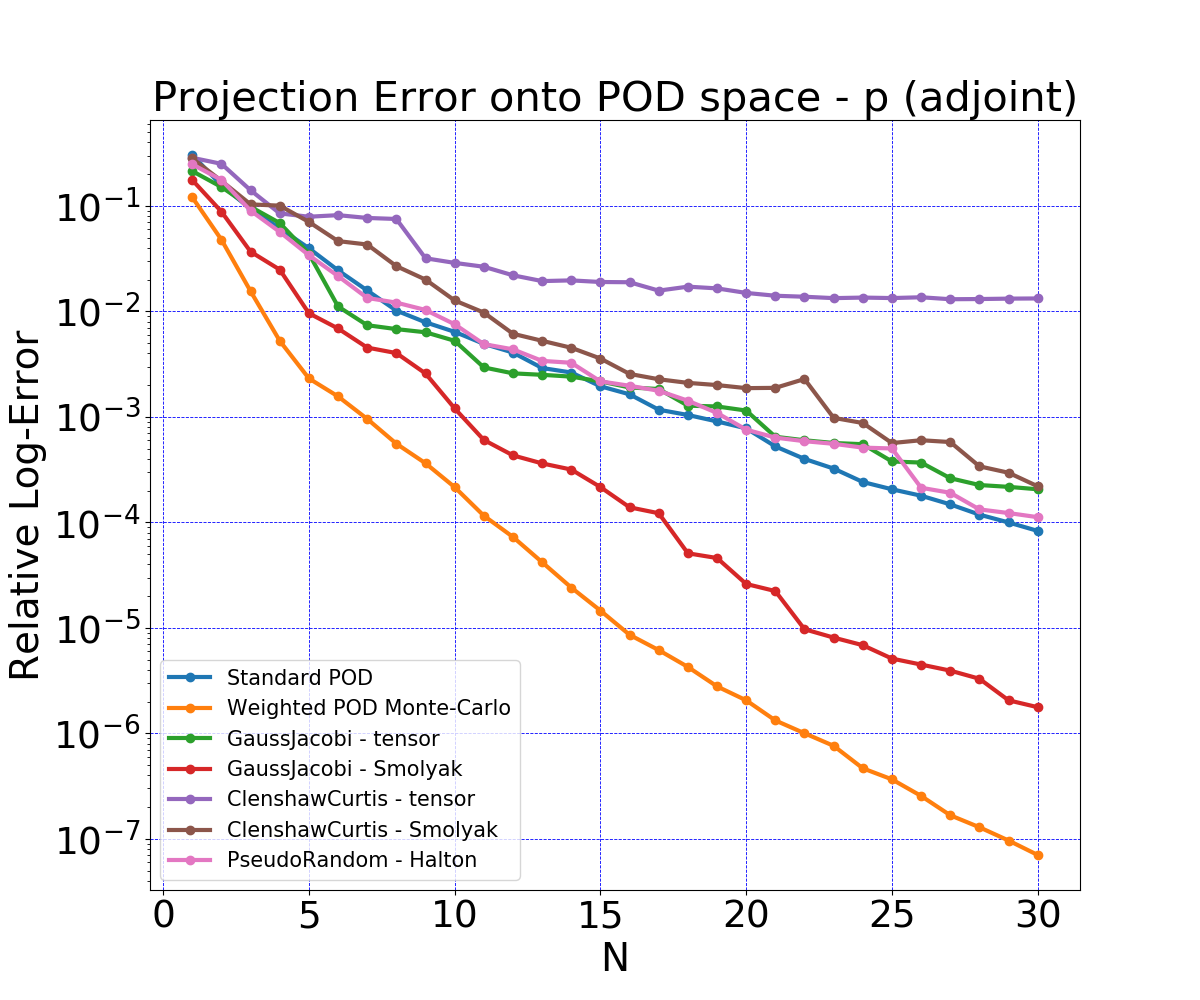} 
        \caption{Projection Errors onto the POD space for the Parabolic Propagating Front in a Problem; State (\underline{left}), Control (\underline{center}), Adjoint (\underline{right}); Standard POD (blue), wPOD Monte-Carlo (orange), Gauss-Jacobi tensor rule (green), Gauss-Jacobi Smolyak grid (red), Clenshaw-Curtis tensor rule (cyan), Clenshaw-Curtis Smolyak grid (dark green), Pseudo-Random based on Halton numbers (pink).}
        \label{fig:plot_Par_Square-proj error}
\end{figure}

\begin{figure}
        \centering
        \includegraphics[scale=0.129]{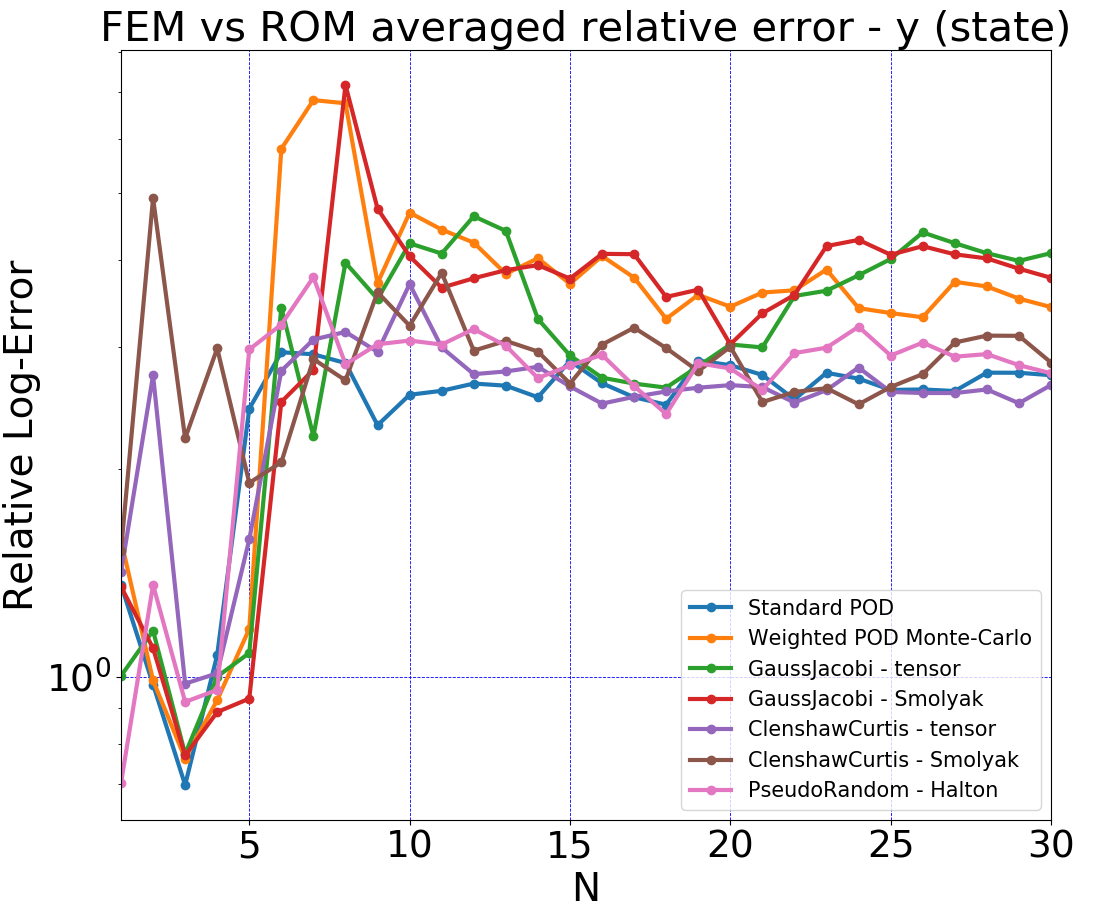}
        \includegraphics[scale=0.129]{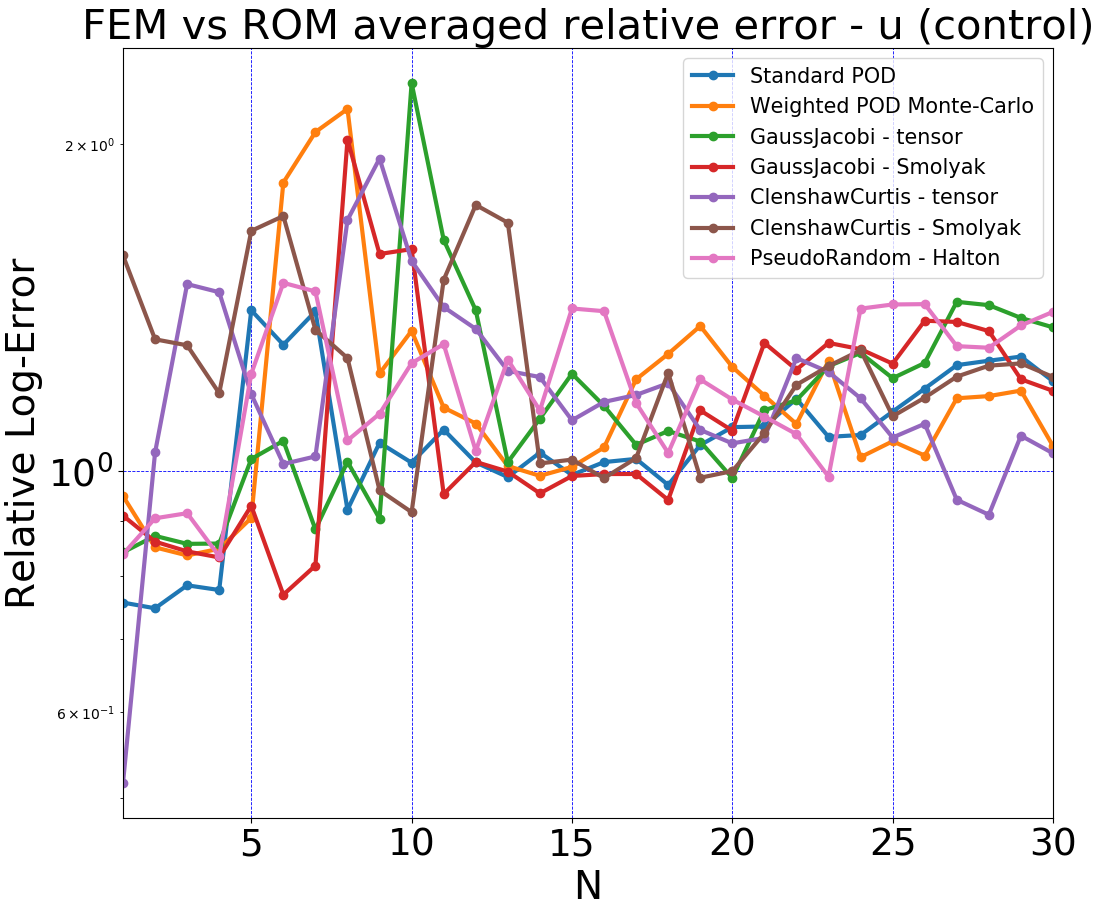}
        \includegraphics[scale=0.129]{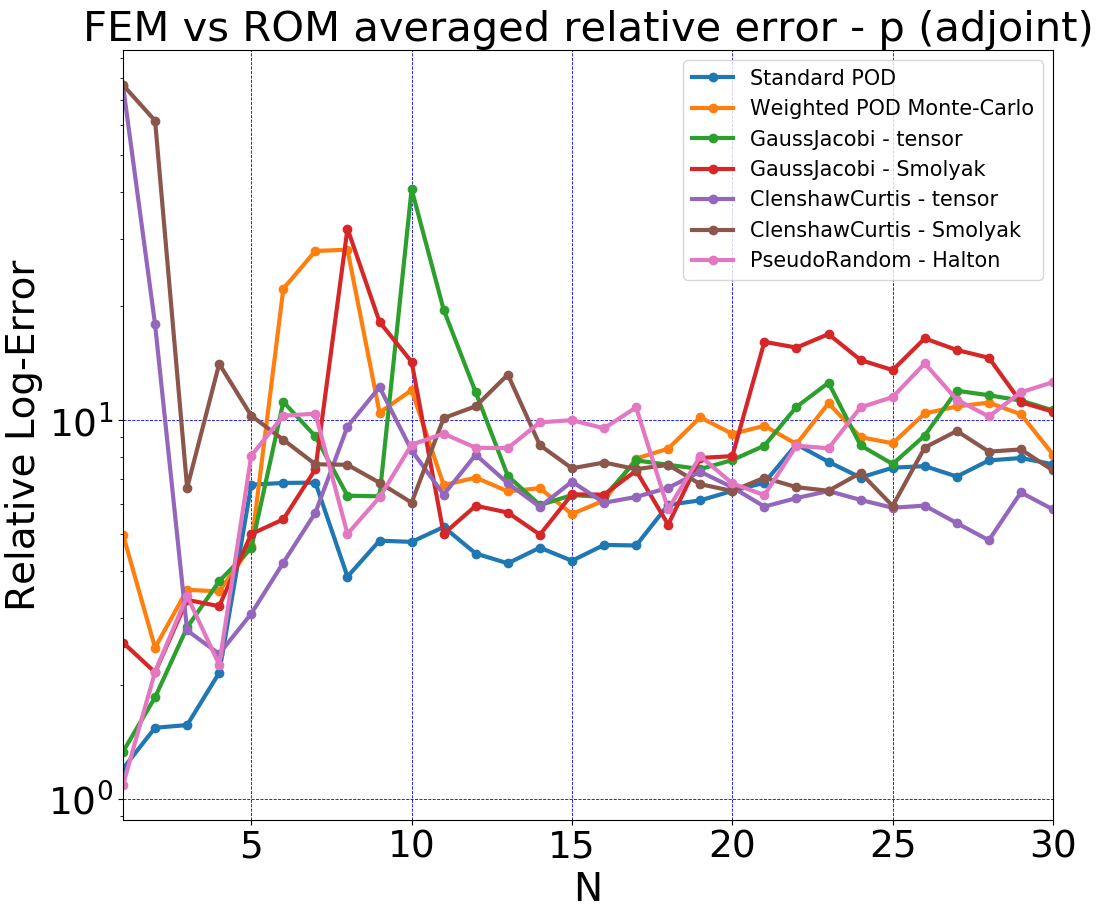}
        \caption{Relative Errors for the Parabolic Propagating Front in a Problem - \textit{Offline-Only} Stabilization; State (\underline{left}), Control (\underline{center}), Adjoint (\underline{right}); Standard POD (blue), wPOD Monte-Carlo (orange), Gauss-Jacobi tensor rule (green), Gauss-Jacobi Smolyak grid (red), Clenshaw-Curtis tensor rule (cyan), Clenshaw-Curtis Smolyak grid (dark green), Pseudo-Random based on Halton numbers (pink).}
        \label{fig:par-off-square-error}
\end{figure}

\begin{figure}
        \centering
        \includegraphics[scale=0.127]{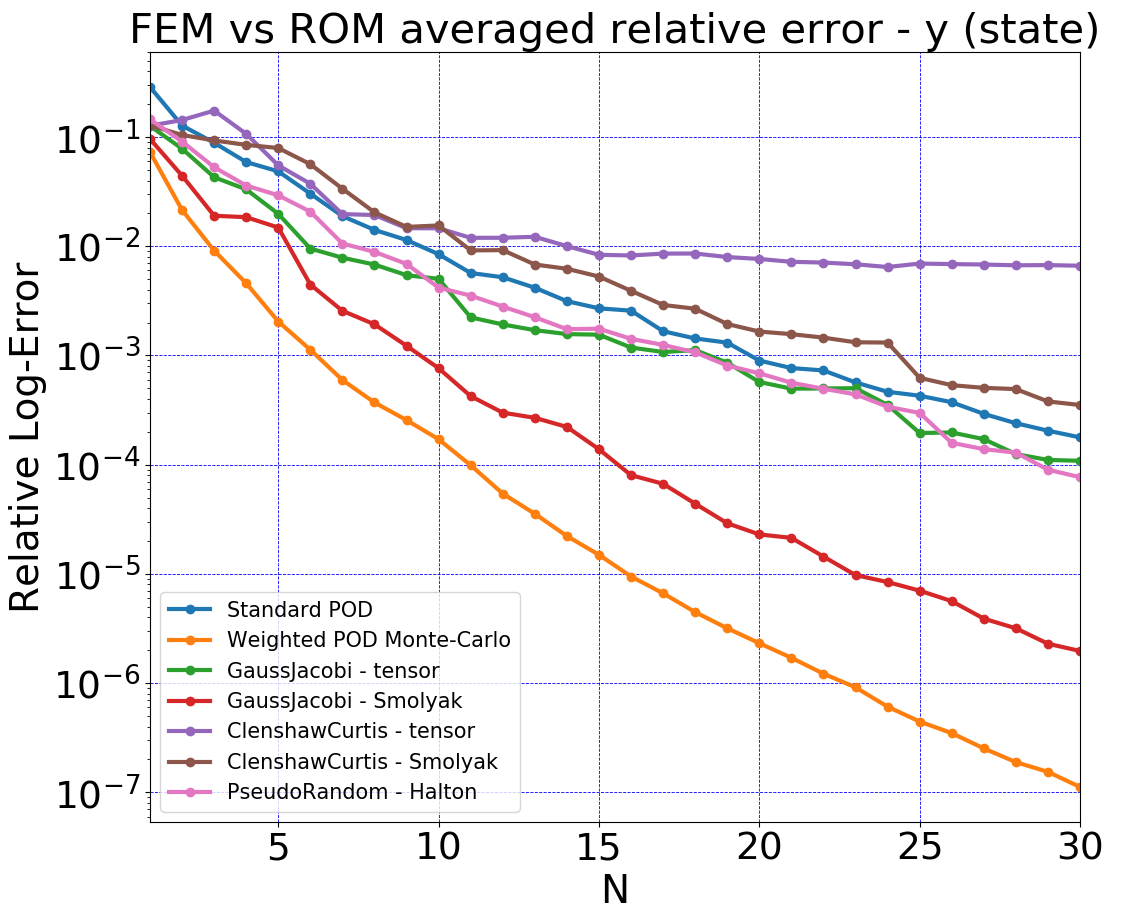}
        \includegraphics[scale=0.127]{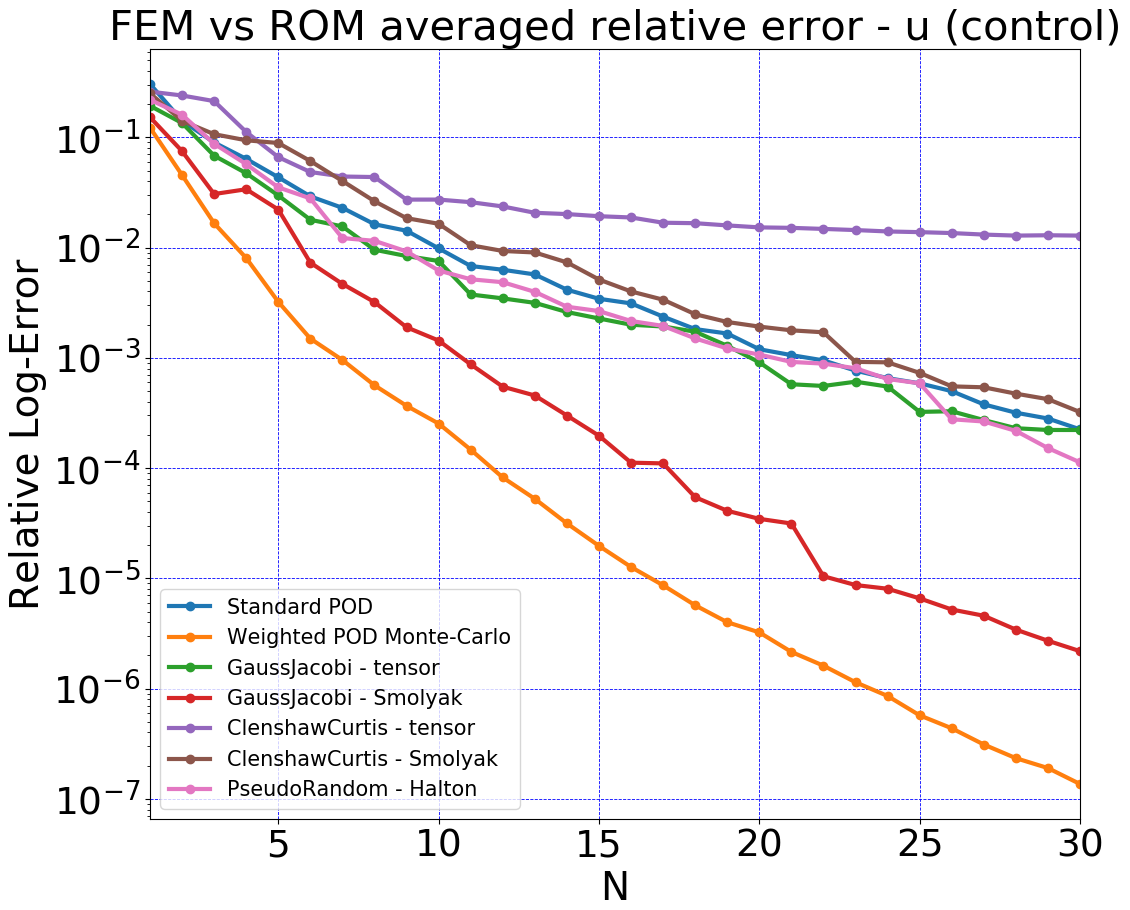}
        \includegraphics[scale=0.127]{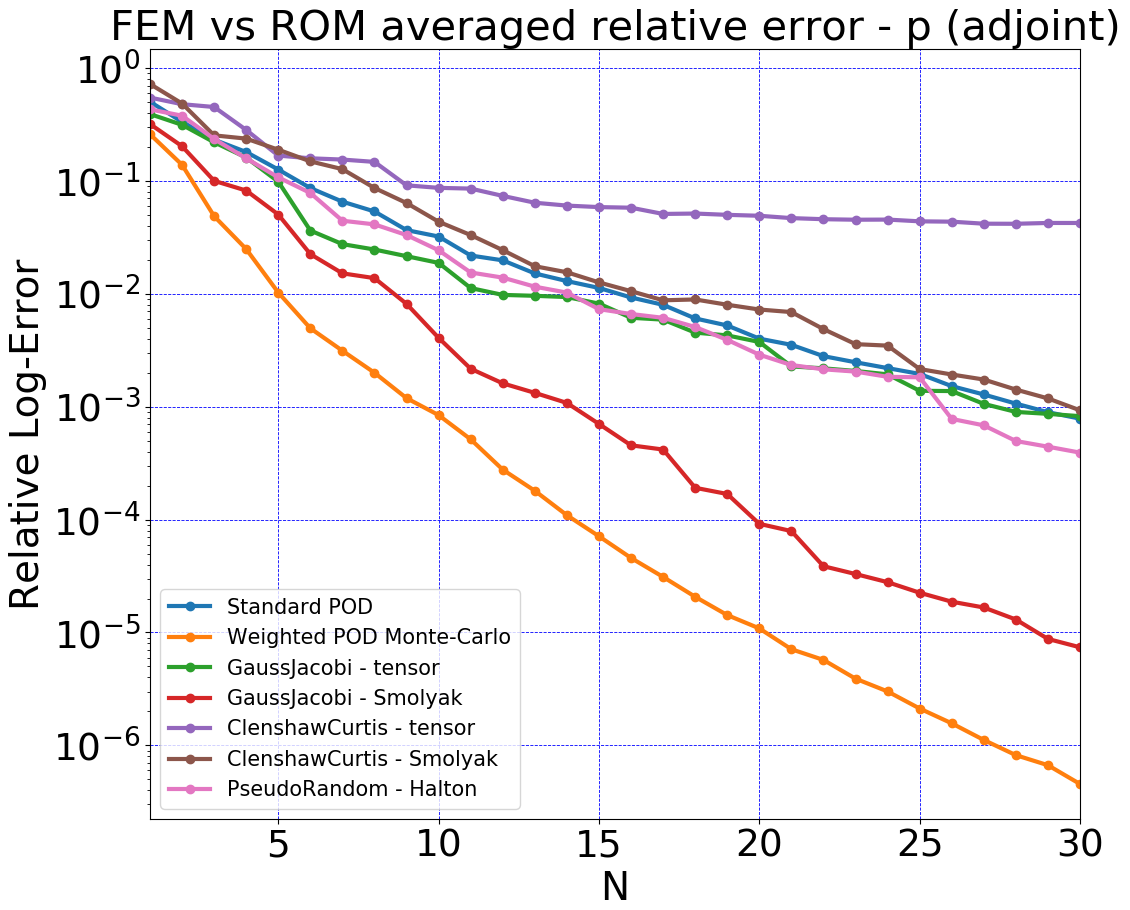}
        \caption{Relative Errors for the Parabolic Propagating Front in a Problem - \textit{Offline-Online} Stabilization; State (\underline{left}), Control (\underline{center}), Adjoint (\underline{right}); Standard POD (blue), wPOD Monte-Carlo (orange), Gauss-Jacobi tensor rule (green), Gauss-Jacobi Smolyak grid (red), Clenshaw-Curtis tensor rule (cyan), Clenshaw-Curtis Smolyak grid (dark green), Pseudo-Random based on Halton numbers (pink).}
        \label{fig:par-onoff-square-error}
\end{figure}

\sloppy
In Figure \ref{fig:par-off-square-error}, we illustrate the relative errors for the \textit{Offline-Only} stabilization. The performance are not satisfactory here, too, where no quantity drops below the accuracy of $10^{-1}$ for all $N$. {Again, all these quantities do not follow any closely the trends of the projection errors in Figure \ref{fig:plot_Par_Square-proj error}.}.
\sloppy
Instead, \textit{Offline-Online} stabilization procedure performs well, as one can notice from Figure \ref{fig:par-onoff-square-error}, {being of the same behavior of the projection ones in Figure \ref{fig:plot_Par_Square-proj error}}. Again, wPOD Monte-Carlo has the best behavior, it reaches $e_{y, 30} = 1.12 \cdot 10^{-7}$ for the state, for the adjoint $e_{p, 30}= 4.55 \cdot 10^{-7}$ and the control $e_{u, 30}=1.36 \cdot 10^{-7}$. Also in this case, isotropic sparse grid techniques are a better choice than tensor rules, both for Gauss-Jacobi and Clenshaw-Curtis approximations.
\sloppy
\sloppy
In Table \ref{speed-par-s}, we compare the speedup-index for all the weighted algorithms: performance are similar for all $N$, for $N=30$ we computed nearly $4000$ \textit{Offline-Online} stabilized reduced solutions in the time of a FEM one.
\sloppy
\begin{table}
\centering

\begin{tabular}{|c|c|c|c|c|c|c|c|}
\hline  \multicolumn{8}{|c|}{Speedup-index Parabolic Propagating front in a Square Problem: \textit{Offline-Online} Stabilization}  \\
\hline $N$ & POD & wPOD & Gauss tensor & Gauss Smolyak & CC tensor & CC Smolyak & Ps. Random \\
\hline $5$ & $6601.9$ & $6503.5$  &  $6702.6$ & $6629.1$ & $6566.6$ & $6605.6$ & $6575.8$\\
\hline $10$ & $6275.9$ &  $6208.0$ & $6336.0$ & $6277.9$ & $6204.4$ & $6293.4$ & $6212.4$\\
\hline $15$ & $5814.3$ & $5702.4$  & $5838.7$ & $5794.3$ & $5699.6$ & $5752.1$ & $5723.7$\\
\hline $20$ & $5327.9$ & $5190.4$  & $5329.8$ & $5270.3$ & $5277.2$ & $5235.9$ & $5197.6$ \\
\hline $25$ & $4465.2$ & $4303.3$  & $4562.6$ & $4422.2$ & $4541.3$ & $4433.1$ & $4479.9$ \\
\hline $30$ & $4061.5$ & $3959.5$  & $4140.3$ & $4026.3$ & $4100.3$ & $4035.6$ & $4043.8$ \\
\hline
\end{tabular}
  \caption{Average Speedup-index of \textit{Offline-Online} Stabilization for the Parabolic Propagating Front in a Square Problem. From left to right: Standard POD, wPOD Monte-Carlo, Gauss-Jacobi tensor, Gauss-Jacobi Smolyak grid, Clenshaw-Curtis tensor, Clenshaw-Curtis Smolyak grid, Pseudo-Random based on Halton numbers. $\mu_1$, $\mu_2$ $\sim $ Beta(10,10)}
  \label{speed-par-s}
\end{table}

\section{Conclusions and Perspectives}
\label{sec:conclusions}
In this work, we illustrated some numerical tests concerning stabilized Parametrized Advection-Dominated OCPs with random parametric inputs in a ROM context. We deal with both steady and unsteady cases and we took advantage of the SUPG stabilization to overcome numerical issues due to high values of the P\'eclet number. Two possibilities of stabilization were analyzed: when SUPG is only used occurs in the offline phase, \textit{Offline-Only stabilization}, or when it is provided in both online and offline phases, \textit{Offline-Online stabilization}. 

In order to deal with the uncertainty quantification caused by random inputs, we consider {wROM}. More precisely, we built our reduced bases using a wPOD in a partitioned approach, using different quadrature rules. We implemented wPOD Monte-Carlo, Gaussian quadrature formulae based on Jacobi polynomials in a tensor rule, approximation related to Clenshaw-Curtis tensor rule, Smolyak isotropic sparse grid approximation of the last two methods, quasi Monte-Carlo method as a Pseudo-Random rule defined on Halton numbers.

We analyzed relative errors between the reduced and the high fidelity solutions and the speedup-index concerning the Graetz-Poiseuille and Propagating Front in a Square Problems, always under a distributed control. 
For the state, control, and adjoint spaces we implemented a $\mathbb{P}^{1}$-FEM approximation in a \textit{optimize-then-discretize} framework as the truth solution. Concerning parabolic problems, a space-time approach is followed applying SUPG in a suitable way. In order to establish which wPOD performs better, we compare them through the same testing set sampled by a Monte-Carlo method according to the probability distribution of the parameter.

\textit{Offline-Only} stabilization technique performed very poorly in terms of errors, this happened for all wROMs considered.
Instead, in all the steady and unsteady experiments, the wROM technique performed excellently in an \textit{Offline-Online} stabilization framework, {where all the relative errors stay very closed to the corresponding projection ones, which represent the best possible reduced approximation that one can achieve with respect to a precise reduced order space}. For parabolic problems, the speedup-index features large values thanks to the space-time formulation. More precisely, {wPOD Monte-Carlo} technique was always the best performer for relative errors, instead, concerning computational efficiency all methods seem equivalent. In addition, the efficiency of the wPOD Monte-Carlo is supported by the fact that after a small number of reduced basis it is nearly $100$ times more accurate than a Standard POD in a deterministic context. Moreover, we notice that {sparse grids perform better than relative tensor ones, although having a bit less number of quadrature nodes.}

Furthermore, in the Graetz-Poiseuille Problem we illustrate that under geometrical parametrization affected by randomness, wROMs still have good performance, despite small fluctuations in the graph of relative errors.

As a first perspective, it might be interesting to create a strongly-consistent stabilization technique that flattens all the fluctuations of geometrical parametrization in a UQ context. Moreover, we want to extend the study to boundary control and convection dominated non-linear PDEs constraints. It might also be interesting to study the performance of other stabilization techniques for the online phases, for instance, of the Online Vanishing Viscosity and the Online Rectification methods \cite{ali2018stabilized, chakir2019non, maday1989analysis} combined with the SUPG technique in the offline phase or with the stabilization strategy used in \cite{torlo2018stabilized}. {In this case, the main advantage to use Vanishing Viscosity or Rectification methods instead of SUPG in the online phases lies in the fact that the former techniques are less expensive than the latter \cite{maday2016online,ali2021reduced}. In case of similar accuracy results, one would gain in the speedup-index. Moreover, these approaches are very versatile as they are independent of the stabilization procedure performed during the offline phase.}

{Finally, a further study of the ``discretize-then-optimize" approach might be beneficial to compare its error results with the ones of this manuscript in order to choose the best numerical strategy.}

\section*{Acknowledgements}
We acknowledge the support by European Union Funding for Research and Innovation -- Horizon 2020 Program -- in the framework of European Research Council Executive Agency: Consolidator Grant H2020 ERC CoG 2015 AROMA-CFD project 681447 ``Advanced Reduced Order Methods with Applications in Computational Fluid Dynamics''. We also acknowledge the PRIN 2017  ``Numerical Analysis for Full and Reduced Order Methods for the efficient and accurate solution of complex systems governed by Partial Differential Equations'' (NA-FROM-PDEs) and the INDAM-GNCS project ``Tecniche Numeriche Avanzate per Applicazioni Industriali''.
The computations in this work have been performed with RBniCS \cite{RBniCS} library, developed at SISSA mathLab, which is an implementation in FEniCS \cite{logg2012automated} of several reduced order modelling techniques; we acknowledge developers and contributors to both libraries.

\bibliographystyle{plain}
\bibliography{BIB}

\end{document}